\pdfoutput=1
\RequirePackage{ifpdf}
\ifpdf 
\documentclass[pdftex]{sigma}
\else
\documentclass{sigma}
\fi

\newcommand{\id}{{\rm id}}

\newcommand{\ext}{{\rm ext}}
\renewcommand{\exp}{{\rm exp}}
\newcommand{\End}{{\rm{End}\ts}}

\newcommand{\diag}{{\rm diag}}

\newcommand{\wh}{\widehat}
\newcommand{\ot}{\otimes}

\newcommand{\al}{\alpha}
\newcommand{\be}{\beta}

\newcommand{\ep}{\epsilon}

\newcommand{\vp}{\varphi}
\newcommand{\de}{\delta}
\newcommand{\ze}{\zeta}

\newcommand{\hra}{\hookrightarrow}

\newcommand{\ts}{\,}
\newcommand{\pr}{^{\tss\prime}}

\newcommand{\tss}{\hspace{1pt}}

\newcommand{\CC}{\mathbb{C}\tss}

\newcommand{\ZZ}{\mathbb{Z}\tss}

\newcommand{\Lc}{\mathcal{L}}

\newcommand{\Rc}{\mathcal{R}}

\newcommand{\Cc}{\mathcal{C}}
\newcommand{\Ec}{\mathcal{E}}
\newcommand{\Fc}{\mathcal{F}}

\newcommand{\Hc}{\mathcal{H}}

\newcommand{\Xc}{\mathcal{X}}

\newcommand{\gl}{\mathfrak{gl}}

\newcommand{\oa}{\mathfrak{o}}

\newcommand{\e}{\mathfrak{e}}
\newcommand{\f}{\mathfrak{f}}
\newcommand{\g}{\mathfrak{g}}
\newcommand{\h}{\mathfrak h}

\newcommand{\z}{\mathfrak{z}}

\newcommand{\tra}{ {\rm t}}

\newcommand{\bi}{\bar{\imath}}
\newcommand{\bj}{\bar{\jmath}}

\newcommand{\Sym}{\mathfrak S}

\numberwithin{equation}{section}

\newtheorem{thm}{Theorem}[section]
\newtheorem{lem}[thm]{Lemma}
\newtheorem{prop}[thm]{Proposition}
\newtheorem{cor}[thm]{Corollary}
\newtheorem{conj}[thm]{Conjecture}
\newtheorem*{mthm}{Main Theorem}

\theoremstyle{definition}
\newtheorem{defin}[thm]{Definition}
\newtheorem{Remark}[thm]{Remark}

\newcommand{\bth}{\begin{thm}}
\renewcommand{\eth}{\end{thm}}
\newcommand{\bpr}{\begin{prop}}
\newcommand{\epr}{\end{prop}}
\newcommand{\ble}{\begin{lem}}
\newcommand{\ele}{\end{lem}}
\newcommand{\bco}{\begin{cor}}
\newcommand{\eco}{\end{cor}}
\newcommand{\bde}{\begin{defin}}
\newcommand{\ede}{\end{defin}}
\newcommand{\bex}{\begin{example}}
\newcommand{\eex}{\end{example}}
\newcommand{\bcj}{\begin{conj}}
\newcommand{\ecj}{\end{conj}}

\begin{document}

\allowdisplaybreaks

\newcommand{\arXivNumber}{1911.03496}

\renewcommand{\PaperNumber}{043}

\FirstPageHeading

\ShortArticleName{Isomorphism between the $R$-Matrix and Drinfeld Presentations}

\ArticleName{Isomorphism between the $\boldsymbol{R}$-Matrix and Drinfeld\\ Presentations of Quantum Affine Algebra:\\ Types $\boldsymbol{B}$ and $\boldsymbol{D}$}

\Author{Naihuan JING~$^{\dag}$, Ming LIU~$^{\ddag\S}$ and Alexander MOLEV~$^\S$}

\AuthorNameForHeading{N.~Jing, M.~Liu and A.~Molev}

\Address{$^\dag$~Department of Mathematics, North Carolina State University, Raleigh, NC 27695, USA}
\EmailD{\href{mailto:jing@math.ncsu.edu}{jing@math.ncsu.edu}}

\Address{$^\ddag$~School of Mathematics, South China University of Technology, Guangzhou, 510640, China}
\EmailD{\href{mailto:mamliu@scut.edu.cn}{mamliu@scut.edu.cn}}

\Address{$^\S$~School of Mathematics and Statistics, University of Sydney, NSW 2006, Australia}
\EmailD{\href{mailto:alexander.molev@sydney.edu.au}{alexander.molev@sydney.edu.au}}

\ArticleDates{Received November 18, 2019, in final form May 10, 2020; Published online May 21, 2020}

\Abstract{Following the approach of Ding and Frenkel~[\textit{Comm. Math. Phys.} \textbf{156} (1993), 277--300] for type $A$, we showed in our previous work [\textit{J.~Math. Phys.} \textbf{61} (2020), 031701, 41~pages] that the Gauss decomposition of the generator matrix in the $R$-matrix presentation of the quantum affine algebra yields the Drinfeld generators in all classical types. Complete details for type~$C$ were given therein, while the present paper deals with types~$B$ and~$D$. The arguments for all classical types are quite similar so we mostly concentrate on necessary additional details specific to the underlying orthogonal Lie algebras.}


\Classification{17B37; 17B69}

\section{Introduction}\label{sec:int}

The quantum affine algebras $U_q(\wh\g)$ associated with simple Lie algebras~$\g$ admit at least three different presentations. The original definition of Drinfeld~\cite{d:ha} and Jimbo~\cite{j:qd} was followed by the {\em new realization} of Drinfeld~\cite{d:nr} which is also known as the {\em Drinfeld presentation}, while the $R$-{\em matrix presentation} was introduced by Reshetikhin and Semenov-Tian-Shansky~\cite{rs:ce} and further developed by Frenkel and Reshetikhin~\cite{fri:qa}. A detailed construction of an isomorphism between the first two presentations was given by Beck~\cite{b:bg}.

An isomorphism between the Drinfeld and $R$-matrix presentations of the algebras $U_q(\wh\g)$
in type $A$ was constructed by Ding and Frenkel~\cite{df:it}. In our previous work~\cite{jlm:ib-c} we were able to extend this construction to the remaining classical types and gave detailed arguments in type~$C$.
The present article is concerned with types $B$ and $D$, where we use the same approach as in~\cite{jlm:ib-c} and mostly concentrate on necessary changes specific to the orthogonal Lie algebras $\oa_N$ and their root systems. In particular, this applies to low rank relations with the underlying Lie algebras~$\oa_3$ and~$\oa_4$, and to formulas for the universal $R$-matrices.

As with the corresponding isomorphisms between the $R$-matrix and Drinfeld presentations
of the Yangians (see respective details in \cite{bk:pp, grw:eb} and \cite{jlm:ib}),
their counterparts in the quantum affine algebra case allow one to
connect two sides of the representation theory in an explicit way:
the parameterization of finite-dimensional irreducible representations via their
Drinfeld polynomials can be translated from one presentation to another;
see \cite[Chapter~12]{cp:gq}, \cite{gm:rt} and \cite{t:im}.
As another consequence of the
isomorphism theorems, one can derive the Poincar\'e--Birkhoff--Witt
theorem for the $R$-matrix presentation of the
quantum affine algebra from the corresponding result
of Beck~\cite{b:cb} for $U_q(\wh\g)$. We will give a more detailed account
of these applications in our forthcoming project.

To work with the quantum affine algebras in types $B$ and $D$, we apply the {\em Gauss decomposition} of the generator matrices in the $R$-matrix presentation in the same way as in types $A$ and $C$; see~\cite{df:it} and~\cite{jlm:ib-c}. We show that the generators arising from the Gauss decomposition satisfy the required relations of the Drinfeld presentation. To demonstrate that the resulting homomorphism is injective we follow the argument of Frenkel and Mukhin~\cite{fm:ha} and rely on the formula for the universal
$R$-matrix due to Khoroshkin and Tolstoy~\cite{kt:ur} and Damiani~\cite{d:rm}.

Similar to the type $C$ case, we
will introduce the {\em extended quantum affine algebra} in types $B$ and $D$
defined by an $R$-matrix presentation. We prove
an embedding theorem which will allow us to regard
the extended algebra of rank $n-1$ as a subalgebra of the
corresponding algebra of rank $n$. We also
produce a Drinfeld-type presentation for the extended quantum affine algebra
and give explicit formulas for generators
of its center. It appears to be very likely that these formulas can be included
in a general scheme as developed by Wendlandt~\cite{w:rm} in the Yangian context.

To state our isomorphism theorem, let $\g=\oa_N$ be the orthogonal Lie algebra, where odd and even values $N=2n+1$ and $N=2n$ respectively correspond to the simple Lie algebras of types~$B_n$ and~$D_n$. Choose their simple roots in the form
\begin{gather*}
\al_i =\ep_i-\ep_{i+1}\qquad\text{for}\quad i=1,\dots,n-1,\\
\al_n =\begin{cases}\ep_n  &\text{if}\quad\g=\oa_{2n+1},\\
\ep_{n-1}+\ep_n &\text{if}\quad\g=\oa_{2n},
\end{cases}
\end{gather*}
where
$\ep_1,\dots,\ep_n$ is an orthonormal basis of a Euclidian space with the inner product
$(\cdot\ts,\cdot)$. The Cartan matrix $[A_{ij}]$ is defined by
\begin{gather}\label{cartan}
A_{ij}=\frac{2(\al_i,\al_j)}{(\al_i,\al_i)}.
\end{gather}
For a variable $q$ we set $q_i=q^{r_i}$ for $i=1,\dots,n$, where $r_i=(\al_i,\al_i)/2$.
We will use the standard notation
\begin{gather}\label{kq}
[k]_q=\frac{q^k-q^{-k}}{q-q^{-1}}
\end{gather}
for a nonnegative integer $k$, and
\begin{gather*}
[k]_q!=\prod_{s=1}^{k}[s]_q,\qquad {k\brack r}_{q}=\frac{[k]_q!}{[r]_q!\ts[k-r]_q!}.
\end{gather*}

We will take $\CC\big(q^{1/2}\big)$ as the base field to define most of our quantum algebras.
In type $B_n$ we will need its extension obtained by adjoining the square root of $[2]_{q_n}=q^{1/2}+q^{-1/2}$.

The {\em quantum affine algebra $U_q(\wh{\oa}_N)$ in its Drinfeld presentation}
is the associative algebra
with generators
$x_{i,m}^{\pm}$, $a_{i,l}$, $k_{i}^{\pm}$ and $q^{\pm c/2}$ for $i=1,\dots,n$ and
$m,l\in\ZZ$ with $l\ne 0$, subject to the following defining relations:
the elements $q^{\pm c/2}$ are central,
\begin{gather*}
 k_{i}k_i^{-1} =k_i^{-1}k_{i}=1, \qquad q^{c/2}q^{-c/2}=q^{-c/2}q^{c/2}=1,\\
 k_ik_j =k_jk_i, \qquad k_i\ts a_{j,k}=a_{j,k}\ts k_i,\qquad
 k_i\ts x_{j,m}^{\pm}\ts k_i^{-1}=q_i^{\pm A_{ij}}x_{j,m}^{\pm},
\\
 [a_{i,m},a_{j,l}] =\delta_{m,-l}\ts
 \frac{[mA_{ij}]_{q_i}}{m}\ts\frac{q^{mc}-q^{-mc}}{q_j-q_j^{-1}},\\
 [a_{i,m}, x_{j,l}^{\pm}] =\pm \frac{[mA_{ij}]_{q_i}}{m}\ts q^{\mp |m|c/2}\ts x^{\pm}_{j,m+l},\\
 x^{\pm}_{i,m+1}x^{\pm}_{j,l}-q_i^{\pm A_{ij}}x^{\pm}_{j,l}x^{\pm}_{i,m+1} =
 q_i^{\pm A_{ij}}x^{\pm}_{i,m}x^{\pm}_{j,l+1}-x^{\pm}_{j,l+1}x^{\pm}_{i,m},\\
 [x_{i,m}^{+},x_{j,l}^{-}] =\delta_{ij}\frac{q^{(m-l)\ts c/2}\ts\psi_{i,m+l}
 -q^{-(m-l)\ts c/2}\ts\varphi_{i,m+l}}{q_i-q_i^{-1}},\\
 \sum_{\pi\in \Sym_{r}}\sum_{l=0}^{r}(-1)^l{{r}\brack{l}}_{q_i}
 x^{\pm}_{i,s_{\pi(1)}}\cdots x^{\pm}_{i,s_{\pi(l)}}
 x^{\pm}_{j,m}x^{\pm}_{i,s_{\pi(l+1)}}\cdots x^{\pm}_{i,s_{\pi(r)}}=0, \qquad i\neq j,
\end{gather*}
where in the last relation we set $r=1-A_{ij}$. The elements
$\psi_{i,m}$ and $\varphi_{i,-m}$ with $m\in \ZZ_+$ are defined by
\begin{gather*}
\psi_i(u) :=\sum_{m=0}^{\infty}\psi_{i,m}u^{-m}=k_i\ts\exp\left(\big(q_i-q_i^{-1}\big)
\sum_{s=1}^{\infty}a_{i,s}u^{-s}\right),\\
\varphi_{i}(u) :=\sum_{m=0}^{\infty}\varphi_{i,-m}u^{m}=k_i^{-1}\exp\left({-}\big(q_i-q_i^{-1}\big)
\sum_{s=1}^{\infty}a_{i,-s}u^{s}\right),
\end{gather*}
whereas $\psi_{i,m}=\varphi_{i,-m}=0$ for $m<0$.

To introduce the $R$-matrix presentation of the quantum affine algebra we will use the
endomorphism algebra $\End\big(\CC^{N}\ot\CC^{N}\big)\cong \End\CC^{N}\ot\End\CC^{N}$.
For $\g=\oa_{2n+1}$ consider the following elements of the endomorphism algebra
(extended over $\CC\big(q^{1/2}\big)$):
\begin{gather*}
P =\sum_{i,j=1}^{N}e_{ij}\otimes e_{ji},\qquad
Q=\sum_{i,j=1}^{N} q^{\bi-\bj}\ts e_{i'j'}\otimes e_{ij}
\end{gather*}
and
\begin{gather*}
R =q\ts\sum_{i=1, i\neq i'}^{N}e_{ii}\otimes e_{ii}+e_{n+1,n+1}\ot e_{n+1,n+1}
+\sum_{i\neq j,j'}e_{ii}\otimes e_{jj}
+q^{-1}\sum_{i\neq i'}e_{ii}\otimes e_{i'i'}\\
\hphantom{R =}{} +\big(q-q^{-1}\big)\sum_{i<j}e_{ij}\otimes
e_{ji}-\big(q-q^{-1}\big)\sum_{i>j}q^{\bi-\bj}\ts e_{i'j'}\otimes e_{ij},
\end{gather*}
where $e_{ij}\in\End\CC^{N}$ are the matrix units, and we used the notation $i'=N+1-i$ and
\begin{gather*}
\big(\ts\overline{1},\overline{2},\dots,\overline{N}\ts\big)
=\left(n-\frac{1}{2},\dots,\frac{3}{2},\frac{1}{2},0,-\frac{1}{2},-\frac{3}{2},\dots,-n+\frac{1}{2}\right).
\end{gather*}
In the case $\g=\oa_{2n}$ we define the elements $P$, $Q$ and $R$ by the same formulas
by taking $N=2n$, except that the second term $e_{n+1,n+1}\ot e_{n+1,n+1}$ in
the expression for $R$ should be omitted, while the barred symbols are now given by
\begin{gather*}
\big(\ts\overline{1},\overline{2},\dots,\overline{N}\ts\big)
= (n-1,\dots,1,0,0,-1,\dots,-n+1).
\end{gather*}

In both cases, following~\cite{fri:qa} consider the formal power series
\begin{gather*}
f(u)=1+\sum_{k=1}^{\infty}f_k\tss u^k,
\end{gather*}
whose coefficients $f_k$ are rational functions in $q$ uniquely determined by the relation
\begin{gather}\label{furel}
f(u)\tss f(u\xi)=\frac{1}{\big(1-u\tss q^{-2}\big)\big(1-u\tss q^{2}\big)(1-u\tss\xi)\big(1-u\tss\xi^{-1}\big)},
\end{gather}
where $\xi=q^{2-N}$. Equivalently, $f(u)$ is given by the infinite product formula
\begin{gather}\label{fu}
f(u)=\prod_{r=0}^{\infty}\frac{\big(1-u\ts\xi^{2r}\big)\big(1-u\ts q^{-2}\ts\xi^{2r+1}\big)
\big(1-u\ts q^2\ts\xi^{2r+1}\big)\big(1-u\ts\xi^{2r+2}\big)}
{\big(1-u\ts\xi^{2r-1}\big)\big(1-u\ts\xi^{2r+1}\big)
\big(1-u\ts q^2\xi^{2r}\big)\big(1-u\ts q^{-2}\xi^{2r}\big)}.
\end{gather}
In accordance with \cite{j:qr}, the $R$-{\em matrix} $R(u)$ given by
\begin{gather}\label{ru}
R(u)=f(u)\big(q^{-1}(u-1)(u-\xi)R-\big(q^{-2}-1\big)(u-\xi)P+\big(q^{-2}-1\big)(u-1)\ts \xi\ts Q\big)
\end{gather}
is a solution of the {\em Yang--Baxter equation}
\begin{gather*}
R_{12}(u)\ts R_{13}(u\tss v)\ts R_{23}(v)=R_{23}(v)\ts R_{13}(u\tss v)\ts R_{12}(u).
\end{gather*}

The associative algebra $U^{R}_q(\wh{\oa}_{N})$
is generated by an invertible central element $q^{c/2}$ and
elements ${l}^{\pm}_{ij}[\mp m]$ with $1\leqslant i,j\leqslant N$ and $m\in \ZZ_{+}$
subject to the following defining relations. We have
\begin{gather*}
{l}_{ij}^{+}[0]={l}^{-}_{ji}[0]=0\quad\text{for}\quad i>j
\qquad \text{and} \qquad  {l}^{+}_{ii}[0]\ts {l}_{ii}^{-}[0]={l}_{ii}^{-}[0]\ts{l}^{+}_{ii}[0]=1,
\end{gather*}
while the remaining relations will be written in terms of the formal power series
\begin{gather}\label{liju}
{l}^{\pm}_{ij}(u)=\sum_{m=0}^{\infty}{l}^{\pm}_{ij}[\mp m]\ts u^{\pm m},
\end{gather}
which we combine into the respective matrices
\begin{gather*}
{L}^{\pm}(u)=\sum\limits_{i,j=1}^{N}{l}_{ij}^{\pm}(u)\ot e_{ij}\in
 U^{R}_q(\wh{\oa}_{N})\big[\big[u,u^{-1}\big]\big]\ot\End\CC^{N}.
\end{gather*}
Consider the tensor product algebra $\End\CC^{N}\ot\End\CC^{N}\ot U^{R}_q(\wh{\oa}_{N})$
and introduce the series with coefficients in this algebra by
\begin{gather}\label{lonetwo}
{L}^{\pm}_1(u)=\sum\limits_{i,j=1}^{N}{l}_{ij}^{\pm}(u)\ot e_{ij}\ot 1
\qquad \text{and}\qquad
{L}^{\pm}_2(u)=\sum\limits_{i,j=1}^{N}{l}_{ij}^{\pm}(u)\ot 1 \ot e_{ij}.
\end{gather}
The defining relations then take the form
\begin{gather}\label{rllss}
{R}(u/v)L^{\pm}_{1}(u)L^{\pm}_2(v)=L^{\pm}_2(v)L^{\pm}_{1}(u){R}(u/v),\\
{R}(uq^c/v)L^{+}_1(u)L^{-}_2(v)=L^{-}_2(v)L^{+}_1(u){R}(uq^{-c}/v),\label{rllmp}
\end{gather}
together with the relations
\begin{gather}\label{unitary}
{L}^{\pm}(u)D{L}^{\pm}(u\ts\xi)^{\tra}D^{-1}=1,
\end{gather}
where $\tra$ denotes the matrix
transposition with $e_{ij}^{\tra}= e_{j',i'}$ and $D$ is the diagonal matrix
\begin{gather}\label{D}
D=\diag\ts\big[q^{\overline{1}},\dots,q^{\overline{N}}\tss\big].
\end{gather}

Now apply the {\em Gauss decomposition} to the matrices ${L}^{+}(u)$ and ${L}^{-}(u)$.
There exist unique matrices of the form
\begin{gather*}
F^{\pm}(u)=\begin{bmatrix}
1&0&\dots&0\ts\\
f^{\pm}_{21}(u)&1&\dots&0\\
\vdots&\vdots&\ddots&\vdots\\
f^{\pm}_{N\ts1}(u)&f^{\pm}_{N\ts2}(u)&\dots&1
\end{bmatrix},
\qquad
E^{\pm}(u)=\begin{bmatrix}
\ts1&e^{\pm}_{12}(u)&\dots&e^{\pm}_{1\ts N}(u)\ts\\
\ts0&1&\dots&e^{\pm}_{2\ts N}(u)\\
\vdots&\vdots&\ddots&\vdots\\
0&0&\dots&1
\end{bmatrix},
\end{gather*}
and $H^{\pm}(u)=\diag\ts\big[h^{\pm}_1(u),\dots,h^{\pm}_{N}(u)\big]$,
such that
\begin{gather}\label{gaussdec}
L^{\pm}(u)=F^{\pm}(u)H^{\pm}(u)E^{\pm}(u).
\end{gather}
Set
\begin{gather*}
X^{+}_i(u)=e^{+}_{i,i+1}\big(u\tss q^{c/2}\big)-e_{i,i+1}^{-}\big(u\tss q^{-c/2}\big),
\qquad X^{-}_i(u)=f^{+}_{i+1,i}\big(u\tss q^{-c/2}\big)-f^{-}_{i+1,i}\big(u\tss q^{c/2}\big),
\end{gather*}
for $i=1,\dots, n-1$, and
\begin{gather*}
X^{+}_n(u) =\begin{cases}
 e^{+}_{n,n+1}\big(u\tss q^{c/2}\big)-e_{n,n+1}^{-}\big(u\tss q^{-c/2}\big) &\text{for type $B_n$}, \\
 e^{+}_{n-1,n+1}\big(u\tss q^{c/2}\big)-e_{n-1,n+1}^{-}\big(u\tss q^{-c/2}\big) &\text{for type $D_n$},
 \end{cases}\\
X^{-}_n(u) =\begin{cases}
 f^{+}_{n+1,n}\big(u\tss q^{-c/2}\big)-f^{-}_{n+1,n}\big(u\tss q^{c/2}\big) &\text{for type $B_n$}, \\
 f^{+}_{n+1,n-1}\big(u\tss q^{-c/2}\big)-f^{-}_{n+1,n-1}\big(u\tss q^{c/2}\big) &\text{for type $D_n$}.
 \end{cases}
\end{gather*}

Combine the generators $x^{\pm}_{i,m}$ of the algebra $U_q(\wh{\oa}_{N})$ into the series
\begin{gather*}
x^{\pm}_{i}(u)=\sum_{m\in\ZZ}x^{\pm}_{i,m}\ts u^{-m}.
\end{gather*}

\begin{mthm}
The maps $q^{c/2}\mapsto q^{c/2}$,
\begin{gather*}
x^{\pm}_{i}(u) \mapsto \big(q_i-q_i^{-1}\big)^{-1}X^{\pm}_i\big(uq^i\big),\\
\psi_{i}(u) \mapsto h^{-}_{i+1}\big(uq^i\big)\ts h^{-}_{i}\big(uq^i\big)^{-1},\\
\varphi_{i}(u) \mapsto h^{+}_{i+1}\big(uq^i\big)\ts h^{+}_{i}\big(uq^i\big)^{-1},
\end{gather*}
for $i=1,\dots,n-1$, and
\begin{gather*}
x^{\pm}_{n}(u)\mapsto
\begin{cases}\big(q_n-q_n^{-1}\big)^{-1}[2]_{q_n}^{-1/2}X^{\pm}_n\big(uq^{n}\big)
 &\text{for type $B_n$}, \\
\big(q_n-q_n^{-1}\big)^{-1}X^{\pm}_n\big(uq^{n-1}\big) &\text{for type $D_n$},
\end{cases}
\\
\psi_{n}(u)\mapsto
\begin{cases}
h^{-}_{n+1}\big(uq^{n}\big)\ts h^{-}_{n}\big(uq^{n}\big)^{-1} &\text{for type $B_n$}, \\
h^{-}_{n+1}\big(uq^{n-1}\big)\ts h^{-}_{n-1}\big(uq^{n-1}\big)^{-1} &\text{for type $D_n$},
\end{cases}
\\
\varphi_{n}(u)\mapsto
\begin{cases}
h^{+}_{n+1}\big(uq^{n}\big)\ts h^{+}_{n}\big(uq^{n}\big)^{-1} &\text{for type $B_n$}, \\
h^{+}_{n+1}\big(uq^{n-1}\big)\ts h^{+}_{n-1}\big(uq^{n-1}\big)^{-1} &\text{for type $D_n$},
\end{cases}
\end{gather*}
define an isomorphism $U_q(\wh{\oa}_{N})\to U^{R}_q(\wh{\oa}_{N})$.
\end{mthm}

To prove the Main Theorem we embed $U_q(\wh{\oa}_{N})$ into an extended
quantum affine algebra $U^{\ext}_q(\wh{\oa}_{N})$ which is defined by
a Drinfeld-type presentation. The next step is to use the Gauss decomposition to
construct a homomorphism from
the extended quantum affine algebra to the
algebra~$U(R)$ which is defined by the same
presentation as the algebra $U^{R}_q(\wh{\oa}_{N})$, except that the relation
\eqref{unitary} is omitted. The expressions on the left hand side of~\eqref{unitary},
considered in the algebra $U(R)$, turn out to be scalar matrices,
\begin{gather*}
{L}^{\pm}(u)D{L}^{\pm}(u\xi)^{\tra}D^{-1}={z}^{\pm}(u)\ts 1,
\end{gather*}
for certain formal series ${z}^{\pm}(u)$. Moreover, all coefficients of these series
are central in~$U(R)$. We will give explicit formulas for ${z}^{\pm}(u)$,
regarded as series with coefficients in the algebra $U^{\ext}_q(\wh{\oa}_{N})$,
in terms of its Drinfeld generators.
The quantum affine algebra $U_q(\wh{\oa}_{N})$ can therefore be considered
as the quotient of $U^{\ext}_q(\wh{\oa}_{N})$ by the relations~${z}^{\pm}(u)=1$.

As a final step, we construct the inverse map $U(R)\to U^{\ext}_q(\wh{\oa}_{N})$
by using the universal $R$-matrix for the quantum affine algebra and producing
the associated $L$-operators corresponding to the vector representation of the
algebra $U_q(\wh{\oa}_{N})$.

\section{Quantum affine algebras}

Recall the original definition of the quantum affine algebra
$U_q(\wh\g)$ as introduced by Drinfeld~\cite{d:ha} and Jimbo~\cite{j:qd}.
We suppose that $\g$ is a simple Lie algebra over $\CC$ of rank $n$
and $\wh\g$ is the corresponding (untwisted) affine Kac--Moody algebra with the affine
Cartan matrix $[A_{ij}]_{i,j=0}^n$. We let $\al_0,\al_1,\dots,\al_n$
denote the simple roots and use the notation of \cite[Sections~9.1 and 12.2]{cp:gq}
so that $q_i=q^{r_i}$ for $r_i=(\al_i,\al_i)/2$.

\subsection{Drinfeld--Jimbo definition and new realization}\label{subsec:isoDJD}

The {\em quantum affine algebra} $U_q(\wh{\g})$ is a unital associative algebra over $\mathbb{C}\big(q^{1/2}\big)$ with generators $E_{\al_i}$, $F_{\al_i}$ and
$k_i^{\pm 1}$ with $i=0,1,\dots,n$, subject to the defining relations:
\begin{gather*}
k_ik_i^{-1}=k_i^{-1}k_i=1,\qquad k_ik_j=k_ik_j,
\\
k_iE_{\al_j}k_i^{-1}=q^{A_{ij}}_iE_{\al_j},\qquad k_iF_{\al_j}k_i^{-1}=q^{-A_{ij}}_iF_{\al_j},
\\
[E_{\al_i},F_{\al_j}]=\delta_{ij}\frac{k_i-k_i^{-1}}{q_i-q_i^{-1}},
\\
\sum_{r=0}^{1-A_{ij}}(-1)^r{{1-A_{ij}}\brack{r}}_{q_i}
(E_{\al_i})^rE_{\al_j}(E_{\al_i})^{1-A_{ij}-r} =0,\qquad\text{if}\quad i\ne j,\\
\sum_{r=0}^{1-A_{ij}}(-1)^r{{1-A_{ij}}\brack{r}}_{q_i}
(F_{\al_i})^rF_{\al_j}(F_{\al_i})^{1-A_{ij}-r} =0,\qquad\text{if}\quad i\ne j.
\end{gather*}

By using the braid group action, the set of generators of the algebra $U_q(\wh\g)$
can be extended to the set of affine root vectors of the form $E_{\al+k\de}$, $F_{\al+k\de}$,
$E_{(k\de,i)}$ and $F_{(k\de,i)}$, where $\al$ runs over the positive roots of $\g$,
and $\de$ is the basic imaginary root; see \cite{b:bg, bcp:ac} for details.
Moreover, we can introduce $k_{\al}=\prod\limits_{i=0}^{n}k_i^{m_i}$ for every $\al=\sum\limits_{i=0}^{n}m_i\al_i$, $m_i\in \mathbb{Z}$. Especially, we denote $q^C=k_{\de}$.
The root vectors are used in the explicit isomorphism between the
Drinfeld--Jimbo presentation of the algebra $U_q(\wh\g)$
and the ``new realization'' of Drinfeld which goes back to~\cite{d:nr}, while detailed arguments
were given by Beck~\cite{b:bg}; see also \cite[Lemma~1.5]{bcp:ac}.
In particular, for the Drinfeld presentation of the algebra $U_q(\wh{\oa}_{N})$
given in the Introduction, we find that the isomorphism between these presentations is given by
\begin{alignat*}{4}
& q^{c/2} \mapsto q^{C/2},\qquad
x^{+}_{ik}\mapsto o(i)^kE_{\al_i+k\de}, \qquad && x^{-}_{i,-k} \mapsto o(i)^kF_{\al_i+k\de},\qquad&& k\geqslant 0,& \\
& x^{+}_{i,-k} \mapsto -o(i)^kF_{-\al_i+k\de}\ts k_i^{-1}q^{kC},\qquad && x^{-}_{i,k}
 \mapsto -o(i)^kq^{-kC}\ts k_i\ts E_{-\al_i+k\de},\qquad&& k> 0,&\\
&a_{i,k}\mapsto o(i)^kq^{-kC/2}E_{(k\de,i)}, \qquad && a_{i,-k}
 \mapsto o(i)^k\tss F_{(k\de,i)}q^{kC/2},\qquad && k> 0,&
\end{alignat*}
where $o\colon \{1,2,\dots,n\}\rightarrow \{\pm 1\}$ is a map such that $o(i)=-o(j)$ whenever $A_{ij}<0$.

\subsection{Extended quantum affine algebra}\label{subsec:eqaa}

We will embed the algebra $U_q(\wh{\oa}_{N})$ into an extended quantum affine algebra which we denote by
$U^{\ext}_q(\wh{\oa}_{N})$; cf.~\cite{df:it,fm:ha} and \cite{jlm:ib-c}.
Recalling the scalar function $f(u)$ defined by \eqref{furel} and \eqref{fu} set
\begin{gather}\label{gu}
g(u)=f(u)\big(u-q^{-2}\big)(u-\xi).
\end{gather}
To make formulas look simpler, for variables of type $u$, $v$, or similar, we will use the notation $u_{\pm}=u\tss q^{\pm c/2}$, $v_{\pm}=v\tss q^{\pm c/2}$, etc.

\bde\label{def:eqaa}
The {\em extended quantum affine algebra} $U^{\ext}_q(\wh{\oa}_{N})$ is an associative algebra
over $\mathbb{C}\big(q^{1/2}\big)$
with generators $X^{\pm}_{i,k}$, $h^{+}_{j,m}$, $h^{-}_{j,-m}$ and $q^{c/2}$, where the subscripts take
values $i=1,\dots,n$ and $k\in\ZZ$, while $j=1,\dots,n+1$ and $m\in \ZZ_{+}$.
The defining relations are written with the use of generating functions in a formal
variable $u$:
\begin{gather*}
X^{\pm}_i(u)=\sum_{k\in\ZZ}X^{\pm}_{i,k}\ts u^{-k},\qquad
h^{\pm}_i(u)=\sum_{m=0}^{\infty}h^{\pm}_{i,\mp m}\ts u^{\pm m},
\end{gather*}
they take the following form. The element $q^{c/2}$ is central and invertible,
\begin{gather*}
h_{i,0}^{+}h_{i,0}^{-}=h_{i,0}^{-}h_{i,0}^{+}=1.
\end{gather*}

{\bf Type $B$:} For the relations involving $h^{\pm}_i(u)$ we have
\begin{gather*}
h^{\pm}_i(u)h^{\pm}_j(v) =h^{\pm}_j(v)h^{\pm}_i(u),\qquad h_{n+1,0}^{\pm}=1,\\
g\big(\big(uq^c/v\big)^{\pm 1}\big)\ts h^{\pm}_i(u)h^{\mp}_i(v)
=g\big(\big(uq^{-c}/v\big)^{\pm 1}\big)\ts h^{\mp}_i(v)h^{\pm}_i(u),\qquad i=1,\dots, n,\\
g\big(\big(uq^c/v\big)^{\pm 1}\big)\ts\frac{u_{\pm}-v_{\mp}}{qu_{\pm}-q^{-1}v_{\mp}}h^{\pm}_i(u)h^{\mp}_j(v) = g\big(\big(uq^{-c}/v\big)^{\pm 1}\big)\ts\frac{u_{\mp}-v_{\pm}}{qu_{\mp}-q^{-1}v_{\pm}}h^{\mp}_j(v)h^{\pm}_i(u),
\end{gather*}
for $i<j$, while
\begin{gather*}
g\big(\big(uq^c/v\big)^{\pm 1}\big) \ts\frac{q^{-1}u_{\pm}-qv_{\mp}}{qu_{\pm}-q^{-1}v_{\mp}}\ts
\frac{q^{1/2}u_{\pm}-q^{-1/2}v_{\mp}}{q^{-1/2}u_{\pm}-q^{1/2}v_{\mp}}\ts
h^{\pm}_{n+1}(u)h^{\mp}_{n+1}(v)\\
 \qquad{} =g\big(\big(uq^{-c}/v\big)^{\pm 1}\big)\ts\frac{q^{-1}u_{\mp}-qv_{\pm}}{qu_{\mp}-q^{-1}v_{\pm}}\ts
\frac{q^{1/2}u_{\mp}-q^{-1/2}v_{\pm}}{q^{-1/2}u_{\mp}-q^{1/2}v_{\pm}}\ts
h^{\mp}_{n+1}(v)h^{\pm}_{n+1}(u).
\end{gather*}
The relations
involving $h^{\pm}_i(u)$ and $X_{j}^{\pm}(v)$ are
\begin{gather*}
h_{i}^{\pm}(u)X_{j}^{+}(v) =\frac{u-v_{\pm}}{q^{(\ep_i,\alpha_j)}u-q^{-(\ep_i,\alpha_j)}v_{\pm}}
X_{j}^{+}(v)h_{i}^{\pm}(u),\\
h_{i}^{\pm}(u)X_{j}^{-}(v) =\frac{q^{-(\ep_i,\alpha_j)}u_{\pm}-q^{(\ep_i,\alpha_j)}v}{u_{\pm}-v}
X_{j}^{-}(v)h_{i}^{\pm}(u)
\end{gather*}
for $i\neq n+1$, together with
\begin{gather*}
h_{n+1}^{\pm}(u)X_n^{+}(v) =
\frac{(qu_{\mp}-v)(u_{\mp}-v)}{(u_{\mp}-qv)\big(qu_{\mp}-q^{-1}v\big)}X_n^{+}(v)h_{n+1}^{\pm}(u),\\
h_{n+1}^{\pm}(u)X_n^{-}(v) =
\frac{(u_{\pm}-qv)\big(qu_{\pm}-q^{-1}v\big)}{(qu_{\pm}-v)(u_{\pm}-v)}X_n^{-}(v)h_{n+1}^{\pm}(u),
\end{gather*}
and
\begin{gather*}
h_{n+1}^{\pm}(u)X_{i}^{+}(v) =X_{i}^{+}(v)h_{n+1}^{\pm}(u),\qquad
h_{n+1}^{\pm}(u)X_{i}^{-}(v) =X_{i}^{-}(v)h_{n+1}^{\pm}(u),
\end{gather*}
for $1\leqslant i\leqslant n-1$. For the relations involving $X^{\pm}_i(u)$ we have
\begin{gather*}
\big(u-q^{\pm (\alpha_i,\alpha_j)}v\big)X_{i}^{\pm}\big(uq^{i}\big)X_{j}^{\pm}\big(vq^{j}\big)
=\big(q^{\pm (\alpha_i,\alpha_j)}u-v\big) X_{j}^{\pm}\big(vq^{j}\big)X_{i}^{\pm}\big(uq^{i}\big)
\end{gather*}
for $1\leqslant i,j\leqslant n$,
and
\begin{gather*}\begin{split}&
[X_i^{+}(u),X_j^{-}(v)]=
\delta_{ij}\big(q-q^{-1}\big)\big(\delta\big(u\ts q^{-c}/v\big)h_i^{-}(v_+)^{-1}h_{i+1}^{-}(v_+)\\
& \hphantom{[X_i^{+}(u),X_j^{-}(v)]=}{}
-\delta\big(u\ts q^{c}/v\big)h_i^{+}(u_+)^{-1}h_{i+1}^{+}(u_+)\big)\end{split}
\end{gather*}
together with the Serre relations
\begin{gather*}
\sum_{\pi\in \Sym_{r}}\sum_{l=0}^{r}(-1)^l{{r}\brack{l}}_{q_i}
 X^{\pm}_{i}(u_{\pi(1)})\cdots X^{\pm}_{i}(u_{\pi(l)})
 X^{\pm}_{j}(v)\tss X^{\pm}_{i}(u_{\pi(l+1)})\cdots X^{\pm}_{i}(u_{\pi(r)})=0,
\end{gather*}
which hold for all $i\neq j$ and we set $r=1-A_{ij}$. Here we used the notation
\begin{gather*}
\de(u)=\sum_{r\in\ZZ}\tss u^r
\end{gather*}
for the {\em formal $\delta$-function}.

{\bf Type $D$:}
For the relations involving $h^{\pm}_i(u)$ we have
\begin{gather*}
h^{\pm}_i(u)h^{\pm}_j(v) =h^{\pm}_j(v)h^{\pm}_i(u),\qquad h^{\pm}_{n,0}h_{n+1,0}^{\pm}=1,\\
g\big((uq^c/v)^{\pm 1}\big)\ts h^{\pm}_i(u)h^{\mp}_i(v)
=g\big(\big(uq^{-c}/v\big)^{\pm 1}\big)\ts h^{\mp}_i(v)h^{\pm}_i(u),\qquad i=1,\dots, n+1,
\end{gather*}
and
\begin{gather*}
g\big(\big(uq^c/v\big)^{\pm 1}\big)\ts \frac{q^{-1}u_{\pm}-qv_{\mp}}{qu_{\pm}-q^{-1}v_{\mp}}\ts
\frac{u_{\pm}-v_{\mp}}{u_{\pm}-q^{-1}v_{\mp}}
h^{\pm}_{n}(u)h^{\mp}_{n+1}(v)\\
\qquad {} =g\big(\big(uq^{-c}/v\big)^{\pm 1}\big)\ts
\frac{q^{-1}u_{\mp}-qv_{\pm}}{qu_{\mp}-q^{-1}v_{\pm}}\ts
\frac{u_{\mp}-v_{\pm}}{u_{\mp}-q^{-1}v_{\pm}}
h^{\mp}_{n+1}(v)h^{\pm}_{n}(u)
\end{gather*}
together with
\begin{gather*}
g\big(\big(uq^c/v\big)^{\pm 1}\big)
\ts\frac{u_{\pm}-v_{\mp}}{qu_{\pm}-q^{-1}v_{\mp}}h^{\pm}_i(u)h^{\mp}_j(v)
  =g\big(\big(uq^{-c}/v\big)^{\pm 1}\big)\ts\frac{u_{\mp}-v_{\pm}}{qu_{\mp}-q^{-1}v_{\pm}}h^{\mp}_j(v)h^{\pm}_i(u)
\end{gather*}
for $i<j$ and $(i,j)\neq (n,n+1)$. The relations
involving $h^{\pm}_i(u)$ and $X_{j}^{\pm}(v)$ are
\begin{gather*}
h_{i}^{\pm}(u)X_{j}^{+}(v)
=\frac{u-v_{\pm}}{q^{(\ep_i,\alpha_j)}u-q^{-(\ep_i,\alpha_j)}v_{\pm}}
X_{j}^{+}(v)h_{i}^{\pm}(u),\\
h_{i}^{\pm}(u)X_{j}^{-}(v)
=\frac{q^{(\ep_i,\alpha_j)}u_{\pm}-q^{-(\ep_i,\alpha_j)}v}{u_{\pm}-v}
X_{j}^{-}(v)h_{i}^{\pm}(u)
\end{gather*}
for $i\neq n+1$, together with
\begin{gather*}
h_{n+1}^{\pm}(u)X_n^{+}(v) =
\frac{u_{\mp}-v}{q^{-1}u_{\mp}-qv}X_n^{+}(v)h_{n+1}^{\pm}(u),\\
h_{n+1}^{\pm}(u)X_n^{-}(v) =
\frac{q^{-1}u_{\pm}-qv}{u_{\pm}-v}X_n^{-}(v)h_{n+1}^{\pm}(u),
\end{gather*}
and
\begin{gather*}
h_{n+1}^{\pm}(u)X_{n-1}^{+}(v) =
\frac{u_{\mp}-v}{qu_{\mp}-q^{-1}v}X_{n-1}^{+}(v)h_{n+1}^{\pm}(u),\\
h_{n+1}^{\pm}(u)X_{n-1}^{-}(v) =
\frac{qu_{\pm}-q^{-1}v}{u_{\pm}-v}X_{n-1}^{-}(v)h_{n+1}^{\pm}(u),
\end{gather*}
while
\begin{gather*}
h_{n+1}^{\pm}(u)X_{i}^{+}(v)
 =X_{i}^{+}(v)h_{n+1}^{\pm}(u),\qquad
h_{n+1}^{\pm}(u)X_{i}^{-}(v) =X_{i}^{-}(v)h_{n+1}^{\pm}(u),
\end{gather*}
for $1\leqslant i\leqslant n-2$. For the relations involving $X^{\pm}_i(u)$ we have
\begin{gather*}
\big(u-q^{\pm (\alpha_i,\alpha_j)}v\big)X_{i}^{\pm}\big(uq^{i}\big)X_{j}^{\pm}\big(vq^{j}\big)
=\big(q^{\pm (\alpha_i,\alpha_j)}u-v\big) X_{j}^{\pm}\big(vq^{j}\big)X_{i}^{\pm}\big(uq^{i}\big)
\end{gather*}
for $i,j=1,\dots,n-1$,
\begin{gather*}
\big(u-q^{\pm (\alpha_i,\alpha_n)}v\big)X_{i}^{\pm}\big(uq^{i}\big)X_{n}^{\pm}\big(vq^{n-1}\big)
=\big(q^{\pm (\alpha_i,\alpha_n)}u-v\big) X_{n}^{\pm}\big(vq^{n-1}\big)X_{i}^{\pm}\big(uq^{i}\big)
\end{gather*}
for $i=1,\dots,n-1$,
\begin{gather*}
\big(u-q^{\pm (\alpha_n,\alpha_n)}v\big)X_{n}^{\pm}(u)X_{n}^{\pm}(v)
=\big(q^{\pm (\alpha_n,\alpha_n)}u-v\big) X_{n}^{\pm}(v)X_{n}^{\pm}(u)
\end{gather*}
and
\begin{gather*}
\big[X_i^{+}(u),X_j^{-}(v)\big]=
\delta_{ij}\big(q-q^{-1}\big)\\
\hphantom{\big[X_i^{+}(u),X_j^{-}(v)\big]=}{}\times \big(\delta\big(u\ts q^{-c}/v\big)h_i^{-}(v_+)^{-1}h_{i+1}^{-}(v_+)
-\delta\big(u\ts q^{c}/v\big)h_i^{+}(u_+)^{-1}h_{i+1}^{+}(u_+)\big)
\end{gather*}
together with the
Serre relations
\begin{gather*}
\sum_{\pi\in \Sym_{r}}\sum_{l=0}^{r}(-1)^l{{r}\brack{l}}_{q_i}
 X^{\pm}_{i}(u_{\pi(1)})\cdots X^{\pm}_{i}(u_{\pi(l)})
 X^{\pm}_{j}(v)\tss X^{\pm}_{i}(u_{\pi(l+1)})\cdots X^{\pm}_{i}(u_{\pi(r)})=0,
\end{gather*}
which hold for all $i\neq j$ and we set $r=1-A_{ij}$.
\ede

Introduce two formal power series $z^+(u)$ and $z^-(u)$ in $u$ and $u^{-1}$, respectively,
with coefficients in the algebra $U^{\ext}_q(\wh{\oa}_{N})$ by
\begin{gather}\label{zpm}
z^{\pm}(u)=\begin{cases}
 \prod\limits_{i=1}^{n}h^{\pm}_{i}\big(u\tss\xi\tss q^{2i}\big)^{-1}
 h^{\pm}_{i}\big(u\tss\xi\tss q^{2i-2}\big)\ts \cdot h^{\pm}_{n+1}(u)\ts h^{\pm}_{n+1}(uq)
 & \text{for type $B$,} \\
 \prod\limits_{i=1}^{n-1}h^{\pm}_{i}\big(u\tss\xi\tss q^{2i}\big)^{-1}
 h^{\pm}_{i}\big(u\tss\xi\tss q^{2i-2}\big)\ts\cdot h^{\pm}_{n}(u)\ts h^{\pm}_{n+1}(u)
 & \text{for type $D$,}
 \end{cases}
\end{gather}
where we keep using the notation $\xi=q^{2-N}$. Note that by the defining relations of Definition~\ref{def:eqaa}, the ordering of the factors in the products is irrelevant.

The following claim is verified in the same way as
for type $C$; see \cite[Section~2.2]{jlm:ib-c}.

\bpr\label{prop:zu}
The coefficients of $z^{\pm}(u)$ are central elements of $U^{\ext}_q(\wh{\oa}_{N})$.
\epr

\bpr\label{prop:embed}
The maps $q^{c/2}\mapsto q^{c/2}$,
\begin{gather*}
x^{\pm}_{i}(u) \mapsto \big(q_i-q_i^{-1}\big)^{-1}X^{\pm}_i\big(uq^i\big),\\
\psi_{i}(u) \mapsto h^{-}_{i+1}\big(uq^i\big)\ts h^{-}_{i}\big(uq^i\big)^{-1},\\
\varphi_{i}(u) \mapsto h^{+}_{i+1}\big(uq^i\big)\ts h^{+}_{i}\big(uq^i\big)^{-1},
\end{gather*}
for $i=1,\dots,n-1$ in both types,
\begin{gather*}
x^{\pm}_{n}(u) \mapsto \big(q_n-q_n^{-1}\big)^{-1}X^{\pm}_n\big(uq^{n}\big),\\
\psi_{n}(u) \mapsto h^{-}_{n+1}\big(uq^{n}\big)\ts h^{-}_{n}\big(uq^{n}\big)^{-1},\\
\varphi_{n}(u) \mapsto h^{+}_{n+1}\big(uq^{n}\big)\ts h^{+}_{n}\big(uq^{n}\big)^{-1},
\end{gather*}
for type $B$, and
\begin{gather*}
x^{\pm}_{n}(u) \mapsto \big(q_n-q_n^{-1}\big)^{-1}X^{\pm}_n\big(uq^{n-1}\big),\\
\psi_{n}(u) \mapsto h^{-}_{n+1}\big(uq^{n-1}\big)\ts h^{-}_{n-1}\big(uq^{n-1}\big)^{-1},\\
\varphi_{n}(u) \mapsto h^{+}_{n+1}\big(uq^{n-1}\big)\ts h^{+}_{n-1}\big(uq^{n-1}\big)^{-1},
\end{gather*}
for type $D$, define an embedding $\varsigma\colon U_q(\wh{\oa}_{N})\hra U^{\ext}_q(\wh{\oa}_{N})$.
\epr

\begin{proof}
As with type $C$ \cite[Section~2.2]{jlm:ib-c},
it is straightforward to check that
the maps define a homomorphism. To show that its kernel
is zero, we extend the algebra $U_q(\wh{\oa}_{N})$ in type $D$ by
adjoining the square roots $(k_{n-1}k_n)^{\pm 1/2}$ and keep using the same notation
for the extended algebra. In both types we
will construct a~homomorphism $\varrho\colon U^{\ext}_q(\wh{\oa}_{N})\to U_q(\wh{\oa}_{N})$ such that
the composition $\varrho\circ\varsigma$ is the identity homomorphism on $U_q(\wh{\oa}_{N})$.

There exist power series
$\ze^{\pm}(u)$ with coefficients in the center of $U^{\ext}_q(\wh{\oa}_{N})$
such that \begin{gather*} \ze^{\pm}(u)\ts\ze^{\pm}(u\tss\xi)=z^{\pm}(u).\end{gather*} Explicitly,
\begin{gather*}
\ze^{\pm}(u)=\prod_{m=0}^{\infty} z^{\pm}\big(u\xi^{-2m-1}\big)z^{\pm}\big(u\xi^{-2m-2}\big)^{-1}.
\end{gather*}
Note that although the formula involves an infinite product,
the coefficients of powers of $u$ turn out to be well-defined
elements of $U^{\ext}_q(\wh{\oa}_{N})$; cf.\ the proof of
Proposition~5.5 in~\cite{jlm:ib-c}. The mappings
$X^{\pm}_i(u)\mapsto X^{\pm}_i(u)$ for $i=1,\dots, n$ and
$h^{\pm}_j(u)\mapsto h^{\pm}_j(u)\ts \ze^{\pm}(u)$
for $j=1,\dots, n+1$
define a homomorphism from the algebra $U^{\ext}_q(\wh{\oa}_{N})$ to itself.
The definition of the series $\ze^{\pm}(u)$ implies that for the images of $h^{\pm}_i(u)$
we have the relation
\begin{gather*}
h^{\pm}_i(u)\ts \ze^{\pm}(u)\ts h^{\pm}_i(u\xi)\ts \ze^{\pm}(u\xi)=
h^{\pm}_i(u)\ts h^{\pm}_i(u\xi)\ts z^{\pm}(u).
\end{gather*}
Hence the property $\varrho\circ\varsigma={\rm id}$ will be satisfied if we define the map $\varrho\colon U^{\ext}_q(\wh{\oa}_{N})\to U_q(\wh{\oa}_{N})$ by
\begin{gather*}
X^{\pm}_i(u)\mapsto \big(q_i-q_i^{-1}\big)\ts x^{\pm}_{i}\big(uq^{-i}\big)\qquad\text{for}\quad i=1,\dots,n-1,
\end{gather*}
and
\begin{gather*}
X^{\pm}_n(u)\mapsto
\big(q_n-q_n^{-1}\big)\ts x^{\pm}_{i}\big(uq^{-n-1}\big),
\end{gather*}
while
\begin{gather*}
h^{\pm}_i(u)\mapsto \al^{\pm}_i(u)\qquad\text{for}\quad i=1,\dots,n+1,
\end{gather*}
where the series $\al^{\pm}_i(u)$ are defined in different ways for types $B$ and $D$ and so we consider these cases separately.

For type $B$ we have
\begin{gather*}
\al^{+}_i(u)\ts \al^{+}_i(u\tss\xi)=
\ts\prod_{k=1}^{n}\vp_{k}\big(u\tss\xi\tss q^k\big)^{-1}
\ts\prod_{k=1}^{i-1}\vp_{k}\big(u\tss\xi\tss q^{-k}\big)\ts\prod_{k=i}^{n}\vp_{k}\big(u\tss q^{-k}\big)^{-1}
\end{gather*}
for $i=1,\dots,n$, and
\begin{gather*}
\al^{+}_{n+1}(u)\ts \al^{+}_{n+1}(u\tss\xi)=
\ts\prod_{k=1}^{n}\vp_{k}\big(u\tss\xi\tss q^k\big)^{-1}
\ts\prod_{k=1}^{n}\vp_{k}\big(u\tss\xi\tss q^{-k}\big).
\end{gather*}
Explicitly, by setting $\tilde{\varphi}_{j}(u)=k_j\varphi_{j}(u)$, we get
\begin{gather*}
\al^{+}_i(u)= \prod_{m=0}^{\infty} \prod_{j=1}^{n}\tilde{\varphi}_{j}\big(u\xi^{-2m}q^j\big)^{-1}
\tilde{\varphi}_{j}\big(u\xi^{-2m-1}q^j\big)
\tilde{\varphi}_{j}\big(u\xi^{-2m-1}q^{-j}\big)^{-1}\tilde{\varphi}_{j}\big(u\xi^{-2m-2}q^{-j}\big)\\
\hphantom{\al^{+}_i(u)=}{} \times \prod_{j=1}^{i-1}\tilde{\varphi}_{j}\big(uq^{-j}\big)\times \prod_{j=i}^{n}k_i
\end{gather*}
for $i=1,\dots,n$, and
\begin{gather*}
\al^{+}_{n+1}(u)=\prod_{m=0}^{\infty} \prod_{j=1}^{n}\tilde{\varphi}_{j}\big(u\xi^{-2m}q^j\big)^{-1}
\tilde{\varphi}_{j}\big(u\xi^{-2m-1}q^j\big)
\tilde{\varphi}_{j}\big(u\xi^{-2m-1}q^{-j}\big)^{-1}\tilde{\varphi}_{j}\big(u\xi^{-2m-2}q^{-j}\big)\\
\hphantom{\al^{+}_{n+1}(u)=}{} \times \prod_{j=1}^{n}\tilde{\varphi}_{j}\big(uq^{-j}\big).
\end{gather*}

In type $D$ we have
\begin{gather*}
\al^{+}_i(u)\ts \al^{+}_i(u\tss\xi)=\vp_{n}\big(u\tss q^{-n+1}\big)^{-1}
\ts\prod_{k=1}^{n-2}\vp_{k}\big(u\tss\xi\tss q^k\big)^{-1}
\ts\prod_{k=1}^{i-1}\vp_{k}\big(u\tss\xi\tss q^{-k}\big)\ts\prod_{k=i}^{n-1}\vp_{k}\big(u\tss q^{-k}\big)^{-1}
\end{gather*}
for $i=1,\dots,n-1$,
\begin{gather*}
\al^{+}_n(u)\ts \al^{+}_n(u\tss\xi)=\vp_{n}\big(u\tss q^{-n+1}\big)^{-1}
\ts\prod_{k=1}^{n-2}\vp_{k}\big(u\tss\xi\tss q^k\big)^{-1}
\ts\prod_{k=1}^{n-1}\vp_{k}\big(u\tss\xi\tss q^{-k}\big)
\end{gather*}
and
\begin{gather*}
\al^{+}_{n+1}(u)\ts \al^{+}_{n+1}(u\tss\xi)=\vp_{n}\big(u\tss\xi\tss q^{-n-1}\big)
\ts\prod_{k=1}^{n-2}\vp_{k}\big(u\tss\xi\tss q^k\big)^{-1}\prod_{k=1}^{n-1}\vp_{k}\big(u\tss\xi\tss q^{-k}\big).
\end{gather*}
Explicitly, by setting $\tilde{\varphi}_{j}(u)=k_j\varphi_{j}(u)$, we get
\begin{gather*}
\al^{+}_i(u)= \prod_{m=0}^{\infty} \prod_{j=1}^{n-2}\tilde{\varphi}_{j}\big(u\xi^{-2m}q^j\big)^{-1}
\tilde{\varphi}_{j}\big(u\xi^{-2m-1}q^j\big)
\tilde{\varphi}_{j}\big(u\xi^{-2m-1}q^{-j}\big)^{-1}\tilde{\varphi}_{j}\big(u\xi^{-2m-2}q^{-j}\big)\\
\hphantom{\al^{+}_i(u)=}{}\times \prod_{m=0}^{\infty} \prod_{j=n-1}^{n}\tilde{\varphi}_{j}\big(u\xi^{-2m}q^{n-1}\big)^{-1}
\tilde{\varphi}_{j}\big(u\xi^{-2m-1}q^{n-1}\big)
 \prod_{j=1}^{i-1}\tilde{\varphi}_{j}\big(uq^{-j}\big) \prod_{j=i}^{n-2}k_j (k_{n-1}k_n)^{1/2}
\end{gather*}
for $i=1,\dots, n-1$,
\begin{gather*}
\al^{+}_n(u)= \prod_{m=0}^{\infty} \prod_{j=1}^{n-2}\tilde{\varphi}_{j}\big(u\xi^{-2m}q^j\big)^{-1}
\tilde{\varphi}_{j}\big(u\xi^{-2m-1}q^j\big)
\tilde{\varphi}_{j}\big(u\xi^{-2m-1}q^{-j}\big)^{-1}\tilde{\varphi}_{j}\big(u\xi^{-2m-2}q^{-j}\big)\\
\hphantom{\al^{+}_n(u)=}{}\times \prod_{m=0}^{\infty} \prod_{j=n-1}^{n}\tilde{\varphi}_{j}\big(u\xi^{-2m}q^{n-1}\big)^{-1}
\tilde{\varphi}_{j}\big(u\xi^{-2m-1}q^{n-1}\big)
  \prod_{j=1}^{n-1}\tilde{\varphi}_{j}\big(uq^{-j}\big) \big(k_{n-1}^{-1}k_n\big)^{1/2},
\end{gather*}
and
\begin{gather*}
\al^{+}_{n+1}(u)= \prod_{m=0}^{\infty} \prod_{j=1}^{n-2}\tilde{\varphi}_{j}\big(u\xi^{-2m}q^j\big)^{-1}
\tilde{\varphi}_{j}\big(u\xi^{-2m-1}q^j\big)
\tilde{\varphi}_{j}\big(u\xi^{-2m-1}q^{-j}\big)^{-1}\tilde{\varphi}_{j}\big(u\xi^{-2m-2}q^{-j}\big)\\
\hphantom{\al^{+}_{n+1}(u)=}{}\times
\prod_{m=0}^{\infty} \prod_{j=n-1}^{n}\tilde{\varphi}_{j}\big(u\xi^{-2m}q^{n-1}\big)^{-1}
\tilde{\varphi}_{j}\big(u\xi^{-2m-1}q^{n-1}\big)\\
\hphantom{\al^{+}_{n+1}(u)=}{}
\times \prod_{j=1}^{n-1}\tilde{\varphi}_{j}\big(uq^{-j}\big) \tilde{\varphi}_{n}\big(uq^{-n+1}\big)
 \big(k_{n-1}^{-1}k_n\big)^{-1/2}.
\end{gather*}
In both types
the relations defining $\al^{-}_i(u)$ are obtained from those above
by the respective replacements $\al^{+}_i(u)\to \al^{-}_i(u)$, $k_i\to k_i^{-1}$
and $\vp_{k}(u)\to\psi_k(u)$. Although the above explicit formulas of $\al^{\pm}_i(u)$
involve infinite products, their coefficients actually belong to $U^{\ext}_q(\wh{\oa}_{N})$.
For instance,
\begin{gather*}
\al_1^{+}(u)=h_{1,0}^{+}  \exp
\left(\sum_{k>0}\sum_{j=1}^n\big(q_j-q_j^{-1}\big)\tilde{B}_{1j}\big(q^k\big)a_{j, -k}u^k\right);
\end{gather*}
see the proof of Proposition 5.5 in~\cite{jlm:ib-c} for more details.

As with type $C$, one can verify directly that
the map $\varrho$ defines a homomorphism or apply the
calculations with Gaussian generators performed below; cf.~\cite[Remark~5.6]{jlm:ib-c}.
\end{proof}

By Proposition~\ref{prop:embed}
we may regard $U_q(\wh{\oa}_{N})$ as a subalgebra
of $U^{\ext}_q(\wh{\oa}_{N})$. In the following corollary
we will keep the same notation for the algebra $U_q(\wh{\oa}_{N})$ in type $D$ extended by
adjoining the square roots $(k_{n-1}k_n)^{\pm 1/2}$ (no extension is needed in type $B$).
Let $\Cc$ be the subalgebra of~$U^{\ext}_q(\wh{\oa}_{N})$ generated by the coefficients
of the series $z^{\pm}(u)$.

\bco\label{cor:decomp}
We have the tensor product decomposition
\begin{gather*}
U^{\ext}_q(\wh{\oa}_{N}) = U_q(\wh{\oa}_{N})\otimes _{\mathbb{C}(q^{1/2})}\Cc.
\end{gather*}
\eco

\begin{proof}The argument is the same as for type $C$ \cite[Section~2.2]{jlm:ib-c}.
\end{proof}

\section[$R$-matrix presentations]{$\boldsymbol{R}$-matrix presentations}\label{sec:nd}

\subsection[The algebras $U(R)$ and $U\big(\overline{R}\big)$]{The algebras $\boldsymbol{U(R)}$ and $\boldsymbol{U\big(\overline{R}\big)}$}

As defined in the introduction, the algebra $U(R)$ is generated by an invertible central element~$q^{c/2}$ and elements ${l}^{\pm}_{ij}[\mp m]$ with $1\leqslant i,j\leqslant N$ and $m\in \ZZ_{+}$ such that
\begin{gather*}
{l}_{ij}^{+}[0]={l}^{-}_{ji}[0]=0\quad\text{for}\quad i>j
\qquad \text{and} \qquad  {l}^{+}_{ii}[0]\ts {l}_{ii}^{-}[0]={l}_{ii}^{-}[0]\ts{l}^{+}_{ii}[0]=1,
\end{gather*}
and the remaining relations \eqref{rllss} and \eqref{rllmp} (omitting \eqref{unitary})
written in terms of the formal power series~\eqref{liju}.
We will need another algebra $U\big(\overline{R}\big)$ which is defined in a very similar way,
except that it is associated with a different $R$-matrix $\overline{R}(u)$ instead of \eqref{ru}.
Namely, the two $R$-matrices are related by $R(u)=g(u)\overline{R}(u)$
with $g(u)$ defined in \eqref{gu}, so that
\begin{gather}\label{rbar}
\overline{R}(u)=\frac{u-1}{u\tss q-q^{-1}}\ts R
+\frac{q-q^{-1}}{u\tss q-q^{-1}}\ts P-\frac{\big(q-q^{-1}\big)(u-1)\xi}{\big(u\tss q-q^{-1}\big)(u-\xi)}\ts Q.
\end{gather}
Note the {\em unitarity property}
\begin{gather}\label{unitarityrbar}
\overline{R}_{12}(u)\ts \overline{R}_{21}\big(u^{-1}\big)=1,
\end{gather}
satisfied by
this $R$-matrix, where $\overline{R}_{12}(u)=\overline{R}(u)$ and
$\overline{R}_{21}(u)=P\overline{R}(u)P$. More explicitly the
$R$-matrix $\overline{R}(u)$ can be written in the form
\begin{gather}
\overline{R}(u) =\sum_{i=1,\ts i\neq i'}^{N}e_{ii}\otimes e_{ii}+\frac{u-1}{qu-q^{-1}}
\sum_{i\neq j,j'}e_{ii}\otimes e_{jj}
+\frac{q-q^{-1}}{qu-q^{-1}}\sum_{i>j,\ts i\neq j'}e_{ij}\otimes e_{ji}\nonumber\\
\hphantom{\overline{R}(u) =}{} +\frac{\big(q-q^{-1}\big)u}{qu-q^{-1}}\sum_{i<j,\ts i\neq j'}e_{ij}\otimes e_{ji}+
\frac{1}{\big(u-q^{-2}\big)(u-\xi)}\sum_{i,j=1}^{N} a_{ij}(u)\ts e_{i'j'}\otimes e_{ij},\label{rbar2}
\end{gather}
where
\begin{gather*}
a_{ij}(u)=
 \begin{cases}
 \big(q^{-2}u-\xi\big)(u-1)  &\text{for}\quad i=j,\quad i\neq i',\\
 q^{-1}(u-\xi)(u-1)+(\xi-1)\big(q^{-2}-1\big)\tss u &\text{for}\quad i=j,\quad i=i',\\
 \big(q^{-2}-1\big)\big(q^{\bi-\bj}\xi(u-1)-\de_{ij'}(u-\xi)\big) &\text{for}\quad i<j, \\
 \big(q^{-2}-1\big)\tss u\big(q^{\bi-\bj}(u-1)-\de_{ij'}(u-\xi)\big) &\text{for}\quad i>j.
 \end{cases}
\end{gather*}

The {\em algebra} $U\big(\overline{R}\big)$
is generated by an invertible central element $q^{c/2}$ and
elements ${\ell}^{\tss\pm}_{ij}[\mp m]$ with $1\leqslant i,j\leqslant N$ and $m\in \ZZ_{+}$
such that
\begin{gather*}
{\ell}_{ij}^{+}[0]={\ell}^{-}_{ji}[0]=0\quad\text{for}\quad i>j
\qquad \text{and} \qquad  {\ell}^{+}_{ii}[0]\ts {\ell}_{ii}^{-}[0]={\ell}_{ii}^{-}[0]\ts{\ell}^{+}_{ii}[0]=1.
\end{gather*}
Introduce the formal power series
\begin{gather*}
{\ell}^{\tss\pm}_{ij}(u)=\sum_{m=0}^{\infty}{\ell}^{\tss\pm}_{ij}[\mp m]\ts u^{\pm m},
\end{gather*}
which we combine into the respective matrices
\begin{gather*}
{\Lc}^{\pm}(u)=\sum\limits_{i,j=1}^{N}{\ell}^{\tss\pm}_{ij}(u)\ot e_{ij}\in
U\big(\overline{R}\big)\big[\big[u,u^{-1}\big]\big]\ot\End\CC^{N}.
\end{gather*}
The remaining defining relations of the algebra $U\big(\overline{R}\big)$ take the form
\begin{gather}\label{gen rel1}
\overline{R}(u/v)\ts\Lc^{\pm}_{1}(u)\ts\Lc^{\pm}_2(v)=\Lc^{\pm}_2(v)\ts\Lc^{\pm}_{1}(u)\ts\overline{R}(u/v),\\
\overline{R}(u\tss q^c/v)\ts\Lc^{+}_1(u)\ts\Lc^{-}_2(v)=\Lc^{-}_2(v)\ts\Lc^{+}_1(u)\ts\overline{R}(u\tss q^{-c}/v),\label{gen rel2}
\end{gather}
where the subscripts have the same meaning as in \eqref{lonetwo}.
The unitarity property~\eqref{unitarityrbar} implies that
relation \eqref{gen rel2} can be written in the equivalent form
\begin{gather*}
\overline{R}(u\tss q^{-c}/v)\ts\Lc^{-}_1(u)\ts\Lc^{+}_2(v)
=\Lc^{+}_2(v)\ts\Lc^{-}_1(u)\ts\overline{R}(u\tss q^{c}/v).
\end{gather*}

\begin{Remark}\label{rem:gln}
The defining relations satisfied by the series
$\ell^{\tss\pm}_{ij}(u)$
with $1\leqslant i,j\leqslant n$ coincide with those for the quantum affine algebra
$U_q\big(\wh{\gl}_n\big)$ in \cite{df:it}.
\end{Remark}

Following \cite{df:it} and \cite{jlm:ib-c} we will relate the algebras
$U(R)$ and $U\big(\overline{R}\big)$ by using
the Heisenberg algebra $\Hc_q(n)$
with generators $q^c$ and $\be_r$ with $r\in\ZZ\setminus \{0\}$.
The defining relations of $\Hc_q(n)$ have the form
\begin{gather*}
\big[\be_r,\be_s\big]=\de_{r, -s}\ts \al_r,\qquad r\geqslant 1,
\end{gather*}
and $q^c$ is central and invertible. The elements $\al_r$ are defined by the expansion
\begin{gather*}
\exp\ts\sum_{r=1}^{\infty}\al_r u^r=\frac{g(u\tss q^{-c})}{g(u\tss q^{c})}.
\end{gather*}
So we have the identity
\begin{gather*}
g\big(u\ts q^{c}/v\big)\ts \exp\ts\sum_{r=1}^{\infty}\be_r u^r\cdot
\ts \exp\ts\sum_{s=1}^{\infty}\be_{-s} v^{-s}
=g\big(u\ts q^{-c}/v\big)\ts \exp\ts\sum_{s=1}^{\infty}\be_{-s} v^{-s}
\cdot \exp\ts\sum_{r=1}^{\infty}\be_r u^r.
\end{gather*}

\bpr\label{prop:homheis}
The mappings
\begin{gather*}
\Lc^{\ts +}(u)\mapsto \exp\ts\sum_{r=1}^{\infty}\be_{-r} u^{-r}\cdot
L^+(u),\qquad
\Lc^{\ts -}(u)\mapsto \exp\ts\sum_{r=1}^{\infty}\be_{r} u^{r}\cdot
L^-(u),
\end{gather*}
define a homomorphism $U\big(\overline{R}\big)\to \Hc_q(n)\ot_{\CC[q^c,\ts q^{-c}]}U(R)$.
\epr

We will use the notation
$\tra_a$ for the
matrix transposition defined in
\eqref{unitary} applied to the $a$-th copy of the endomorphism algebra
$\End\CC^{N}$ in a multiple tensor product.
Note the following {\em crossing symmetry} relations
satisfied by the $R$-matrices:
\begin{gather*}
 \overline{R}(u)D_1\overline{R}(u\tss\xi)^{\tra_1}D_1^{-1}
  =\frac{\big(u-q^2\big)(u\tss\xi-1)}{(1-u)\big(1-u\tss\xi q^2\big)},\\
 {R}(u)D_1{R}(u\tss\xi)^{\tra_1}D_1^{-1} =\xi^2 q^{-2},
\end{gather*}
where the diagonal matrix $D$ is defined in \eqref{D} and the meaning of the subscripts
is the same as in \eqref{lonetwo}. The next two propositions are verified in the same way as
for type $C$; see \cite[Section~3.1]{jlm:ib-c}.

\bpr\label{prop:central}
In the algebras $U(R)$ and $U\big(\overline{R}\big)$ we have
the relations
\begin{gather*}
D{L}^{\pm}(u\tss\xi)^{\tra}D^{-1}{L}^{\pm}(u)
={L}^{\pm}(u)D{L}^{\pm}(u\tss\xi)^{\tra}D^{-1}={z}^{\pm}(u)\ts1,
\end{gather*}
and
\begin{gather}\label{DLbarDLbar}
D\Lc^{\pm}(u\tss\xi)^{\tra}D^{-1}\Lc^{\pm}(u)
=\Lc^{\pm}(u)D\Lc^{\pm}(u\tss\xi)^{\tra}D^{-1}=\z^{\pm}(u)\ts1,
\end{gather}
for certain series ${z}^{\pm}(u)$ and $\z^{\pm}(u)$ with coefficients in the respective algebra.
\epr

\bpr\label{prop:centr} All coefficients of the series $z^{+}(u)$ and $z^{-}(u)$
belong to the center of the algebra $U(R)$.
\epr

\begin{Remark}\label{rem:noncent}
Although the coefficients of the series $\z^{+}(u)$ and $\z^{-}(u)$ are central
in the respective subalgebras of $U\big(\overline{R}\big)$ generated by the coefficients of the
series $\ell^+_{ij}(u)$ and $\ell^-_{ij}(u)$, they are not central in the entire
algebra $U\big(\overline{R}\big)$.
\end{Remark}

\subsection{Homomorphism theorems}

Now we aim to make a connection between the algebras $U\big(\overline{R}\big)$
associated with the Lie algebras~$\oa_{N-2}$ and~$\oa_{N}$.
We will use quasideterminants as defined in~\cite{gr:dm} and~\cite{gr:tn}.
Let $A=[a_{ij}]$ be a~square matrix over a~ring with~$1$.
Denote by $A^{ij}$ the matrix obtained from $A$
by deleting the $i$-th row and $j$-th column. Suppose that the matrix $A^{ij}$ is invertible. The $ij$-{\em th quasideterminant of}~$A$ is defined by the formula
\begin{gather*}
|A|_{ij}=a_{ij}-r^{\tss j}_i\big(A^{ij}\big)^{-1}\ts c^{\tss i}_j,
\end{gather*}
where $r^{\tss j}_i$ is the row matrix obtained from the $i$-th
row of $A$ by deleting the element~$a_{ij}$, and~$c^{\tss i}_j$
is the column matrix obtained from the $j$-th
column of $A$ by deleting the element $a_{ij}$.
The quasideterminant $|A|_{ij}$ is also denoted
by boxing the entry $a_{ij}$ in the matrix~$A$.

The rank $n$ of the Lie algebra $\oa_N$ with $N=2n+1$ or $N=2n$
will vary so we will indicate the dependence on $n$ by
adding a subscript $[n]$ to the $R$-matrices.
Consider the algebra $U\big(\overline{R}^{\tss[n-1]}\big)$ and let
the indices of the generators $\ell^{\tss\pm}_{ij}[\mp m]$ range over the sets
$2\leqslant i,j\leqslant 2\pr$ and $m=0,1,\dots$, where $i\pr=N-i+1$, as before.

Proofs of the following theorems are not different from those in type $C$;
see \cite[Section~3.3]{jlm:ib-c}.

\bth\label{thm:embed}
The mappings $q^{\pm c/2}\mapsto q^{\pm c/2}$ and
\begin{gather*}
\ell^{\tss\pm}_{ij}(u)\mapsto \left|\begin{matrix}
\ell^{\tss\pm}_{11}(u)&\ell^{\tss\pm}_{1j}(u)\vspace{1mm}\\
\ell^{\tss\pm}_{i1}(u)&\boxed{\ell^{\tss\pm}_{ij}(u)}
\end{matrix}\right|,\qquad 2\leqslant i,j\leqslant 2\pr,
\end{gather*}
define a homomorphism ${U}\big(\overline{R}^{\tss[n-1]}\big)\to {U}\big(\overline{R}^{\tss[n]}\big)$.
\eth

Fix a positive integer $m$ such that
$m< n$. Suppose that the generators $\ell_{ij}^{\tss\pm}(u)$ of the algebra $U\big(\overline{R}^{\tss[n-m]}\big)$ are labelled by the indices
$m+1\leqslant i,j\leqslant (m+1)\pr$.

\bth\label{thm:red}
For $m\leqslant n-1$,
the mapping
\begin{gather}\label{redu}
\ell^{\tss\pm}_{ij}(u)\mapsto \left|\begin{matrix}
\ell^{\tss\pm}_{11}(u)&\dots&\ell^{\tss\pm}_{1m}(u)&\ell^{\tss\pm}_{1j}(u)\\
\dots&\dots&\dots&\dots\\
\ell^{\tss\pm}_{m1}(u)&\dots&\ell^{\tss\pm}_{mm}(u)&\ell^{\tss\pm}_{mj}(u)\vspace{1mm}\\
\ell^{\tss\pm}_{i1}(u)&\dots&\ell^{\tss\pm}_{im}(u)&\boxed{\ell^{\tss\pm}_{ij}(u)}
\end{matrix}\right|,\qquad m+1\leqslant i,j\leqslant (m+1)\pr,
\end{gather}
defines a homomorphism $\psi_m\colon U\big(\overline{R}^{\tss[n-m]}\big)\to U\big(\overline{R}^{\tss[n]}\big)$.
\eth

We also point out a consistence property of the homomorphisms~\eqref{redu}.
Write $\psi_m=\psi^{(n)}_m$ to indicate the dependence of~$n$.
For a parameter $l$ we have the corresponding homomorphism
\begin{gather*}
\psi^{(n-l)}_m\colon \ U\big(\overline{R}^{\tss[n-l-m]}\big)\to U\big(\overline{R}^{\tss[n-l]}\big)
\end{gather*}
provided by \eqref{redu}. Then we have the equality of maps
$\psi^{(n)}_l\circ\psi^{(n-l)}_m=\psi^{(n)}_{l+m}$.

\bco\label{cor:commu}
Under the assumptions of Theorem~{\rm \ref{thm:red}} we have
\begin{gather*}
\big[\ell^{\tss\pm}_{ab}(u),\psi_m\big(\ell^{\tss\pm}_{ij}(v)\big)\big] =0,\\
\frac{u_{\pm}-v_{\mp}}{qu_{\pm}-q^{-1}v_{\mp}}\ell^{\tss\pm}_{ab}(u)\psi_m\big(\ell^{\mp}_{ij}(v)\big)
 =\frac{u_{\mp}-v_{\pm}}{qu_{\mp}-q^{-1}v_{\pm}}\psi_m\big(\ell^{\mp}_{ij}(v)\big)\ell^{\tss\pm}_{ab}(u),
\end{gather*}
for all $1\leqslant a,b\leqslant m$ and $m+1\leqslant i,j\leqslant (m+1)\pr$.
\eco

\section{Gauss decomposition}\label{sec:gauss decom}

Apply the Gauss decompositions~\eqref{gaussdec} to the matrices
$L^{\pm}(u)$ and $\Lc^{\pm}(u)$ associated with the respective algebras $U\big(R^{[n]}\big)$
and $U\big(\overline{R}^{\tss[n]}\big)$. These algebras
are generated by the coefficients of the matrix elements of the triangular and diagonal matrices
which we will refer to as the {\em Gaussian generators}. Here
we produce necessary relations satisfied by these generators to be able to get presentations
of the $R$-matrix algebras $U\big(R^{[n]}\big)$ and $U\big(\overline{R}^{\tss[n]}\big)$.

\subsection{Gaussian generators}

The entries of the matrices $F^{\pm}(u)$, $H^{\pm}(u)$ and $E^{\pm}(u)$ occurring
in the decompositions \eqref{gaussdec} can be described by the universal quasideterminant
formulas \cite{gr:dm,gr:tn}:
\begin{gather}\label{hmqua}
h^{\pm}_i(u)=\begin{vmatrix} l^{\pm}_{1\tss 1}(u)&\dots&l^{\pm}_{1\ts i-1}(u)&l^{\pm}_{1\tss i}(u)\\
 \vdots&\ddots&\vdots&\vdots\\
 l^{\pm}_{i-1\ts 1}(u)&\dots&l^{\pm}_{i-1\ts i-1}(u)&l^{\pm}_{i-1\ts i}(u)\vspace{1mm}\\
 l^{\pm}_{i\tss 1}(u)&\dots&l^{\pm}_{i\ts i-1}(u)&\boxed{l^{\pm}_{i\tss i}(u)}\\
 \end{vmatrix},\qquad i=1,\dots,N,
\end{gather}
whereas
\begin{gather}\label{eijmlqua}
e^{\pm}_{ij}(u)=h^{\pm}_i(u)^{-1}\ts\begin{vmatrix} l^{\pm}_{1\tss 1}(u)&\dots&
l^{\pm}_{1\ts i-1}(u)&l^{\pm}_{1\ts j}(u)\\
 \vdots&\ddots&\vdots&\vdots\\
 l^{\pm}_{i-1\ts 1}(u)&\dots&l^{\pm}_{i-1\ts i-1}(u)&l^{\pm}_{i-1\ts j}(u)\vspace{1mm}\\
 l^{\pm}_{i\tss 1}(u)&\dots&l^{\pm}_{i\ts i-1}(u)&\boxed{l^{\pm}_{i\tss j}(u)}\\
 \end{vmatrix}
\end{gather}
and
\begin{gather}\label{fijlmqua}
f^{\pm}_{ji}(u)=\begin{vmatrix} l^{\pm}_{1\tss 1}(u)&\dots&l^{\pm}_{1\ts i-1}(u)&l^{\pm}_{1\tss i}(u)\\
 \vdots&\ddots&\vdots&\vdots\\
 l^{\pm}_{i-1\ts 1}(u)&\dots&l^{\pm}_{i-1\ts i-1}(u)&l^{\pm}_{i-1\ts i}(u)\vspace{1mm}\\
 l^{\pm}_{j\ts 1}(u)&\dots&l^{\pm}_{j\ts i-1}(u)&\boxed{l^{\pm}_{j\tss i}(u)}\\
 \end{vmatrix}\ts h^{\pm}_i(u)^{-1}
\end{gather}
for $1\leqslant i<j\leqslant N$. The same formulas hold for the expressions
of the entries of the respective triangular matrices $\Fc^{\pm}(u)$ and
$\Ec^{\pm}(u)$ and the diagonal matrices $\Hc^{\pm}(u)=\diag\ts[\h^{\pm}_{i}(u)]$
in terms of the formal series $\ell^{\tss\pm}_{ij}(u)$, which arise from the Gauss decomposition
\begin{gather*}
\Lc^{\pm}(u)=\Fc^{\pm}(u)\ts\Hc^{\pm}(u)\ts\Ec^{\pm}(u)
\end{gather*}
for the algebra $U\big(\overline{R}^{\tss[n]}\big)$. We will denote by
$\e_{ij}(u)$ and $\f_{ji}(u)$ the entries of the respective matrices~$\Ec^{\pm}(u)$ and $\Fc^{\pm}(u)$ for $i<j$.

The following Laurent series with coefficients in
the respective algebras $U\big(R^{[n]}\big)$ and $U\big(\overline{R}^{\tss[n]}\big)$
will be used frequently:
\begin{alignat}{3}
\label{Xi} & X^{+}_i(u)=e^{+}_{i,i+1}(u_{+})-e_{i,i+1}^{-}(u_{-}),
\qquad && X^{-}_i(u)=f^{+}_{i+1,i}(u_{-})-f^{-}_{i+1,i}(u_{+}), & \\
& \Xc^{+}_i(u)=\e^{+}_{i,i+1}(u_{+})-\e_{i,i+1}^{-}(u_{-}),
\qquad && \Xc^{-}_i(u)=\f^{+}_{i+1,i}(u_{-})-\f^{-}_{i+1,i}(u_{+})& \label{Xibar}
\end{alignat}
for $i=1,\dots,n-1$,
 and
\begin{gather}
\label{Xnbarp}
X^{+}_n(u) =\begin{cases}
 e^{+}_{n,n+1}(u_{+})-e_{n,n+1}^{-}(u_{-}) &\text{for type $B$}, \\
 e^{+}_{n-1,n+1}(u_{+})-e_{n-1,n+1}^{-}(u_{-}) &\text{for type $D$},
 \end{cases}\\
\label{Xnbarm}
X^{-}_n(u) =\begin{cases}
 f^{+}_{n+1,n}(u_{-})-f^{-}_{n+1,n}(u_{+}) &\text{for type $B$}, \\
 f^{+}_{n+1,n-1}(u_{-})-f^{-}_{n+1,n-1}(u_{+}) &\text{for type $D$},
 \end{cases}
\end{gather}
while
\begin{gather*}
\Xc^{+}_n(u) =\begin{cases}
 \e^{+}_{n,n+1}(u_{+})-\e_{n,n+1}^{-}(u_{-}) &\text{for type $B$}, \\
 \e^{+}_{n-1,n+1}(u_{+})-\e_{n-1,n+1}^{-}(u_{-}) &\text{for type $D$},
 \end{cases}\\
 \Xc^{-}_n(u) =\begin{cases}
 \f^{+}_{n+1,n}(u_{-})-\f^{-}_{n+1,n}(u_{+}) &\text{for type $B$}, \\
 \f^{+}_{n+1,n-1}(u_{-})-\f^{-}_{n+1,n-1}(u_{+}) &\text{for type $D$}.
 \end{cases}
\end{gather*}

\bpr\label{prop:corrgauss}Under the homomorphism $U\big(\overline{R}\big)\to \Hc_q(n)\ot_{\CC[q^c,\ts q^{-c}]}U(R)$
provided by Proposition~{\rm \ref{prop:homheis}} we have
\begin{gather*}
 \e^{\pm}_{ij}(u)   \mapsto e^{\pm}_{ij}(u), \\
 \f^{\pm}_{ij}(u)   \mapsto f^{\pm}_{ij}(u), \\
 \h^{\pm}_{i}(u)  \mapsto \exp\sum\limits _{k=1}^{\infty}\be_{\mp k}u^{\mp k}\cdot h^{\pm}_{i}(u).
\end{gather*}
\epr

\begin{proof}This is immediate from the formulas for the Gaussian generators.
\end{proof}

Suppose that $0\leqslant m< n$.
We will use the superscript $[n-m]$ to indicate
square submatrices corresponding to rows and columns labelled by
$m+1,m+2,\dots,(m+1)'$. In particular, we set
\begin{gather*}
\Fc^{\pm[n-m]}(u)=\begin{bmatrix}
1&0&\dots&0\ts\\
\f^{\pm}_{m+2\ts m+1}(u)&1&\dots&0\\
\vdots&\ddots&\ddots&\vdots\\
\f^{\pm}_{(m+1)'\tss m+1}(u)&\dots&\f^{\pm}_{(m+1)'\ts (m+2)'}(u)&1
\end{bmatrix},
\\
\Ec^{\pm[n-m]}(u)=\begin{bmatrix} 1&\e^{\pm}_{m+1\tss m+2}(u)&\ldots&\e^{\pm}_{m+1\tss (m+1)'}(u)\\
 0&1&\ddots &\vdots\\
 \vdots&\vdots&\ddots&\e^{\pm}_{(m+2)'\tss(m+1)'}(u)\\
 0&0&\ldots&1\\
 \end{bmatrix}
\end{gather*}
and $\Hc^{\pm[n-m]}(u)=\diag\ts\big[\h^{\pm}_{m+1}(u),\dots,\h^{\pm}_{(m+1)'}(u)\big]$. Furthermore,
introduce the products of these matrices by
\begin{gather*}
\Lc^{\pm[n-m]}(u)=\Fc^{\pm[n-m]}(u)\ts \Hc^{\pm[n-m]}(u)\ts \Ec^{\pm[n-m]}(u).
\end{gather*}
The entries of $\Lc^{\pm[n-m]}(u)$ will be denoted by $\ell^{\tss\pm[n-m]}_{ij}(u)$.

The next series of relations are $B$ and $D$ type counterparts of the corresponding relations
in type $C$ and verified by the same calculations; see \cite[Section~4.2]{jlm:ib-c}.

\bpr\label{prop:gauss-consist}
The series $\ell^{\tss\pm[n-m]}_{ij}(u)$ coincides with the image of
the generator series $\ell^{\tss\pm}_{ij}(u)$
of the extended quantum affine algebra $U\big(\overline{R}^{\tss[n-m]}\big)$
under the homomorphism \eqref{redu},
\begin{gather*}
\ell^{\pm[n-m]}_{ij}(u)=\psi_m\big(\ell^{\pm}_{ij}(u)\big),\qquad m+1\leqslant i,j\leqslant (m+1)'.
\end{gather*}
\epr

\bco\label{cor:guass-embed}
The following relations hold in $U\big(\overline{R}^{\tss[n]}\big)$:
\begin{gather*}
\overline{R}^{\tss[n-m]}_{12}(u/v)\ts \Lc^{\pm[n-m]}_1(u)\ts \Lc^{\pm[n-m]}_2(v)
=\Lc^{\pm[n-m]}_2(v)\ts \Lc^{\pm[n-m]}_1(u)\ts \overline{R}^{\tss[n-m]}_{12}(u/v),
\\
\overline{R}^{\tss[n-m]}_{12}(u_{+}/v_{-})\ts \Lc^{+[n-m]}_1(u)\ts \Lc^{-[n-m]}_2(v)
=\Lc^{-[n-m]}_2(v)\ts \Lc^{+[n-m]}_1(u)\ts \overline{R}^{\tss[n-m]}_{12}(u_{-}/v_{+}).
\end{gather*}
\eco

\bpr\label{prop:relmone}
Suppose that
$m+1\leqslant j,k,l\leqslant (m+1)'$ and $j\neq l'$. Then the following relations hold
in $U\big(\overline{R}^{\tss[n]}\big)$: if $j=l$ then
\begin{gather}
\e_{mj}^{\pm}(u)\ell^{\mp [n-m]}_{kl}(v)
 =\frac{qu_{\mp}-q^{-1}v_{\pm}}{u_{\mp}-v_{\pm}}\ell^{\mp [n-m]}_{kj}(v)\e_{ml}^{\pm}(u)
-\frac{\big(q-q^{-1}\big)u_{\mp}}{u_{\mp}-v_{\pm}}\ell^{\mp [n-m]}_{kj}(v)\e_{mj}^{\mp}(v),\nonumber\\
\e_{mj}^{\pm}(u)\ell^{\tss\pm [n-m]}_{kl}(v)
 =\frac{qu-q^{-1}v}{u-v}\ell^{\tss\pm [n-m]}_{kj}(v)\e_{ml}^{\pm}(u)
-\frac{\big(q-q^{-1}\big)u}{u-v}\ell^{\tss\pm [n-m]}_{kj}(v)\e_{mj}^{\pm}(v);\label{ELMPj=l}
\end{gather}
if $j<l$ then
\begin{gather}
\big[\e^{\pm}_{mj}(u),\ell^{\mp [n-m]}_{kl}(v)\big]
 =\frac{\big(q-q^{-1}\big)v_{\pm}}{u_{\mp}-v_{\pm}}\ell^{\mp [n-m]}_{kj}(v)\e^{\pm}_{ml}(u)-
\frac{\big(q-q^{-1}\big)u_{\mp}}{u_{\mp}-v_{\pm}}\ell^{\mp [n-m]}_{kj}(v)\e^{\mp}_{ml}(v),\nonumber\\
\big[\e^{\pm}_{mj}(u),\ell^{\tss\pm [n-m]}_{kl}(v)\big]
 =\frac{\big(q-q^{-1}\big)v}{u-v}\ell^{\tss\pm [n-m]}_{kj}(v)\e^{\pm}_{ml}(u)-
\frac{\big(q-q^{-1}\big)u}{u-v}\ell^{\tss\pm [n-m]}_{kj}(v)\e^{\pm}_{ml}(v);\label{ELMPj<l}
\end{gather}
if $j>l$ then
\begin{gather*}
\big[\e^{\pm}_{mj}(u),\ell^{\mp [n-m]}_{kl}(v)\big]
 =\frac{\big(q-q^{-1}\big)u_{\mp}}{u_{\mp}-v_{\pm}}\ell^{\mp [n-m]}_{kj}(v)
(\e^{\pm}_{ml}(u)-\e_{ml}^{\mp}(v)),\\
\big[\e^{\pm}_{mj}(u),\ell^{\tss\pm [n-m]}_{kl}(v)\big]
 =\frac{\big(q-q^{-1}\big)u}{u-v}\ell^{\tss\pm [n-m]}_{kj}(v)(\e^{\pm}_{ml}(u)-\e_{ml}^{\pm}(v)).
\end{gather*}
\epr

\bpr\label{prop:relmonf}
Suppose that $m+1\leqslant j,k,l\leqslant (m+1)'$ and $j\neq k'$. Then the following relations hold
in $U\big(\overline{R}^{\tss[n]}\big)$: if $j=k$ then
\begin{gather*}
\f_{jm}^{\pm}(u)\ell_{jl}^{\mp [n-m]}(v)
 =\frac{u_{\pm}-v_{\mp}}{qu_{\pm}-q^{-1}v_{\mp}}\ell_{jl}^{\mp [n-m]}(v)\f_{jm}^{\pm}(u)
+\frac{\big(q-q^{-1}\big)v_{\mp}}{qu_{\pm}-q^{-1}v_{\mp}}\ts\f_{jm}^{\mp}(v)\ell_{jl}^{\mp [n-m]}(v),\\
\f_{jm}^{\pm}(u)\ell_{jl}^{\pm [n-m]}(v)
 =\frac{u-v}{qu-q^{-1}v}\ell_{jl}^{\pm [n-m]}(uv)\ts\f_{jm}^{\pm}(u)
+\frac{\big(q-q^{-1}\big)v}{qu-q^{-1}v}\f_{jm}^{\pm}(v)\ell_{jl}^{\pm [n-m]}(v);
\end{gather*}
if $j<k$ then
\begin{gather*}
\big[\f_{jm}^{\pm}(u),\ell_{kl}^{\mp [n-m]}(v)\big]
 =\frac{\big(q-q^{-1}\big)v_{\mp}}{u_{\pm}-v_{\mp}}\f_{km}^{\mp}(v)\ell_{jl}^{\mp [n-m]}(v)
-\frac{\big(q-q^{-1}\big)u_{\pm}}{u_{\pm}-v_{\mp}}\f_{km}^{\pm}(u)\ell_{jl}^{\mp [n-m]}(v),\\
\big[\f_{jm}^{\pm}(u),\ell_{kl}^{\pm [n-m]}(v)\big]
 =\frac{\big(q-q^{-1}\big)v}{u-v}\f_{km}^{\pm}(v)\ell_{jl}^{\pm [n-m]}(v)
-\frac{\big(q-q^{-1}\big)u}{u-v}\f_{km}^{\pm}(u)\ell_{jl}^{\pm [n-m]}(v);
\end{gather*}
if $j>k$ then
\begin{gather*}
\big[\f_{jm}^{\pm}(u),\ell_{kl}^{\mp [n-m]}(v)\big]
 =\frac{\big(q-q^{-1}\big)v_{\mp}}{u_{\pm}-v_{\mp}}(\f_{km}^{\mp}(v)-\f_{km}^{\pm}(u))
\ell_{jl}^{\mp [n-m]}(v),\\
\big[\f_{jm}^{\pm}(u),\ell_{kl}^{\pm [n-m]}(v)\big]
 =\frac{\big(q-q^{-1}\big)v}{u-v}(\f_{km}^{\pm}(v)-\f_{km}^{\pm}(u))\ell_{jl}^{\pm [n-m]}(v).
\end{gather*}
\epr

\subsection[Type $A$ relations]{Type $\boldsymbol{A}$ relations}\label{subsec:typea}

Due to the observation made in Remark~\ref{rem:gln}
and the quasideterminant formulas \eqref{hmqua}, \eqref{eijmlqua}
and \eqref{fijlmqua}, some of the relations between the Gaussian generators will follow
from those for the quantum affine algebra
$U_q\big(\wh{\mathfrak{gl}}_n\big)$; see \cite{df:it}. To reproduce them, set
\begin{gather*}
\Lc^{A\pm}(u)=\sum_{i,j=1}^n \ell^{\tss\pm}_{ij}(u)\otimes e_{ij}
\end{gather*}
and consider the $R$-matrix used in \cite{df:it} which is given by
\begin{gather*}
R_A(u)=\sum_{i=1}^n e_{ii}\otimes e_{ii}+\frac{u-1}{qu-q^{-1}}
\sum_{i\neq j} e_{ii}\otimes e_{jj}\\
\hphantom{R_A(u)=}{} +\frac{q-q^{-1}}{qu-q^{-1}}\sum_{i>j}e_{ij}\otimes e_{ji}
+\frac{\big(q-q^{-1}\big)u}{qu-q^{-1}}\sum_{i<j}e_{ij}\otimes e_{ji}.
\end{gather*}
By comparing it with the $R$-matrix \eqref{rbar}, we come to the relations in the algebra
$U\big(\overline{R}^{\tss[n]}\big)$:
\begin{gather*}
R_A(u/w)\Lc^{A\pm}_{1}(u)\Lc^{A\pm}_2(v) =\Lc^{A\pm}_2(v)\Lc^{A\pm}_{1}(u){R}_A(u/v),\\
R_A(uq^c/v)\Lc^{A+}_1(u)\Lc^{A-}_2(v) =\Lc^{A-}_2(v)\Lc^{A+}_1(u){R}_A(uq^{-c}/v).
\end{gather*}
Hence we get the following relations for the Gaussian generators which were verified in~\cite{df:it}, where we use notation \eqref{Xibar}.

\bpr\label{TypeArelation}
In the algebra $U\big(\overline{R}^{\tss[n]}\big)$ we have
\begin{gather*}
\h^{\pm}_i(u)\h^{\pm}_j(v)=\h^{\pm}_j(v)\h^{\pm}_i(u),\qquad
\h^{\pm}_i(u)\h^{\mp}_i(v)=\h^{\mp}_i(v)\h^{\pm}_i(u) \qquad\text{for}\quad 1\leqslant i,j\leqslant n,
\\
\frac{u_{\pm}-v_{\mp}}{qu_{\pm}-q^{-1}v_{\mp}}\h^{\pm}_i(u)\h^{\mp}_j(v)=
\frac{u_{\mp}-v_{\pm}}{qu_{\mp}-q^{-1}v_{\pm}}\h^{\mp}_j(v)\h^{\pm}_i(u)
\qquad\text{for}\quad 1\leqslant i<j\leqslant n.
\end{gather*}
Moreover,
\begin{gather*}
\h_{i}^{\pm}(u)\Xc_{j}^{+}(v)
 =\frac{u_{\mp}-v}{q^{(\ep_i,\alpha_j)}u_{\mp}-q^{-(\ep_i,\alpha_j)}v}
\Xc_{j}^{+}(v)\h_{i}^{\pm}(u),\\
\h_{i}^{\pm}(u)\Xc_{j}^{-}(v)
 =\frac{q^{(\ep_i,\alpha_j)}u_{\pm}-q^{-(\ep_i,\alpha_j)}v}{u_{\pm}-v}
\Xc_{j}^{-}(v)\h_{i}^{\pm}(u) \qquad\text{for}
\quad 1\leqslant i\leqslant n,\quad 1\leqslant j< n,
\end{gather*}
while
\begin{gather*}
\big(u-q^{\pm (\alpha_i,\alpha_j)}v\big)\Xc_{i}^{\pm}\big(uq^{i}\big)\Xc_{j}^{\pm}(vq^{j})
=\big(q^{\pm (\alpha_i,\alpha_j)}u-v\big) \Xc_{j}^{\pm}(vq^{j})\Xc_{i}^{\pm}\big(uq^{i}\big),
\end{gather*}
and
\begin{gather*}
\big[\Xc_i^{+}(u),\Xc_j^{-}(v)\big]=\delta_{ij}\big(q-q^{-1}\big)\\
\hphantom{\big[\Xc_i^{+}(u),\Xc_j^{-}(v)\big]=}{}\times
\big(\delta\big(u\ts q^{-c}/v\big)\ts\h_i^{-}(v_+)^{-1}\h_{i+1}^{-}(v_+)
-\delta\big(u\ts q^{c}/v\big)\ts\h_i^{+}(u_+)^{-1}\h_{i+1}^{+}(u_+)\big)
\end{gather*}
for $1\leqslant i,j<n$, together with the Serre relations for the series $\Xc_1^{\pm}(u),\dots,\Xc_{n-1}^{\pm}(u)$.
\epr

\begin{Remark}\label{re:anotherA}
Consider the inverse matrices $\Lc^{\pm}(u)^{-1}=[\ell^{\tss\pm}_{ij}(u)']_{i,j=1}^N$.
By the defining relations~\eqref{gen rel1} and \eqref{gen rel2}, we have
\begin{gather*}
\Lc^{\pm}_{1}(u)^{-1}\Lc^{\pm}_2(v)^{-1}\overline{R}^{\tss[n]}(u/v)
 =\overline{R}^{\tss[n]}(u/v)\Lc^{\pm}_2(v)^{-1}\Lc^{\pm}_{1}(u)^{-1},\\
\Lc^{-}_2(v)^{-1}\Lc^{+}_1(u)^{-1}\overline{R}^{\tss[n]}\big({u}q^{c}/v\big)
 =\overline{R}^{\tss[n]}({u}q^{-c}/v)\Lc^{-}_2(v)^{-1}\Lc^{+}_1(u)^{-1}.
\end{gather*}
So we can get another family of generators of the algebra $U\big(\overline{R}^{\tss[n]}\big)$
which satisfy the defining relations of $U_q\big(\wh{\gl}_n\big)$. Namely, these relations are satisfied by
the coefficients of the series $\ell^{\tss\pm}_{ij}(u)'$ with $i,j=n',\dots,1'$. In particular,
by taking the inverse matrices, we get a Gauss decomposition for the matrix $[\ell^{\tss\pm}_{ij}(u)']_{i,j=n',\dots,1'}$ from the Gauss decomposition of the matrix $\Lc^{\pm}(u)$.
\end{Remark}

\subsection[Relations for low rank algebras: type $B$]{Relations for low rank algebras: type $\boldsymbol{B}$}\label{subsec:lowrankB}

In view of Theorem~\ref{thm:red}, a significant part of relations between
the Gaussian generators is implied by those in low rank algebras.
In this section we describe them for the algebra $U\big(\overline{R}^{\tss[1]}\big)$ in type~$B$
associated with the Lie algebra~$\oa_3$.

\ble\label{lem:rel_h1e1B} The following relations hold in the algebra $U\big(\overline{R}^{\tss[1]}\big)$.
For the diagonal generators we have
\begin{gather}
\h_1^{\pm}(u)\h_1^{\pm}(v)=\h_1^{\pm}(v)\h_1^{\pm}(u),
\qquad\h_1^{\pm}(u)\h_1^{\mp}(v)=\h_1^{\mp}(v)\h_1^{\pm}(u),\nonumber\\
\label{h1h2}
\h_1^{\pm}(u)\h_{2}^{\pm}(v)=\h_{2}^{\pm}(v)\h_1^{\pm}(u),\\
\label{h1h2pm}
\frac{u_{\pm}-v_{\mp}}{qu_{\pm}-q^{-1}v_{\mp}}\h_1^{\pm}(u)\h_{2}^{\mp}(v)=
\frac{u_{\mp}-v_{\pm}}{qu_{\mp}-q^{-1}v_{\pm}}\h_{2}^{\mp}(v)\h_1^{\pm}(u).
\end{gather}
Moreover,
\begin{gather}
\h_1^{\pm}(u)\tss\e_{1,2}^{\mp}(v) =
\frac{u_{\mp}-v_{\pm}}{qu_{\mp}-q^{-1}v_{\pm}}\ts\e_{1,2}^{\mp}(v)\h_1^{\pm}(u)
+\frac{\big(q-q^{-1}\big)v_{\pm}}{qu_{\mp}-q^{-1}v_{\pm}}\ts\h_1^{\pm}(u)\e_{1,2}^{\pm}(u),
\nonumber\\
\h_1^{\pm}(u)\tss\e_{1,2}^{\pm}(v) =
\frac{u-v}{qu-q^{-1}v}\ts\e_{1,2}^{\pm}(v)\h_1^{\pm}(u)
+\frac{\big(q-q^{-1}\big)v}{qu-q^{-1}v}\ts\h_1^{\pm}(u)\e_{1,2}^{\pm}(u),
\nonumber\\
\label{f21h1mp}
\f_{2,1}^{\mp}(v)\tss\h_1^{\pm}(u) =
\frac{u_{\pm}-v_{\mp}}{qu_{\pm}-q^{-1}v_{\mp}}\ts\h_1^{\pm}(u)\f_{2,1}^{\mp}(v)
+\frac{\big(q-q^{-1}\big)u_{\pm}}{qu_{\pm}-q^{-1}v_{\mp}}\ts\f_{2,1}^{\pm}(u)\h_1^{\pm}(u),
\\
\label{f21h1}
\f_{2,1}^{\pm}(v)\tss\h_1^{\pm}(u) =
\frac{u-v}{qu-q^{-1}v}\ts\h_1^{\pm}(u)\f_{2,1}^{\pm}(v)
+\frac{\big(q-q^{-1}\big)u}{qu-q^{-1}v}\ts\f_{2,1}^{\pm}(u)\h_1^{\pm}(u),
\end{gather}
and
\begin{gather}
\big[\e_{1,2}^{\pm}(u),\f_{2,1}^{\mp}(v)\big] =
\frac{\big(q-q^{-1}\big)u_{\mp}}{qu_{\mp}-q^{-1}v_{\pm}}\ts\h_{2}^{\mp}(v)\h_{1}^{\mp}(v)^{-1}
-\frac{\big(q-q^{-1}\big)u_{\pm}}{qu_{\pm}-q^{-1}v_{\mp}}\ts\h_{2}^{\pm}(u)\h_{1}^{\pm}(u)^{-1},
\nonumber\\
\label{e12f21}
\big[\e_{1,2}^{\pm}(u),\f_{2,1}^{\pm}(v)\big] =
\frac{\big(q-q^{-1}\big)u}{qu-q^{-1}v}\big(\h_{2}^{\pm}(v)\h_{1}^{\pm}(v)^{-1}
-\h_{2}^{\pm}(u)\h_{1}^{\pm}(u)^{-1}\big).
\end{gather}
\ele

\begin{proof} All relations in the lemma are consequences of those between the series $\ell_{ij}^{\pm}(u)$ and $\ell_{kl}^{\pm}(v)$ with $i\neq k'$ and $j\neq l'$ in the algebra $U\big(\overline{R}^{\tss[1]}\big)$.
Therefore, they are essentially relations occurring in type $A$ and verified in the
same way; cf.\ Proposition~\ref{TypeArelation}.
\end{proof}

Now we turn to the $B$-type-specific relations.

\ble\label{lem:e1e1B}
In the algebra $U\big(\overline{R}^{\tss[1]}\big)$ we have
\begin{gather*}
\e_{1,2}^{\pm}(u)\e_{1,2}^{\mp}(v)
=\frac{u_{\mp}-q^{-1}v_{\pm}}{q^{-1}u_{\mp}-v_{\pm}}\e_{1,2}^{\mp}(v)\e_{1,2}^{\pm}(u)
-\frac{\big(q-q^{-1}\big)u_{\mp}}{q^{-1}u_{\mp}-qv_{\pm}}\e_{1,2}^{\mp}(v)^2\\
\hphantom{\e_{1,2}^{\pm}(u)\e_{1,2}^{\mp}(v) =}{}
-\frac{\big(u_{\mp}-q^{-1}v_{\pm}\big)\big(1-q^{-2}\big)v_{\pm}}{\big(q^{-1}u_{\mp}-v_{\pm}\big)
\big(u_{\mp}-q^{-2}v_{\pm}\big)}\e_{1,2}^{\pm}(u)^2 \\
\hphantom{\e_{1,2}^{\pm}(u)\e_{1,2}^{\mp}(v) =}{}
 +\frac{(u_{\mp}-v_{\pm})q^{-1/2}\big(q^{-2}-1\big)v_{\pm}}{\big(q^{-1}u_{\mp}-v_{\pm}\big)
\big(u_{\mp}-q^{-2}v_{\pm}\big)}\e_{1,3}^{\pm}(u)\\
\hphantom{\e_{1,2}^{\pm}(u)\e_{1,2}^{\mp}(v) =}{}
+\frac{(u_{\mp}-v_{\pm})q^{-1/2}\big(q^{-1}-q\big)u_{\mp}}{\big(q^{-1}u_{\mp}-v_{\pm}\big)
\big(q^{-1}u_{\mp}-qv_{\pm}\big)}\e_{1,3}^{\mp}(v)
\end{gather*}
and
\begin{gather*}
\e_{1,2}^{\pm}(u)\e_{1,2}^{\pm}(v) =\frac{u-q^{-1}v}{q^{-1}u-v}\e_{1,2}^{\pm}(v)\e_{1,2}^{\pm}(u)
-\frac{\big(q-q^{-1}\big)u}{q^{-1}u-qv}\e_{1,2}^{\pm}(v)^2\\
\hphantom{\e_{1,2}^{\pm}(u)\e_{1,2}^{\pm}(v) =}{} -\frac{\big(u-q^{-1}v\big)\big(1-q^{-2}\big)v}{\big(q^{-1}u-v\big)\big(u-q^{-2}v\big)}\e_{1,2}^{\pm}(u)^2\\
\hphantom{\e_{1,2}^{\pm}(u)\e_{1,2}^{\pm}(v) =}{}
+\frac{(u-v)q^{-1/2}\big(q^{-2}-1\big)v}{\big(q^{-1}u-v\big)\big(u-q^{-2}v\big)}\e_{1,3}^{\pm}(u)
+\frac{(u-v)q^{-1/2}\big(q^{-1}-q\big)u}{\big(q^{-1}u-v\big)\big(q^{-1}u-qv\big)}\e_{1,3}^{\pm}(v).
\end{gather*}
Moreover,
\begin{gather*}
\f_{2,1}^{\mp}(v)\f_{2,1}^{\pm}(u) =\frac{u_{\pm}-q^{-1}
v_{\mp}}{q^{-1}u_{\pm}-v_{\mp}}\f_{2,1}^{\pm}(u)\f_{2,1}^{\mp}(v)
-\frac{\big(q-q^{-1}\big)v_{\mp}}{q^{-1}u_{\pm}-qv_{\mp}}\f_{2,1}^{\mp}(v)^2\\
\hphantom{\f_{2,1}^{\mp}(v)\f_{2,1}^{\pm}(u) =}{}
-\frac{\big(u_{\pm}-q^{-1}v_{\mp}\big)\big(1-q^{-2}\big)u_{\pm}}{(q^{-1}
u_{\pm}-v_{\mp})\big(u_{\pm}-q^{-2}v_{\mp}\big)}\f_{2,1}^{\pm}(u)^2\\
\hphantom{\f_{2,1}^{\mp}(v)\f_{2,1}^{\pm}(u) =}{}
+\frac{(u_{\pm}-v_{\mp})q^{-1/2}(q^{-2}-1)u_{\pm}}{\big(q^{-1}
u_{\pm}-v_{\mp}\big)\big(u_{\pm}-q^{-2}v_{\mp}\big)}\f_{3,1}^{\pm}(u)\\
\hphantom{\f_{2,1}^{\mp}(v)\f_{2,1}^{\pm}(u) =}{}
+\frac{(u_{\pm}-v_{\mp})q^{-1/2}\big(q^{-1}-q\big)v_{\mp}}{(q^{-1}
u_{\pm}-v_{\mp})\big(q^{-1}u_{\pm}-qv_{\mp}\big)}\f_{3,1}^{\mp}(v)
\end{gather*}
and
\begin{gather*}
\f_{2,1}^{\pm}(v)\f_{2,1}^{\pm}(u) =\frac{u-q^{-1}v}{q^{-1}u-v}\f_{2,1}^{\pm}(u)\f_{2,1}^{\pm}(v)
-\frac{\big(q-q^{-1}\big)v}{q^{-1}u-qv}\f_{2,1}^{\pm}(v)^2 \\
\hphantom{\f_{2,1}^{\pm}(v)\f_{2,1}^{\pm}(u) =}{}
 -\frac{\big(u-q^{-1}v\big)\big(1-q^{-2}\big)u}{\big(q^{-1}u-v\big)\big(u-q^{-2}v\big)}\f_{2,1}^{\pm}(u)^2\\
\hphantom{\f_{2,1}^{\pm}(v)\f_{2,1}^{\pm}(u) =}{}
+\frac{(u-v)q^{-1/2}\big(q^{-2}-1\big)u}{\big(q^{-1}u-v\big)\big(u-q^{-2}v\big)}\f_{3,1}^{\pm}(u)
+\frac{(u-v)q^{-1/2}\big(q^{-1}-q\big)v}{\big(q^{-1}u-v\big)\big(q^{-1}u-qv\big)}\f_{3,1}^{\pm}(v).
\end{gather*}
\ele

\begin{proof}
By using the expression \eqref{rbar2} for the $R$-matrix, we obtain from
\eqref{gen rel2} that
\begin{gather}
\ell_{12}^{+}(u)\ts\ell_{12}^{-}(v)=
\frac{1}{\big(y-q^{-2}\big)\big(y-q^{-1}\big)}\ts\big(a_{12}(y)\tss\ell_{11}^{-}(v)\ell_{13}^{+}(u)\nonumber\\
\hphantom{\ell_{12}^{+}(u)\ts\ell_{12}^{-}(v)=}{} +a_{22}(y)\tss\ell_{12}^{-}(v)\ell_{12}^{+}(u)+a_{32}(y)\tss\ell_{13}^{-}(v)\ell_{11}^{+}(u)\big),\label{eq:l12pm}
\end{gather}
where we set $y=u_-/u_+$. Similarly, we also have
\begin{gather}
\ell_{12}^{+}(u)\tss\ell_{11}^{-}(v) =\frac{y-1}{qy-q^{-1}}\ell_{11}^{-}(v)\ell_{12}^{+}(u)
+\frac{\big(q-q^{-1}\big)y}{qy-q^{-1}}\ell_{12}^{-}(v)\ell_{11}^{+}(u),
\nonumber\\
\label{l11l12pm}
\ell_{11}^{+}(u)\ell_{12}^{-}(v) =\frac{y-1}{qy-q^{-1}}\ell_{12}^{-}(v)\ell_{11}^{+}(u)
+\frac{q-q^{-1}}{qy-q^{-1}}\ell_{11}^{-}(v)\ell_{12}^{+}(u).
\end{gather}
In terms of Gaussian generators the left hand side of~\eqref{eq:l12pm} can now be written as
\begin{gather*}
\frac{y-1}{qy-q^{-1}}\ts\h_1^{-}(v)\h_1^{+}(u)\e_{12}^{+}(u)\e_{12}^{-}(v)
+\frac{\big(q-q^{-1}\big)y}{qy-q^{-1}}\ts\h_1^{-}(v)\e_{12}^{-}(v)\h_1^{+}(u)\e_{12}^{-}(v),
\end{gather*}
which equals
\begin{gather*}
\frac{q^{-1}y-q}{y-1}\h_1^{-}(v)\h_1^{+}(u)\e_{12}^{+}(u)\e_{12}^{-}(v)
+\frac{\big(q-q^{-1}\big)y}{y-1}\h_1^{+}(u)\h_1^{-}(v)\e_{12}^{-}(v)^2.
\end{gather*}
Now use another consequence of \eqref{gen rel2},
\begin{gather*}
\ell_{11}^{+}(u)\tss\ell_{13}^{-}(v)=
\frac{1}{\big(y-q^{-2}\big)\big(y-q^{-1}\big)}\big(a_{13}(y)\ell_{11}^{-}(v)\ell_{13}^{+}(u)\\
\hphantom{\ell_{11}^{+}(u)\tss\ell_{13}^{-}(v)=}{} +a_{23}(y)\ell_{12}^{-}(v)\ell_{12}^{+}(u)+a_{33}(y)\ell_{13}^{-}(v)\ell_{11}^{+}(u)\big),
\end{gather*}
which together with \eqref{l11l12pm} brings \eqref{eq:l12pm} to the form
\begin{gather*}
\frac{q^{-1}y-q}{y-1}\h_1^{-}(v)\h_1^{+}(u)\e_{12}^{+}(u)\e_{12}^{-}(v)
 +\frac{\big(q-q^{-1}\big)y}{y-1}\h_1^{+}(u)\h_1^{-}(v)\e_{12}^{-}(v)^2\\
\qquad{} =\frac{q^{-1/2}\big(q^{-2}-1\big)\big(q^{-1}y-q\big)}{\big(y-q^{-2}\big)\big(q^{-1}y-1\big)} \h_1^{-}(v)\h_1^{+}(u)\e_{13}^{+}(u)\\
\qquad\quad{} +\frac{\big(y-q^{-1}\big)\big(q^{-1}y-q\big)}{(y-1)\big(q^{-1}y-1\big)}\h_1^{+}(u)\h_1^{-}(v)\e_{12}^{-}(v)\e_{12}^{+}(u)\\
\qquad\quad{} -\frac{\big(y-q^{-1}\big)\big(q^{-2}y-1\big)\big(q-q^{-1}\big)}{\big(y-q^{-2}\big)\big(q^{-1}y-1\big) (y-1)}\h_1^{-}(v)\h_1^{+}(u)\e_{12}^{+}(u)^2\\
\qquad\quad{} +\frac{q^{1/2}\big(q^{-2}-1\big)y}{q^{-1}y-1}\h_1^{+}(u)\h_1^{-}(v)\e_{13}^{-}(v).
\end{gather*}
Since the series $\h_1^{+}(u)$ and $\h_1^{-}(u)$ are invertible and their coefficients pairwise commute,
we arrive at one case of the first relation of the lemma. The remaining relations
are verified by quite a~similar calculation.
\end{proof}

Now we will be concerned with relations in the algebra $U\big(\overline{R}^{\tss[1]}\big)$ involving the diagonal generators $\h^{\pm}_{2}(u)$.

\ble\label{lem:h2e12B} We have the relations
\begin{gather}
 \h^{\mp}_{2}(v)\f_{21}^{\pm}(u)
+\frac{\big(q-q^{-1}\big)v_{\mp}}{u_{\pm}-v_{\mp}}\f^{\mp}_{21}(v)\h^{\mp}_{2}(v)\nonumber\\
\qquad{} =\frac{\big(q^{-1}u_{\pm}-qv_{\mp}\big)\big(u_{\pm}-q^{-1}v_{\mp}\big)}{\big(u_{\pm}-v_{\mp}\big)
\big(q^{-1}u_{\pm}-v_{\mp}\big)}\f^{\pm}_{21}(u)\h^{\mp}_{2}(v)
+\frac{\big(q^{-2}-1\big)q^{1/2}v_{\mp}}{q^{-1}u_{\pm}-v_{\mp}}\f^{\mp}_{32}(v)\h^{\mp}_{2}(v)\label{f21h2B}
\end{gather}
and
\begin{gather*}
 \h^{\pm}_{2}(v)\f_{21}^{\pm}(u)
+\frac{\big(q-q^{-1}\big)v}{u-v}\f^{\pm}_{21}(v)\h^{\pm}_{2}(v)\\
 \qquad{} =\frac{\big(q^{-1}u-qv\big)\big(u-q^{-1}v\big)}{(u-v)\big(q^{-1}u-v\big)}\f^{\pm}_{21}(u)\h^{\pm}_{2}(v)
+\frac{\big(q^{-2}-1\big)q^{1/2}v}{q^{-1}u-v}\f^{\pm}_{32}(v)\h^{\pm}_{2}(v).
\end{gather*}
Moreover,
\begin{gather*}
 \e_{12}^{\pm}(u)\h^{\mp}_{2}(v)
+\frac{\big(q-q^{-1}\big)u_{\mp}}{u_{\mp}-v_{\pm}}\h^{\mp}_{2}(v)\e^{\mp}_{12}(v)\\
\qquad{}=\frac{\big(q^{-1}u_{\mp}-qv_{\pm}\big)\big(u_{\mp}-q^{-1}v_{\pm}\big)}{(u_{\mp}-v_{\pm})
\big(q^{-1}u_{\mp}-v_{\pm}\big)}\h^{\mp}_{2}(v)\e^{\pm}_{12}(u)
+\frac{\big(q^{-2}-1\big)q^{1/2}u_{\mp}}{q^{-1}u_{\mp}-v_{\pm}}\h^{\mp}_{2}(v)\e^{\mp}_{32}(v)
\end{gather*}
and
\begin{gather*}
 \e_{12}^{\pm}(u)\h^{\pm}_{2}(v)
+\frac{\big(q-q^{-1}\big)u}{u-v}\h^{\pm}_{2}(v)\e^{\pm}_{12}(v)\\
\qquad{} =\frac{\big(q^{-1}u-qv\big)\big(u-q^{-1}v\big)}{(u-v)\big(q^{-1}u-v\big)}\h^{\pm}_{2}(v)\e^{\pm}_{12}(u)
+\frac{(q^{-2}-1)q^{1/2}u}{q^{-1}u-v}\h^{\pm}_{2}(v)\e^{\pm}_{32}(v).
\end{gather*}
\ele

\begin{proof}All eight relations are verified in the same way so we only give full details
to check one case of \eqref{f21h2B}, where the top signs are chosen.
The defining relations \eqref{gen rel2} imply
\begin{gather}
\frac{1}{\big(x-q^{-2}\big)\big(x-q^{-1}\big)} \big(a_{21}(x)\ell_{31}^{+}(u)\ell_{12}^{-}(v)
+a_{22}(x)\ell_{21}^{+}(u)\ell_{22}^{-}(v)+a_{23}(x)\ell_{11}^{+}(u)\ell_{32}^{-}(v)\big)\nonumber\\
\qquad{} =\frac{y-1}{qy-q^{-1}}\ell^{-}_{22}(v)\ell^{+}_{21}(u)
+\frac{q-q^{-1}}{qy-q^{-1}}\ell^{-}_{21}(v)\ell^{+}_{22}(u),
\label{l21l22}
\end{gather}
where $x=u_{+}/v_{-}$ and $y=u_{-}/v_{+}$. In terms of the Gaussian generators the right hand side can be written as
\begin{gather*}
\frac{y-1}{qy-q^{-1}}\h^{-}_{2}(v)\ell^{+}_{21}(u)
+\frac{y-1}{qy-q^{-1}}\f^{-}_{21}(v)\ell^{-}_{12}(v)\ell^{+}_{21}(u)
+\frac{q-q^{-1}}{qy-q^{-1}}\f^{-}_{21}(v)\ell^{-}_{11}(v)\ell^{+}_{22}(u).
\end{gather*}
Applying \eqref{gen rel2} again we get the relations
\begin{gather*}
\frac{x-1}{qx-q^{-1}}\ell^{+}_{21}(u)\ell^{-}_{12}(v)
+\frac{q-q^{-1}}{qx-q^{-1}}\ell^{+}_{11}(u)\ell^{-}_{22}(v)\\
\qquad{} =\frac{y-1}{qy-q^{-1}}\ell^{-}_{12}(v)\ell^{+}_{21}(u)
+\frac{q-q^{-1}}{qy-q^{-1}}\ell^{-}_{11}(v)\ell^{+}_{22}(u)
\end{gather*}
and
\begin{gather*}
\frac{x-1}{qx-q^{-1}}\ell^{+}_{21}(u)\ell^{-}_{11}(v)
+\frac{q-q^{-1}}{qx-q^{-1}}\ell^{+}_{11}(u)\ell^{-}_{21}(v)=\ell^{-}_{11}(v)\ell^{+}_{21}(u).
\end{gather*}
They allow us to bring
the right hand side of \eqref{l21l22} to the form
\begin{gather*}
 \frac{y-1}{qy-q^{-1}}\h^{-}_{2}(v)\ell^{+}_{21}(u)
+\frac{x-1}{qx-q^{-1}}\f^{-}_{21}(v)\ell^{+}_{21}(u)\ell^{-}_{12}(v)
+\frac{q-q^{-1}}{qx-q^{-1}}\f^{-}_{21}(v)\ell^{+}_{11}(u)\ell^{-}_{22}(v)\\
 {} =\frac{y-1}{qy-q^{-1}}\h^{-}_{2}(v)\ell^{+}_{21}(u)+\f^{-}_{21}(v)
\left(\frac{x-1}{qx-q^{-1}}\ell^{+}_{21}(u)\ell^{-}_{11}(v)
+\frac{q-q^{-1}}{qx-q^{-1}}\ell^{+}_{11}(u)\ell^{-}_{21}(v)\right)\e^{-}_{12}(v)\\
\quad
{}+\frac{q-q^{-1}}{qx-q^{-1}}\f^{-}_{21}(v)\h^{+}_{1}(u)\h^{-}_{2}(v)
\end{gather*}
which is equal to
\begin{gather*}
\frac{y-1}{qy-q^{-1}}\h^{-}_{2}(v)\ell^{+}_{21}(u)+\ell^{-}_{21}(v)\ell^{+}_{21}(u)\e^{-}_{12}(v)
+\frac{q-q^{-1}}{qx-q^{-1}}\f^{-}_{21}(v)\h^{+}_{1}(u)\h^{-}_{2}(v).
\end{gather*}
Due to \eqref{gen rel2} the expression
\begin{gather*}
\frac{1}{\big(x-q^{-2}\big)\big(x-q^{-1}\big)}\big(a_{21}(x)\ell_{31}^{+}(u)\ell_{11}^{-}(v)
+a_{22}(x)\ell_{21}^{+}(u)\ell_{21}^{-}(v)+a_{23}(x)\ell_{11}^{+}(u)\ell_{31}^{-}(v)\big)
\end{gather*}
coincides with $\ell^{-}_{21}(v)\ell^{+}_{21}(u)$ so that the right hand side of~\eqref{l21l22} equals
\begin{gather*}
\frac{1}{\big(x-q^{-2}\big)\big(x-q^{-1}\big)} \big(a_{21}(x)\ell_{31}^{+}(u)\ell_{11}^{-}(v)
+a_{22}(x)\ell_{21}^{+}(u)\ell_{21}^{-}(v)+a_{23}(x)\ell_{11}^{+}(u)\ell_{31}^{-}(v)\big)\e^{-}_{12}(v)\\
\qquad{} +\frac{y-1}{qy-q^{-1}}\h^{-}_{2}(v)\ell^{+}_{21}(u)
+\frac{q-q^{-1}}{qx-q^{-1}}\f^{-}_{21}(v)\h^{+}_{1}(u)\h^{-}_{2}(v).
\end{gather*}
Hence we can write \eqref{l21l22} in the form
\begin{gather*}
\frac{1}{\big(x-q^{-2}\big)\big(x-q^{-1}\big)} \big(
a_{22}(x)\ell_{21}^{+}(u)\h_{2}^{-}(v)
+a_{23}(x)\ell_{11}^{+}(u)\f_{32}^{-}(v)\h_{2}^{-}(v)\big)\\
\qquad{} =\frac{y-1}{qy-q^{-1}}\h^{-}_{2}(v)\ell^{+}_{21}(u)
+\frac{q-q^{-1}}{qx-q^{-1}}\f^{-}_{21}(v)\h^{+}_{1}(u)\h^{-}_{2}(v).
\end{gather*}
Together with \eqref{f21h1mp} this leads to the relation
\begin{gather*}
\left(\frac{a_{22}(x)}{\big(x-q^{-2}\big)\big(x-q^{-1}\big)}
 -\frac{\big(q-q^{-1}\big)^2x}{\big(qx-q^{-1}\big)^2}\right)\ts\ell_{21}^{+}(u)\h_{2}^{-}(v) \\
 \qquad\quad{}
 +\frac{a_{23}(x)}{\big(x-q^{-2}\big)\big(x-q^{-1}\big)}\h_{1}^{+}(u)\f_{32}^{-}(v)\h_{2}^{-}(v)\\
\qquad {} =\frac{y-1}{qy-q^{-1}}\h^{-}_{2}(v)\ell^{+}_{21}(u)
+\frac{\big(q-q^{-1}\big)(x-1)}{\big(qx-q^{-1}\big)^2}\h^{+}_{1}(u)\f^{-}_{21}(v)\h^{-}_{2}(v).
\end{gather*}
By the following consequence of \eqref{f21h1},
\begin{gather*}
\ell_{21}^{+}(u)=\f_{21}^{+}(u)\h_{1}^{+}(u)=q\tss\h_{1}^{+}(u)\f_{21}^{+}\big(uq^2\big),
\end{gather*}
the relation takes the form
\begin{gather*}
 \left(\frac{a_{22}(x)}{\big(x-q^{-2}\big)\big(x-q^{-1}\big)}
-\frac{\big(q-q^{-1}\big)^2x}{\big(qx-q^{-1}\big)^2}\right)q\h_{1}^{+}(u)\f_{21}^{+}\big(uq^2\big)\h_{2}^{-}(v)\\
 \qquad\quad{} +\frac{a_{23}(x)}{\big(x-q^{-2}\big)\big(x-q^{-1}\big)}\h_{1}^{+}(u)\f_{32}^{-}(v)\h_{2}^{-}(v)\\
 \qquad{}=\frac{q(y-1)}{qy-q^{-1}}\h^{-}_{2}(v)\h_{1}^{+}(u)\f_{21}^{+}\big(uq^2\big)
+\frac{\big(q-q^{-1}\big)(x-1)}{\big(qx-q^{-1}\big)^2}\h^{+}_{1}(u)\f^{-}_{21}(v)\h^{-}_{2}(v).
\end{gather*}
Finally, apply relations \eqref{h1h2pm} between $\h^{+}_{1}(u)$ and $\h^{-}_{2}(v)$ and use
the invertibility of $\h^{+}_{1}(u)$ to come to the relation
\begin{gather}
\left(\frac{a_{22}(x)}{\big(x-q^{-2}\big)\big(x-q^{-1}\big)}
 -\frac{\big(q-q^{-1}\big)^2x}{\big(qx-q^{-1}\big)^2}\right)\ts q\tss\f_{21}^{+}\big(uq^2\big)\h_{2}^{-}(v)
 +\frac{a_{23}(x)}{\big(x-q^{-2}\big)\big(x-q^{-1}\big)}\f_{32}^{-}(v)\h_{2}^{-}(v)\nonumber\\
\qquad{} =\frac{q(x-1)}{qx-q^{-1}}\h^{-}_{2}(v)\f_{21}^{+}\big(uq^2\big)
+\frac{\big(q-q^{-1}\big)(x-1)}{\big(qx-q^{-1}\big)^2}\f^{-}_{21}(v)\h^{-}_{2}(v).\label{l21l22sixth}
\end{gather}
It remains to use the formulas for $a_{ij}(u)$ to see that \eqref{l21l22sixth} is equivalent to
the considered case of~\eqref{f21h2B}.
\end{proof}

\ble\label{lem:h2h2B}
In the algebra $U(\overline{R}^{\tss[1]})$ we have
\begin{gather*}
\frac{\big(q^{-1}u_{\pm}-qv_{\mp}\big)\big(u_{\pm}-q^{-1}v_{\mp}\big)}{\big(qu_{\pm}-q^{-1}v_{\mp}\big) \big(q^{-1}u_{\pm}-v_{\mp}\big)}\h_{2}^{\pm}(u)\h_{2}^{\mp}(v)
=\frac{\big(q^{-1}u_{\mp}-qv_{\pm}\big)\big(u_{\mp}-q^{-1}v_{\pm}\big)}{\big(qu_{\mp}-q^{-1}v_{\pm}\big) \big(q^{-1}u_{\mp}-v_{\pm}\big)} \h_{2}^{\mp}(v)\h_{2}^{\pm}(u)
\end{gather*}
and
\begin{gather*}
\h_{2}^{\pm}(u)\h_{2}^{\pm}(v) =\h_{2}^{\pm}(v)\h_{2}^{\pm}(u).
\end{gather*}
\ele

\begin{proof}We only give details for one case of the more complicated first relation by choosing the top signs; the remaining cases are considered in a similar way. We begin with the following
consequence of~\eqref{gen rel2},
\begin{gather}
 \frac{1}{\big(x-q^{-2}\big)\big(x-q^{-1}\big)}\big(a_{21}(x)\ell_{32}^{+}(u)\ell_{12}^{-}(v)
 +a_{22}(x)\ell_{22}^{+}(u)\ell_{22}^{-}(v)
 +a_{23}(x)\ell_{12}^{+}(u)\ell_{32}^{-}(v)\big)\label{l22l22}\\
\qquad{} =\frac{1}{\big(y-q^{-2}\big)\big(y-q^{-1}\big)}\big(a_{12}(y)\ell_{21}^{-}(v)\ell_{23}^{+}(u)
 +a_{22}(y)\ell_{22}^{-}(v)\ell_{22}^{+}(u)
 +a_{32}(y)\ell_{23}^{-}(v)\ell_{21}^{+}(u)\big),\nonumber
\end{gather}
where $x=u_{+}/v_{-}$ and $y=u_{-}/v_{+}$, and then express both sides in terms of the Gaussian generators.
The left hand side takes the form
\begin{gather*}
 \frac{1}{\big(x-q^{-2}\big)\big(x-q^{-1}\big)}\big(a_{21}(x)\ell_{32}^{+}(u)\ell_{11}^{-}(v)
 +a_{22}(x)\ell_{22}^{+}(u)\ell_{21}^{-}(v)
 +a_{23}(x)\ell_{12}^{+}(u)\ell_{31}^{-}(v)\big)\e_{12}^{-}(v)\\
\qquad{} +\frac{1}{\big(x-q^{-2}\big)\big(x-q^{-1}\big)}\big(
 a_{22}(x)\ell_{22}^{+}(u)\h_{2}^{-}(v)
 +a_{23}(x)\ell_{12}^{+}(u)\f_{32}^{-}(v)\h_{2}^{-}(v)\big).
\end{gather*}
The defining relations \eqref{gen rel2} also give
\begin{gather*}
 \frac{1}{\big(x-q^{-2}\big)\big(x-q^{-1}\big)}\big(a_{21}(x)\ell_{32}^{+}(u)\ell_{11}^{-}(v)
 +a_{22}(x)\ell_{22}^{+}(u)\ell_{21}^{-}(v)
 +a_{23}(x)\ell_{12}^{+}(u)\ell_{31}^{-}(v)\big)\\
\qquad{} =\frac{y-1}{qy-q^{-1} }\ell_{21}^{-}(v)\ell_{22}^{+}(u)
+\frac{\big(q-q^{-1}\big)y}{qy-q^{-1}}\ell_{22}^{-}(v)\ell_{21}^{+}(u)
\end{gather*}
so that the left hand side of \eqref{l22l22} takes the form
\begin{gather*}
 \frac{y-1}{qy-q^{-1}}\ell_{21}^{-}(v)\ell_{22}^{+}(u)\e_{12}^{-}(v)
+\frac{\big(q-q^{-1}\big)y}{qy-q^{-1}}\ell_{22}^{-}(v)\ell_{21}^{+}(u)\e_{12}^{-}(v)\\
 \qquad{} +\frac{1}{\big(x-q^{-2}\big)\big(x-q^{-1}\big)}\big(
 a_{22}(x)\ell_{22}^{+}(u)\h_{2}^{-}(v)
 +a_{23}(x)\ell_{12}^{+}(u)\f_{32}^{-}(v)\h_{2}^{-}(v)\big).
\end{gather*}
A similar calculation shows that the right hand side of \eqref{l22l22} equals
\begin{gather*}
 \frac{x-1}{qx-q^{-1}}\f_{21}^{-}(v)\ell_{22}^{+}(u)\ell_{12}^{-}(v)
+\frac{q-q^{-1}}{qx-q^{-1}}\f_{21}^{-}(v)\ell_{12}^{+}(u)\ell_{22}^{-}(v)\\
\qquad{} +\frac{1}{\big(y-q^{-2}\big)\big(y-q^{-1}\big)}\big(
 a_{22}(y)\h_{2}^{-}(v)\ell_{22}^{+}(u)
 +a_{23}(x)\h_{2}^{-}(v)\e_{23}^{-}(v)\ell_{21}^{+}(u)\big).
\end{gather*}
Therefore, by rearranging \eqref{l22l22} we come to the relation
\begin{gather}
 \frac{1}{\big(x-q^{-2}\big)\big(x-q^{-1}\big)}\big(a_{22}(x)\ell_{22}^{+}(u)\h_{2}^{-}(v)
 +a_{23}(x)\ell_{12}^{+}(u)\f_{32}^{-}(v)\h_{2}^{-}(v)\big)\nonumber\\
 \qquad\quad{}-\frac{x-1}{qx-q^{-1}}\f_{21}^{-}(v)\ell_{22}^{+}(u)\ell_{12}^{-}(v)
-\frac{q-q^{-1}}{qx-q^{-1}}\f_{21}^{-}(v)\ell_{12}^{+}(u)\ell_{22}^{-}(v)\nonumber\\
 \qquad{}=\frac{1}{\big(y-q^{-2}\big)\big(y-q^{-1}\big)}\big(a_{22}(y)\h_{2}^{-}(v)\ell_{22}^{+}(u)
 +a_{23}(y)\h_{2}^{-}(v)\e_{23}^{-}(v)\ell_{21}^{+}(u)\big)\nonumber\\
 \qquad\quad{}-\frac{y-1}{qy-q^{-1}}\ell_{21}^{-}(v)\ell_{22}^{+}(u)\e_{12}^{-}(v)
-\frac{\big(q-q^{-1}\big)y}{qy-q^{-1}}\ell_{22}^{-}(v)\ell_{21}^{+}(u)\e_{12}^{-}(v).\label{l22l22second}
\end{gather}
Furthermore, by \eqref{gen rel2} we have
\begin{gather*}
\frac{x-1}{qx-q^{-1}}\ell_{22}^{+}(u)\ell_{11}^{-}(v)
+\frac{q-q^{-1}}{qx-q^{-1}}\ell_{12}^{+}(u)\ell_{21}^{-}(v)\\
\qquad{} =\frac{y-1}{qy-q^{-1}}\ell_{11}^{-}(v)\ell_{22}^{+}(u)
+\frac{\big(q-q^{-1}\big)y}{qy-q^{-1}}\ell_{12}^{-}(v)\ell_{21}^{+}(u),
\end{gather*}
which allows us to write \eqref{l22l22second} in the form
\begin{gather}
 \frac{1}{\big(x-q^{-2}\big)\big(x-q^{-1}\big)}\big(a_{22}(x)\ell_{22}^{+}(u)\h_{2}^{-}(v)
 +a_{23}(x)\ell_{12}^{+}(u)\f_{32}^{-}(v)\h_{2}^{-}(v)\big)\nonumber\\
\qquad\quad{}
-\frac{q-q^{-1}}{qx-q^{-1}}\f_{21}^{-}(v)\ell_{12}^{+}(u)\h_{2}^{-}(v)\nonumber\\
\qquad{} =\frac{1}{\big(y-q^{-2}\big)\big(y-q^{-1}\big)}\big(a_{22}(y)\h_{2}^{-}(v)\ell_{22}^{+}(u)
 +a_{23}(y)\h_{2}^{-}(v)\e_{23}^{-}(v)\ell_{21}^{+}(u)\big)\nonumber\\
\qquad\quad{}
-\frac{\big(q-q^{-1}\big)y}{qy-q^{-1}}\h_{2}^{-}(v)\ell_{21}^{+}(u)\e_{12}^{-}(v).\label{l22l22third}
\end{gather}
Now transform the left hand side of this relation.
Since
\begin{gather*}
\ell_{12}^{+}(u)\ell_{11}^{-}(v)=\frac{y-1}{qy-q^{-1}}\ell_{11}^{-}(v)\ell_{12}^{+}(u)
+\frac{\big(q-q^{-1}\big)y}{qy-q^{-1}}\ell_{12}^{-}(v)\ell_{11}^{+}(u),
\end{gather*}
we have
\begin{gather*}
\f_{21}^{-}(v)\ell_{12}^{+}(u)\h_{2}^{-}(v)=\ell_{21}^{-}(v)
\ell_{11}^{-}(v)^{-1}\ell_{12}^{+}(u)\h_{2}^{-}(v)\\
\qquad{}=\frac{y-1}{qy-q^{-1}}\ell_{21}^{-}(v)\ell_{12}^{+}(u)
\ell_{11}^{-}(v)^{-1}\h_{2}^{-}(v)+\frac{\big(q-q^{-1}\big)y}{qy-q^{-1}}\ell_{21}^{-}(v)
\e_{12}^{-}(v)\ell_{11}^{+}(u)\ell_{11}^{-}(v)^{-1}\h_{2}^{-}(v),
\end{gather*}
which equals
\begin{gather*}
\frac{y-1}{qy-q^{-1}}\ell_{21}^{-}(v)\ell_{12}^{+}(u)\ell_{11}^{-}(v)^{-1}\h_{2}^{-}(v)
+\frac{\big(q-q^{-1}\big)y}{qy-q^{-1}}\ell_{22}^{-}(v)\ell_{11}^{+}(u)\ell_{11}^{-}(v)^{-1}\h_{2}^{-}(v)\\
\qquad{} -\frac{\big(q-q^{-1}\big)y}{qy-q^{-1}}\h_{2}^{-}(v)\ell_{11}^{+}(u)\ell_{11}^{-}(v)^{-1}\h_{2}^{-}(v).
\end{gather*}
Furthermore, by \eqref{gen rel2} we have
\begin{gather*}
\frac{x-1}{qx-q^{-1}}\ell_{12}^{+}(u)\ell_{21}^{-}(v)
+\frac{\big(q-q^{-1}\big)x}{qx-q^{-1}}\ell_{22}^{+}(u)\ell_{11}^{-}(v)\\
\qquad{} =
\frac{y-1}{qy-q^{-1}}\ell_{21}^{-}(v)\ell_{12}^{+}(u)
+\frac{\big(q-q^{-1}\big)y}{qy-q^{-1}}\ell_{22}^{-}(v)\ell_{11}^{+}(u)
\end{gather*}
and so
\begin{gather*}
\f_{21}^{-}(v)\ell_{12}^{+}(u)\h_{2}^{-}(v)
=\frac{x-1}{qx-q^{-1}}\ell_{12}^{+}(u)\f_{21}^{-}(v)\h_{2}^{-}(v)
+\frac{\big(q-q^{-1}\big)x}{qx-q^{-1}}\ell_{22}^{+}(u)\h_{2}^{-}(v)\\
\hphantom{\f_{21}^{-}(v)\ell_{12}^{+}(u)\h_{2}^{-}(v)=}{} -\frac{\big(q-q^{-1}\big)y}{qy-q^{-1}}\h_{2}^{-}(v)\h_{1}^{+}(u)\h_{1}^{-}(v)^{-1}\h_{2}^{-}(v).
\end{gather*}
Therefore, the left hand side of \eqref{l22l22third} is equal to
\begin{gather*}
\frac{1}{\big(x-q^{-2}\big)\big(x-q^{-1}\big)}\big(a_{22}(x)\ell_{22}^{+}(u)\h_{2}^{-}(v)
 +a_{23}(x)\ell_{12}^{+}(u)\f_{32}^{-}(v)\h_{2}^{-}(v)\big)\\
\qquad{} -\frac{\big(q-q^{-1}\big)(x-1)}{\big(qx-q^{-1}\big)^2}\ell_{12}^{+}(u)\f_{21}^{-}(v)\h_{2}^{-}(v)
-\frac{\big(q-q^{-1}\big)^2x}{\big(qx-q^{-1}\big)^2}\ell_{22}^{+}(u)\h_{2}^{-}(v)\\
\qquad{} +\frac{\big(q-q^{-1}\big)^2y}{\big(qx-q^{-1}\big)\big(qy-q^{-1}\big)} \h_{2}^{-}(v)\h_{1}^{+}(u)\h_{1}^{-}(v)^{-1}\h_{2}^{-}(v).
\end{gather*}
Similarly, the right hand side of \eqref{l22l22third} takes the form
\begin{gather*}
\frac{1}{\big(y-q^{-2}\big)\big(y-q^{-1}\big)}\big(a_{22}(y)\h_{2}^{-}(v)\ell_{22}^{+}(u)
 +a_{23}(y)\h_{2}^{-}(v)\e_{23}^{-}(v)\ell_{21}^{+}(u)\big)\\
\qquad{} -\frac{\big(q-q^{-1}\big)y(y-1)}{\big(qy-q^{-1}\big)^2}\h_{2}^{-}(v)\e_{12}^{-}(v)\ell_{21}^{+}(u)
-\frac{\big(q-q^{-1}\big)^2y}{\big(qy-q^{-1}\big)^2}\h_{2}^{-}(v)\ell_{22}^{+}(u)\\
\qquad{} +\frac{\big(q-q^{-1}\big)^2y}{\big(qx-q^{-1}\big)\big(qy-q^{-1}\big)} \h_{2}^{-}(v)\h_{1}^{-}(v)^{-1}\h_{1}^{+}(u)\h_{2}^{-}(v),
\end{gather*}
and taking into account the relation $\h_{2}^{-}(v)\h_{1}^{-}(v)=\h_{1}^{-}(v)\h_{2}^{-}(v)$,
we get
\begin{gather}
 \frac{1}{\big(x-q^{-2}\big)\big(x-q^{-1}\big)}\big(a_{22}(x)\ell_{22}^{+}(u)\h_{2}^{-}(v)
 +a_{23}(x)\ell_{12}^{+}(u)\f_{32}^{-}(v)\h_{2}^{-}(v)\big)\nonumber\\
\qquad\quad{}-\frac{\big(q-q^{-1}\big)(x-1)}{\big(qx-q^{-1}\big)^2}\ell_{12}^{+}(u)\f_{21}^{-}(v)\h_{2}^{-}(v)
-\frac{\big(q-q^{-1}\big)^2x}{\big(qx-q^{-1}\big)^2}\ell_{22}^{+}(u)\h_{2}^{-}(v)\nonumber\\
 \qquad{}=\frac{1}{\big(y-q^{-2}\big)\big(y-q^{-1}\big)}\big(a_{22}(y)\h_{2}^{-}(v)\ell_{22}^{+}(u)
 +a_{23}(y)\h_{2}^{-}(v)\e_{23}^{-}(v)\ell_{21}^{+}(u)\big)\nonumber\\
 \qquad\quad {}-\frac{\big(q-q^{-1}\big)y(y-1)}{\big(qy-q^{-1}\big)^2}\h_{2}^{-}(v)\e_{12}^{-}(v)\ell_{21}^{+}(u)
-\frac{\big(q-q^{-1}\big)^2y}{\big(qy-q^{-1}\big)^2}\h_{2}^{-}(v)\ell_{22}^{+}(u).\label{l22l22fourth}
\end{gather}
It follows from \eqref{f21h1} that $\f_{21}^{+}(u)\h_1^{+}(u)=q\ts\h_1^{+}(u)\f_{21}^{+}\big(uq^2\big)$ which implies
\begin{gather*}
\f_{21}^{+}(u)\h_1^{+}(u)\e_{12}^{+}(u)=q\h_1^{+}(u)\f_{21}^{+}\big(uq^2\big)\e_{12}^{+}(u)
=q\h_1^{+}(u)\e_{12}^{+}(u)\f_{21}^{+}\big(uq^2\big)\!
-q\h_1^{+}(u)\tss\e_{12}^{+}(u)\f_{21}^{+}\big(uq^2\big).
\end{gather*}
Furthermore, \eqref{e12f21} gives
\begin{gather*}
\f_{21}^{+}(u)\h_1^{+}(u)\e_{12}^{+}(u)
=q\h_1^{+}(u)\e_{12}^{+}(u)\f_{21}^{+}\big(uq^2\big)
+\h_1^{+}(u)\h_{1}^{+}\big(uq^2\big)^{-1}\h_{2}^{+}\big(uq^2\big)
-\h_{2}^{+}(u)
\end{gather*}
and so
\begin{gather*}
\ell_{22}^{+}(u)=q\h_1^{+}(u)\e_{12}^{+}(u)\f_{21}^{+}\big(uq^2\big)
+\h_1^{+}(u)\h_{1}^{+}\big(uq^2\big)^{-1}\h_{2}^{+}\big(uq^2\big).
\end{gather*}
Therefore, the left hand side of \eqref{l22l22fourth} takes the form
\begin{gather*}
 \left(\frac{a_{22}(x)}{\big(x-q^{-2}\big)\big(x-q^{-1}\big)}
-\frac{\big(q-q^{-1}\big)^2x}{\big(qx-q^{-1}\big)^2}\right)q\ell_{12}^{+}(u)\f_{21}^{+}\big(uq^2\big)\h_{2}^{-}(v)\\
\qquad{}  +\frac{a_{23}(x)}{\big(x-q^{-2}\big)\big(x-q^{-1}\big)}\ell_{12}^{+}(u)\f_{32}^{-}(v)\h_{2}^{-}(v)
-\frac{\big(q-q^{-1}\big)(x-1)}{\big(qx-q^{-1}\big)^2}\ell^{+}_{12}(u)\f^{-}_{21}(v)\h^{-}_{2}(v)\\
\qquad{}+ \left(\frac{a_{22}(x)}{\big(x-q^{-2}\big)\big(x-q^{-1}\big)}
+\frac{\big(q-q^{-1}\big)^2x}{\big(qx-q^{-1}\big)^2}\right) \h_{1}^{+}(u)\h_{1}^{+}\big(uq^2\big)^{-1}\h_{2}^{+}\big(uq^2\big)\h_{2}^{-}(v).
\end{gather*}
Finally, by using the \eqref{l21l22sixth} we can write the left hand side of \eqref{l22l22fourth} as
\begin{gather*}
 \left(\frac{a_{22}(x)}{\big(x-q^{-2}\big)\big(x-q^{-1}\big)}
+\frac{\big(q-q^{-1}\big)^2x}{\big(qx-q^{-1}\big)^2}\right) \h_{1}^{+}(u)\h_{1}^{+}\big(uq^2\big)^{-1}\h_{2}^{+}\big(uq^2\big)\h_{2}^{-}(v)\\
\qquad{} +\frac{q(x-1)}{qx-q^{-1}}\ell_{12}^{+}(u)\h_2^{-}(v)\f_{21}^{+}\big(uq^2\big).
\end{gather*}
Similarly, the right hand side of \eqref{l22l22fourth} equals
\begin{gather*}
 \left(\frac{a_{22}(y)}{\big(y-q^{-2}\big)\big(y-q^{-1}\big)}
+\frac{\big(q-q^{-1}\big)^2y}{\big(qy-q^{-1}\big)^2}\right)\h_{2}^{-}(v)\h_{2}^{+}\big(uq^2\big)\h_{1}^{+} \big(uq^2\big)^{-1}\h_{1}^{+}(u)\\
\qquad{} +\frac{q(x-1)}{qx-q^{-1}}\ell_{12}^{+}(u)\h_2^{-}(v)\f_{21}^{+}\big(uq^2\big).
\end{gather*}
Cancelling equal terms on both sides and applying~\eqref{h1h2} and~\eqref{h1h2pm} we get
\begin{gather*}
 \left(\!\frac{a_{22}(x)}{\big(x-q^{-2}\big)\big(x-q^{-1}\big)}
+\frac{\big(q-q^{-1}\big)^2x}{\big(qx-q^{-1}\big)^2}\right)\!\frac{q^2x-1}{q^3x-q^{-1}}\frac{q^3y-q^{-1}}{q^2y-1}
\h_{1}^{+}(u)\h_{2}^{+}\big(uq^2\big)\h_{2}^{-}(v)\h_{1}^{+}\big(uq^2\big)^{-1}\\
 =\left(\!\frac{a_{22}(y)}{\big(y-q^{-2}\big)\big(y-q^{-1}\big)}
+\frac{\big(q-q^{-1}\big)^2y}{\big(qy-q^{-1}\big)^2}\right)\!\frac{x-1}{qx-q^{-1}}\frac{qy-q^{-1}}{y-1}
\h_{1}^{+}(u)\h_{2}^{-}(v)\h_{2}^{+}\big(uq^2\big)\h_{1}^{+}\big(uq^2\big)^{-1}.
\end{gather*}
Recalling the formula for $a_{22}(u)$ and using the invertibility of $\h_{1}^{+}(u)$, we come to the
relation
\begin{gather*}
\frac{(x-1)\big(qx-q^{-2}\big)}{\big(x-q^{-1}\big)\big(q^2x-q^{-2}\big)}\h_{2}^{+}\big(uq^2\big)\h_{2}^{-}(v)
=\frac{(y-1)\big(qy-q^{-2}\big)}{\big(y-q^{-1}\big)\big(q^2y-q^{-2}\big)}\h_{2}^{-}(v)\h_{2}^{+}\big(uq^2\big),
\end{gather*}
which is equivalent to the considered case of the first relation in the lemma.
\end{proof}

\subsection[Relations for low rank algebras: type $D$]{Relations for low rank algebras: type $\boldsymbol{D}$}\label{subsec:lowrankD}

As with the case of type $B$, a key role in deriving relations in $U\big(\overline{R}^{\tss[n]}\big)$ between the Gaussian generators will be played by Theorem~\ref{thm:red} and Proposition~\ref{prop:gauss-consist}. This time we will need relations in the algebra $U\big(\overline{R}^{\tss[2]}\big)$ in type~$D$ associated with the semisimple Lie algebra~$\oa_4$.

\ble\label{lem:n-1throotD}
The following relations hold in the algebra $U\big(\overline{R}^{\tss[2]}\big)$. For the diagonal generators we have
\begin{gather*}
\h_{i}^{\pm}(u)\h_{i}^{\pm}(v)=\h_{i}^{\pm}(v)\h_{i}^{\pm}(u),\qquad
\h_{i}^{\pm}(u)\h_{i}^{\mp}(v)=\h_{i}^{\mp}(v)\h_{i}^{\pm}(u),\qquad i=1,2,
\\
\h_{1}^{\pm}(u)\h_{2}^{\pm}(v)=\h_{2}^{\pm}(v)\h_{1}^{\pm}(u),
\qquad
\frac{u_{\pm}-v_{\mp}}{qu_{\pm}-q^{-1}v_{\mp}}\h_{1}^{\pm}(u)\h_{2}^{\mp}(v)
=\frac{u_{\mp}-v_{\pm}}{qu_{\mp}-q^{-1}v_{\pm}}\h_{2}^{\mp}(v)\h_{1}^{\pm}(u).
\end{gather*}
Moreover,
\begin{gather*}
\h_{1}^{\pm}(u)\e_{12}^{\mp}(v) =\frac{u_{\mp}-v_{\pm}}
{qu_{\mp}-q^{-1}v_{\pm}}\e_{12}^{\mp}(v)\h_{1}^{\pm}(u)+
\frac{\big(q-q^{-1}\big)v_{\pm}}{qu_{\mp}-q^{-1}v_{\pm}}\h_{1}^{\pm}(u)\e_{12}^{\pm}(u),\\
\h_{1}^{\pm}(u)\e_{12}^{\pm}(v) =\frac{u-v}{qu-q^{-1}v}\e_{12}^{\pm}(v)\h_{1}^{\pm}(u)+
\frac{\big(q-q^{-1}\big)v}{qu-q^{-1}v}\h_{1}^{\pm}(u)\e_{12}^{\pm}(u),\\
\f_{21}^{\mp}(v)\h_{1}^{\pm}(u) =\frac{u_{\pm}-v_{\mp}}
{qu_{\pm}-q^{-1}v_{\mp}}\h_{1}^{\pm}(u)\f_{21}^{\mp}(v)+
\frac{\big(q-q^{-1}\big)u_{\pm}}{qu_{\pm}-q^{-1}v_{\mp}}\f_{21}^{\pm}(u)\h_{1}^{\pm}(u),\\
\f_{21}^{\pm}(v)\h_{1}^{\pm}(u) =\frac{u-v}{qu-q^{-1}v}\h_{1}^{\pm}(u)\f_{21}^{\pm}(v)+
\frac{\big(q-q^{-1}\big)u}{qu-q^{-1}v}\f_{21}^{\pm}(u)\h_{1}^{\pm}(u),
\end{gather*}
and
\begin{gather*}
\e_{12}^{\pm}(u)\h_{2}^{\mp}(v) =
\frac{qu_{\mp}-q^{-1}v_{\pm}}{u_{\mp}-v_{\pm}}\h_{2}^{\mp}(v)\e_{12}^{\pm}(u)
-\frac{\big(q-q^{-1}\big)u_{\mp}}{u_{\mp}-v_{\pm}}\h_{2}^{\mp}(v)\e_{12}^{\mp}(v),\\
\e_{12}^{\pm}(u)\h_{2}^{\pm}(v) =
\frac{qu-q^{-1}v}{u-v}\h_{2}^{\pm}(v)\e_{12}^{\pm}(u)
-\frac{\big(q-q^{-1}\big)u}{u-v}\h_{2}^{\pm}(v)\e_{12}^{\pm}(v),\\
\h_{2}^{\mp}(v)\f_{21}^{\pm}(u) =
\frac{qu_{\pm}-q^{-1}v_{\mp}}{u_{\pm}-v_{\mp}}\f_{21}^{\pm}(u)\h_{2}^{\mp}(v)
-\frac{\big(q-q^{-1}\big)v_{\mp}}{u_{\pm}-v_{\mp}}\f_{21}^{\mp}(v)\h_{2}^{\mp}(v),\\
\h_{2}^{\pm}(v)\f_{21}^{\pm}(u) =
\frac{qu-q^{-1}v}{u-v}\f_{21}^{\pm}(u)\h_{2}^{\pm}(v)
-\frac{\big(q-q^{-1}\big)v}{u-v}\f_{21}^{\pm}(v)\h_{2}^{\pm}(v).
\end{gather*}
For the off-diagonal generators we have
\begin{gather*}
\e_{12}^{\pm}(u)\e_{12}^{\mp}(v) =-\frac{\big(q-q^{-1}\big)u_{\mp}}{q^{-1}u_{\mp}-qv_{\pm}}\e_{12}^{\mp}(v)^2
-\frac{\big(q-q^{-1}\big)v_{\pm}}{q^{-1}u_{\mp}-qv_{\pm}}\e_{12}^{\pm}(u)^2
+\frac{qu_{\mp}-q^{-1}v_{\pm}}{q^{-1}u_{\mp}-qv_{\pm}}\e_{12}^{\mp}(v)\e_{12}^{\pm}(u),\\
\e_{12}^{\pm}(u)\e_{12}^{\pm}(v) =-\frac{\big(q-q^{-1}\big)u}{q^{-1}u-qv}\e_{12}^{\pm}(v)^2
-\frac{\big(q-q^{-1}\big)v}{q^{-1}u-qv}\e_{12}^{\pm}(u)^2
+\frac{qu-q^{-1}v}{q^{-1}u-qv}\e_{12}^{\pm}(v)\e_{12}^{\pm}(u),\\
\f_{21}^{\pm}(u)\f_{21}^{\mp}(v) =\frac{\big(q-q^{-1}\big)u_{\pm}}{qu_{\pm}-q^{-1}v_{\mp}}\f_{21}^{\pm}(u)^2
+\frac{\big(q-q^{-1}\big)v_{\mp}}{qu_{\pm}-q^{-1}v_{\mp}}\f_{21}^{\mp}(v)^2
+\frac{q^{-1}u_{\pm}-qv_{\mp}}{qu_{\pm}-q^{-1}v_{\mp}}\f_{21}^{\mp}(v)\f_{21}^{\pm}(u),\\
\f_{21}^{\pm}(u)\f_{21}^{\pm}(v) =\frac{\big(q-q^{-1}\big)u}{qu-q^{-1}v}\f_{21}^{\pm}(u)^2
+\frac{\big(q-q^{-1}\big)v}{qu-q^{-1}v}\f_{21}^{\pm}(v)^2
+\frac{q^{-1}u-qv}{qu-q^{-1}v}\f_{21}^{\pm}(v)\f_{21}^{\pm}(u),
\end{gather*}
together with
\begin{gather*}
\big[\e_{12}^{\pm}(u),\f_{21}^{\mp}(v)\big] =
\frac{\big(q-q^{-1}\big)u_{\mp}}{u_{\mp}-v_{\pm}}\h_{2}^{\mp}(v)\h_{1}^{\mp}(v)^{-1}
-\frac{\big(q-q^{-1}\big)u_{\pm}}{u_{\pm}-v_{\mp}}\h_{2}^{\pm}(u)\h_{1}^{\pm}(u)^{-1},\\
\big[\e_{12}^{\pm}(u),\f_{21}^{\pm}(v)\big] =
\frac{\big(q-q^{-1}\big)u}{u-v}\big(\h_{2}^{\pm}(v)\h_{1}^{\pm}(v)^{-1}
-\h_{2}^{\pm}(u)\h_{1}^{\pm}(u)^{-1}\big).
\end{gather*}
\ele

\begin{proof} The generating series $\ell_{ij}^{\pm}(u)$ with $i,j=1,2$ satisfy the same relations as those in the algebra $U_q(\widehat{\mathfrak{gl}}_2)$; cf.\ Section~\ref{subsec:typea}. Therefore, all relations follow by the same calculations as in~\cite{df:it}.
\end{proof}

\ble\label{lem:e23=0D} In the algebra $U\big(\overline{R}^{\tss[2]}\big)$ we have
\begin{gather*}
\e_{23}^{\pm}(u)=\f_{32}^{\pm}(u)=0.
\end{gather*}
\ele

\begin{proof}
By Corollary~\ref{cor:guass-embed},
\begin{gather*}
\ell_{22}^{\pm\ts [1]}(u)\ell_{23}^{\pm\ts [1]}(v)=
\frac{\big(q^{-1}u-qv\big)(u-v)}{\big(qu-q^{-1}v\big)\big(u-q^{-1}v\big)}\ell_{23}^{\pm\ts [1]}(v)\ell_{22}^{\pm\ts [1]}(u).
\end{gather*}
Hence $\ell_{22}^{\pm\ts [1]}(u)\ell_{23}^{\pm\ts [1]}(v)=0$. Since the series $\h_{2}^{\pm}(u)$ is invertible, we get $\e_{23}^{\pm}(u)=0$. The second relation follows by a similar argument.
\end{proof}

\ble\label{lem:nthrootD} All relations of Lemma~{\rm \ref{lem:n-1throotD}} remain valid after the replacements
\begin{gather*}
\h_{2}^{\pm}(u)\mapsto \h_{3}^{\pm}(u),\qquad \e_{12}^{\pm}(u)\mapsto \e_{13}^{\pm}(u),
\qquad \e_{21}^{\pm}(u)\mapsto \e_{31}^{\pm}(u),\\
\f_{12}^{\pm}(u)\mapsto \f_{13}^{\pm}(u),
\qquad \f_{21}^{\pm}(u)\mapsto \f_{31}^{\pm}(u).
\end{gather*}
\ele

\begin{proof}In view of Lemma~\ref{lem:e23=0D}, this holds because the series $\ell_{ij}^{\pm}(u)$ with $i,j=1,3$ satisfy the same relations as in the algebra $U_q\big(\widehat{\mathfrak{gl}}_2\big)$.
\end{proof}

\ble\label{lem:h2h3D} In the algebra $U\big(\overline{R}^{\tss[2]}\big)$ we have
\begin{gather*}
\h_2^{\pm}(u)\h_3^{\pm}(v)=\h_3^{\pm}(v)\h_2^{\pm}(u),
\\
\frac{(q^{-1}u_{\pm}-qv_{\mp})(u_{\pm}-v_{\mp})}{\big(qu_{\pm}-q^{-1}v_{\mp}\big)
\big(u_{\pm}-q^{-1}v_{\mp}\big)}\h_2^{\pm}(u)\h_3^{\mp}(v)
=\frac{\big(q^{-1}u_{\pm}-qv_{\mp}\big)(u_{\pm}-v_{\mp})}{\big(qu_{\pm}-q^{-1}v_{\mp}\big)
\big(u_{\pm}-q^{-1}v_{\mp}\big)}\h_3^{\mp}(v)\h_2^{\pm}(u).
\end{gather*}
\ele

\begin{proof}By Corollary~\ref{cor:guass-embed} we have
\begin{gather*}
\frac{\big(q^{-1}u_{\pm}-qv_{\mp}\big)(u_{\pm}-v_{\mp})}{\big(qu_{\pm}-q^{-1}v_{\mp}\big) \big(u_{\pm}-q^{-1}v_{\mp}\big)}
\ell_{22}^{\pm\ts [1]}(u)\ell_{33}^{\pm\ts [1]}(v)\\
\qquad{} =
\frac{\big(q^{-1}u_{\pm}-qv_{\mp}\big)(u_{\pm}-v_{\mp})}{\big(qu_{\pm}-q^{-1}v_{\mp}\big) \big(u_{\pm}-q^{-1}v_{\mp}\big)}
\ell_{33}^{\pm\ts [1]}(v)\ell_{22}^{\pm\ts [1]}(u).
\end{gather*}
Writing this in terms of the Gaussian generators and using Lemma~\ref{lem:e23=0D} we get the second relation. The first relation is verified in the same way.
\end{proof}

\ble\label{lem:eeqee}In the algebra $U\big(\overline{R}^{\tss[2]}\big)$ we have
\begin{gather}
\e_{14}^{\pm}(u)=-\e_{12}^{\pm}(u)\e_{13}^{\pm}(u)=-\e_{13}^{\pm}(u)\e_{12}^{\pm}(u),\nonumber\\
\f_{41}^{\pm}(u)=-\f_{21}^{\pm}(u)\f_{31}^{\pm}(u)=-\f_{31}^{\pm}(u)\f_{21}^{\pm}(u),\label{e14=-e12e13}
\end{gather}
and
\begin{alignat}{3}
& \e_{12}^{\pm}(u)\e_{13}^{\mp}(v)
=\e_{12}^{\mp}(v)\e_{13}^{\pm}(u),\qquad&&
\e_{12}^{\pm}(u)\e_{13}^{\pm}(v)=\e_{12}^{\pm}(v)\e_{13}^{\pm}(u),&\nonumber\\
& \f_{21}^{\pm}(u)\f_{31}^{\mp}(v)
=\f_{21}^{\mp}(v)\f_{31}^{\pm}(u),\qquad &&
\f_{21}^{\pm}(u)\f_{31}^{\pm}(v)=\f_{21}^{\pm}(v)\f_{31}^{\pm}(u).& \label{e12e13=e12e13}
\end{alignat}
\ele

\begin{proof} The arguments are similar for all relations so we only give details for the first equality in
\eqref{e14=-e12e13} and the first part of \eqref{e12e13=e12e13}. The defining relations \eqref{gen rel1} give
\begin{gather*}
\ell_{12}^{\pm}(u)\ell_{13}^{\pm}(v)=
\frac{1}{\big(u/v-q^{-2}\big)^2}\sum_{i=1}^{4}a_{i2}(u/v)\ell_{1i}^{\pm}(v)\ell_{1i'}^{\pm}(u)
\end{gather*}
and
\begin{gather*}
\ell_{11}^{\pm}(u)\ell_{14}^{\pm}(v)=
\frac{1}{\big(u/v-q^{-2}\big)^2}\sum_{i=1}^{4}a_{i4}(u/v)\ell_{1i}^{\pm}(v)\ell_{1i'}^{\pm}(u).
\end{gather*}
Hence we can write
\begin{gather*}
\ell_{12}^{\pm}(u)\ell_{13}^{\pm}(v) =
\frac{\big(q^{-2}-1\big)\big(q^{-1}u/v-q\big)}{\big(u/v-q^{-2}\big)(u/v-1)}\ell_{11}^{\pm}(v)\ell_{14}^{\pm}(u)\\
 \qquad{} +\frac{q^{-2}u/v-1}{u/v-q^{-2}}\ell_{12}^{\pm}(v)\ell_{13}^{\pm}(u)
+\frac{\big(q^{-1}-q\big)u/v}{u/v-1}\ell_{11}^{\pm}(u)\ell_{14}^{\pm}(v).\label{l12l13}
\end{gather*}
Using again \eqref{gen rel1}, we get
\begin{gather*}
\ell_{11}^{\pm}(u)\ell_{12}^{\pm}(v)=\frac{u/v-1}{qu/v-q^{-1}}\ell_{12}^{\pm}(v)\ell_{11}^{\pm}(u)
+\frac{q-q^{-1}}{qu/v-q^{-1}}\ell_{11}^{\pm}(v)\ell_{12}^{\pm}(u).
\end{gather*}
Therefore, \eqref{l12l13} is equivalent to
\begin{gather*}
\frac{q^{-1}u/v-q}{u/v-1}\big(\h_{1}^{\pm}(v)\h_{1}^{\pm}(u)\e_{12}^{\pm}(u) \e_{13}^{\pm}(v)
-\h_{1}^{\pm}(u)\h_{1}^{\pm}(v)\e_{12}^{\pm}(v)\e_{13}^{\pm}(u)\big)\\
\qquad{} =
\frac{\big(q^{-2}-1\big)\big(q^{-1}u/v-q\big)}{\big(u/v-q^{-2}\big)(u/v-1)}\h_{1}^{\pm}(v)\h_{1}^{\pm}(u)
(\e_{14}^{\pm}(u)+\e_{12}^{\pm}(u)\e_{13}^{\pm}(u)\\
\qquad\quad{} +\frac{\big(q^{-1}-q\big)u/v}{u/v-1}\h_{1}^{\pm}(u)\h_{1}^{\pm}(v)
\big(\e_{14}^{\pm}(v)+\e_{12}^{\pm}(v)\e_{13}^{\pm}(v)\big).
\end{gather*}
Since $\h_{1}^{\pm}(v)\h_{1}^{\pm}(u)=\h_{1}^{\pm}(u)\h_{1}^{\pm}(v)$ and the series $\h_{1}^{\pm}(u)$
is invertible, we come to the relation
\begin{gather*}
\frac{q^{-1}u/v-q}{u/v-1}\big(\e_{12}^{\pm}(u)\e_{13}^{\pm}(v)
-\e_{12}^{\pm}(v)\e_{13}^{\pm}(u)\big) =
\frac{\big(q^{-2}-1\big)\big(q^{-1}u/v-q\big)}{\big(u/v-q^{-2}\big)(u/v-1)}
\big(\e_{14}^{\pm}(u)+\e_{12}^{\pm}(u)\e_{13}^{\pm}(u)\big)\\
\hphantom{\frac{q^{-1}u/v-q}{u/v-1}\big(\e_{12}^{\pm}(u)\e_{13}^{\pm}(v)
-\e_{12}^{\pm}(v)\e_{13}^{\pm}(u)\big) =}{}  +\frac{\big(q^{-1}-q\big)u/v}{u/v-1}\big(\e_{14}^{\pm}(v)+\e_{12}^{\pm}(v)\e_{13}^{\pm}(v)\big).
\end{gather*}
Setting $u/v=q^2$, we get $\e_{14}^{\pm}(v)+\e_{12}^{\pm}(v)\e_{13}^{\pm}(v)=0$
which is the first relation in \eqref{e14=-e12e13}.

For the proof of the first part of \eqref{e12e13=e12e13},
consider the relations
\begin{gather*}
\ell_{12}^{\pm}(u)\ell_{13}^{\mp}(v)=
\frac{1}{\big(u_{\mp}/v_{\pm}-q^{-2}\big)^2}
\sum_{i=1}^{4}a_{i2}(u_{\mp}/v_{\pm})\ell_{1i}^{\mp}(v)\ell_{1i'}^{\pm}(u)
\end{gather*}
and
\begin{gather*}
\ell_{11}^{\pm}(u)\ell_{14}^{\mp}(v)=
\frac{1}{\big(u_{\mp}/v_{\pm}-q^{-2}\big)^2}
\sum_{i=1}^{4}a_{i4}(u_{\mp}/v_{\pm})\ell_{1i}^{\mp}(v)\ell_{1i'}^{\pm}(u),
\end{gather*}
which hold by \eqref{gen rel2}. As with the above argument, they imply
\begin{gather*}
\frac{q^{-1}u_{\mp}/v_{\pm}-q}{u_{\mp}/v_{\pm}-1}\big(\e_{12}^{\pm}(u)\e_{13}^{\mp}(v)
 -\e_{12}^{\mp}(v)\e_{13}^{\pm}(u)\big)\\
\qquad{} =\frac{\big(q^{-2}-1\big)\big(q^{-1}u_{\mp}/v_{\pm}-q\big)}{\big(u_{\mp}/v_{\pm}-q^{-2}\big) \big(u_{\mp}/v_{\pm}-1\big)} \big(\e_{14}^{\pm}(u)+\e_{12}^{\pm}(u)\e_{13}^{\pm}(u)\big)\\
\qquad\quad{}+\frac{\big(q^{-1}-q\big)u_{\mp}/v_{\pm}}{u_{\mp}/v_{\pm}-1}
\big(\e_{14}^{\mp}(v)+\e_{12}^{\mp}(v)\e_{13}^{\mp}(v)\big).
\end{gather*}
Using \eqref{e14=-e12e13}, we get
$\e_{12}^{\pm}(u)\e_{13}^{\mp}(v)
-\e_{12}^{\mp}(v)\e_{13}^{\pm}(u)=0$.
\end{proof}

\ble\label{lem:e12e13} In the algebra $U\big(\overline{R}^{\tss[2]}\big)$ we have
\begin{alignat*}{3}
& \e_{12}^{\pm}(u)\e_{13}^{\mp}(v)=\e_{13}^{\mp}(v)\e_{12}^{\pm}(u),
\qquad && \e_{12}^{\pm}(u)\e_{13}^{\pm}(v)=\e_{13}^{\pm}(v)\e_{12}^{\pm}(u),&\\
& \f_{21}^{\pm}(u)\f_{31}^{\mp}(v)=\f_{31}^{\mp}(v)\f_{21}^{\pm}(u),
\qquad &&
\f_{21}^{\pm}(u)\f_{31}^{\pm}(v)=\f_{31}^{\pm}(v)\f_{21}^{\pm}(u).&
\end{alignat*}
\ele

\begin{proof}All relations are verified in the same way so we only give details for
the first one with the top signs.
By the defining relations \eqref{gen rel2}, we have
\begin{gather*}\label{l12l13D}
\ell_{12}^{\pm}(u)\ell_{13}^{\mp}(v)=
\frac{1}{\big(u_{\mp}/v_{\pm}-q^{-2}\big)^2}
\sum_{i=1}^{4}a_{i2}(u_{\mp}/v_{\pm})\ell_{1i}^{\mp}(v)\ell_{1i'}^{\pm}(u).
\end{gather*}
Using the Gauss decomposition and \eqref{e14=-e12e13},
we can write the right hand side of~\eqref{l12l13D} as
\begin{gather*}
\frac{1}{\big(x-q^{-2}\big)^2}\big( {-}a_{12}(x)\h_{1}^{-}(v)\h_{1}^{+}(u)\e_{12}^{+}(u)\e_{13}^{+}(u)
+a_{22}(x)\h_{1}^{-}(v)\e_{12}^{-}(v)\h_{1}^{+}(u)\e_{13}^{+}(u)\\
\qquad{} +a_{32}(x)\h_{1}^{-}(v)\e_{13}^{-}(v)\h_{1}^{+}(u)\e_{12}^{+}(u)
-a_{42}(x)\h_{1}^{-}(v)\e_{12}^{-}(v)\e_{13}^{-}(v)\h_{1}^{+}(u)
\big),
\end{gather*}
where $x=u_{-}/v_{+}$. Note that $\e_{12}^{+}(u)\e_{13}^{+}(u)=\e_{13}^{+}(u)\e_{12}^{+}(u)$ by~\eqref{e14=-e12e13}.
Hence, using the relations between $\h_{1}^{+}(u)$ and the series
$\e_{12}^{-}(v)$ and $\e_{13}^{-}(v)$, provided by Lemmas~\ref{lem:n-1throotD}
and \ref{lem:nthrootD}, we can write the right hand side of~\eqref{l12l13D} in the form
\begin{gather*}
\h_{1}^{-}(v)\h_{1}^{+}(u)  \bigg(\frac{q^{-1}x-1}{x-1}\e_{12}^{-}(v)\e_{13}^{+}(u)
+\frac{\big(q-q^{-1}\big)x}{x-1}\e_{12}^{-}(v)\e_{13}^{-}(v)\\
\qquad{} +\frac{\big(1-q^{-2}\big)^2 qx}{(x-1)\big(x-q^{-2}\big)}
\big(\e_{13}^{-}(v)\e_{12}^{+}(u)-\e_{12}^{+}(u)\e_{13}^{-}(v)\big)
\bigg).
\end{gather*}
On the other hand, by the relations between $\e_{12}^{+}(u)$
and $\h_{1}^{-}(v)$, the left hand side of~\eqref{l12l13D} can be written as
\begin{gather*}
\h_{1}^{+}(u)\h_{1}^{-}(v) \left(\frac{q^{-1}x-1}{x-1}\e_{12}^{+}(u)\e_{13}^{-}(v)
+\frac{\big(q-q^{-1}\big)x}{x-1}\e_{12}^{-}(v)\e_{13}^{-}(v)\right).
\end{gather*}
Hence, due to \eqref{e12e13=e12e13} and the property
$\h_{1}^{-}(v)\h_{1}^{+}(u)=\h_{1}^{+}(u)\h_{1}^{-}(v)$ we get
\[ \e_{13}^{-}(v)\e_{12}^{+}(u)=\e_{12}^{+}(u)\e_{13}^{-}(v),\]
as required.
\end{proof}

\ble\label{lem:e12=-e34} In the algebra $U\big(\overline{R}^{\tss[2]}\big)$ we have
\begin{gather*}
\e_{24}^{\pm}(u)=-\e_{13}^{\pm}(u),\qquad \e_{34}^{\pm}(u)=-\e_{12}^{\pm}(u),
\qquad
\f_{43}^{\pm}(u)=-\f_{21}^{\pm}(u),
\qquad
\f_{42}^{\pm}(u)=-\f_{31}^{\pm}(u).
\end{gather*}
\ele

\begin{proof}We only verify the first relation. By Proposition~\ref{prop:central}, we have the matrix relation
\begin{gather*}
\Lc^{\pm}(u)^{-1}\z^{\pm\ts[2]}(u)=D^{[2]}
\Lc^{\pm}\big(uq^{-2}\big)^{\tra}\big(D^{[2]}\big)^{-1}.
\end{gather*}
Take $(4,4)$ and $(2,4)$-entries on both sides and use the property $\e_{23}^{\pm}(u)=0$, which holds by Lemma~\ref{lem:e23=0D}, to get
\begin{gather*}
\h_{1}^{\pm}\big(uq^{-2}\big)=\h_{4}^{\pm}(u)^{-1}\z^{\pm\ts[2]}(u)
\end{gather*}
and
\begin{gather*}
q\ts\h_{1}^{\pm}\big(uq^{-2}\big)\e_{13}^{\pm}\big(uq^{-2}\big) =-\e_{24}^{\pm}(u)\h_{4}^{\pm}(u)^{-1}\z^{\pm\ts[2]}(u).
\end{gather*}
This implies
\begin{gather*}
q\ts\h_{1}^{\pm}\big(uq^{-2}\big)\e_{13}^{\pm}\big(uq^{-2}\big)=-\e_{24}^{\pm}(u)\h_{1}^{\pm}\big(uq^{-2}\big).
\end{gather*}
By the relations between $\h_{1}^{\pm}(u)$ and $\e_{13}^{\pm}(v)$ from Lemma~\ref{lem:nthrootD},
we also have
\begin{gather*}
q\ts\h_{1}^{\pm}\big(uq^{-2}\big)\e_{13}^{\pm}\big(uq^{-2}\big)=\e_{13}^{\pm}(u)\h_{1}^{\pm}\big(uq^{-2}\big).
\end{gather*}
By comparing the two formulas we conclude that $\e_{13}^{\pm}(u)=-\e_{24}^{\pm}(u)$.
\end{proof}

\ble\label{lem:e12f31D} We have the relations
\begin{alignat}{3}&
\big[\e_{12}^{\pm}(u),\f_{31}^{\mp}(v)\big]=0,\qquad && \big[\e_{12}^{\pm}(u),\f_{31}^{\pm}(v)\big]=0,& \label{e12f31}
\\
\label{e13f21}
& \big[\e_{13}^{\pm}(u),\f_{21}^{\mp}(v)\big]=0,\qquad && \big[\e_{13}^{\pm}(u),\f_{21}^{\pm}(v)\big]=0,&
\end{alignat}
and
\begin{gather}
\e_{12}^{\pm}(u)\h_{3}^{\mp}(v)
=\frac{q^{-1}u_{\mp}-qv_{\pm}}{u_{\mp}-v_{\pm}}\h_{3}^{\mp}(v)\e_{12}^{\pm}(u)
+\frac{\big(q-q^{-1}\big)u_{\mp}}{u_{\mp}-v_{\pm}}\h_{3}^{\mp}(v)\e_{12}^{\mp}(v),\nonumber\\
\e_{12}^{\pm}(u)\h_{3}^{\pm}(v) =\frac{q^{-1}u-qv}{u-v}\h_{3}^{\pm}(v)\e_{12}^{\pm}(u)
+\frac{\big(q-q^{-1}\big)u}{u-v}\h_{3}^{\pm}(v)\e_{12}^{\pm}(v),
\nonumber\\
\h_{3}^{\mp}(v)\f_{21}^{\pm}(u)
=\frac{q^{-1}u_{\pm}-qv_{\mp}}{u_{\pm}-v_{\mp}}\f_{21}^{\pm}(u)\h_{3}^{\mp}(v)
+\frac{\big(q-q^{-1}\big)v_{\mp}}{u_{\pm}-v_{\mp}}\f_{21}^{\mp}(v)\h_{3}^{\mp}(v),
\nonumber\\
\h_{3}^{\pm}(v)\f_{21}^{\pm}(u) =\frac{q^{-1}u-qv}{u-v}\f_{21}^{\pm}(u)\h_{3}^{\pm}(v)
+\frac{\big(q-q^{-1}\big)v}{u-v}\f_{21}^{\pm}(v)\h_{3}^{\pm}(v).\label{e12h3pm}
\end{gather}
\ele

\begin{proof}We only give a proof of one case of \eqref{e12f31} and \eqref{e12h3pm}, the remaining
relations are verified in a similar way.
As before, we set $x=u_{+}/v_{-}$ and $y=u_{-}/v_{+}$.
The defining relations \eqref{gen rel2} imply
\begin{gather}
\frac{x-1}{qx-q^{-1}}\ell_{12}^{+}(u)\ell_{31}^{-}(v)+
\frac{\big(q-q^{-1}\big)x}{qx-q^{-1}}\ell_{32}^{+}(u)\ell_{11}^{-}(v)\nonumber\\
\qquad{} =\frac{y-1}{qy-q^{-1}}\ell_{31}^{-}(v)\ell_{12}^{+}(u)+
\frac{\big(q-q^{-1}\big)y}{qy-q^{-1}}\ell_{32}^{-}(v)\ell_{11}^{+}(u).\label{l12l31}
\end{gather}
Taking into account Lemma~\ref{lem:e23=0D}, we can write the right hand side as
\begin{gather*}
\f_{31}^{-}(v)\left(\frac{y-1}{qy-q^{-1}}\ell_{11}^{-}(v)\ell_{12}^{+}(u)+
\frac{\big(q-q^{-1}\big)y}{qy-q^{-1}}\ell_{12}^{-}(v)\ell_{11}^{+}(u)\right).
\end{gather*}
Using again \eqref{gen rel2}, we get
\begin{gather*}
\ell_{12}^{+}(u)\ell_{11}^{-}(v)=\frac{y-1}{qy-q^{-1}}\ell_{11}^{-}(v)\ell_{12}^{+}(u)+
\frac{\big(q-q^{-1}\big)y}{qy-q^{-1}}\ell_{12}^{-}(v)\ell_{11}^{+}(u).
\end{gather*}
Therefore, \eqref{l12l31} is equivalent to
\begin{gather}\label{l12l31second}
\frac{x-1}{qx-q^{-1}}\ell_{12}^{+}(u)\ell_{31}^{-}(v)+
\frac{\big(q-q^{-1}\big)x}{qx-q^{-1}}\ell_{32}^{+}(u)\ell_{11}^{-}(v)
=\f_{31}^{-}(v)\ell_{12}^{+}(u)\ell_{11}^{-}(v).
\end{gather}
By using the relation between $\f_{31}^{-}(v)$ and $\h_{1}^{+}(u)$ from Lemma~\ref{lem:nthrootD}
bring the right hand side to the form
\begin{gather*}
\frac{x-1}{qx-q^{-1}}\h_{1}^{+}(u)\f_{31}^{-}(v)\e_{12}^{+}(u)\h_{1}^{-}(v)+
\frac{\big(q-q^{-1}\big)x}{qx-q^{-1}}\f_{31}^{+}(u)\h_{1}^{+}(u)\e_{12}^{+}(u)\h_{1}^{-}(v).
\end{gather*}
On the other hand, Lemma~\ref{lem:e23=0D} implies that the left hand side of~\eqref{l12l31second} equals
\begin{gather*}
\frac{x-1}{qx-q^{-1}}\h_{1}^{+}(u)\e_{12}^{+}(u)\f_{31}^{-}(v)\h_{1}^{-}(v)+
\frac{\big(q-q^{-1}\big)x}{qx-q^{-1}}\f_{31}^{+}(u)\h_{1}^{+}(u)\e_{12}^{+}(u)\h_{1}^{-}(v),
\end{gather*}
thus proving that $[\e_{12}^{+}(u),\f_{31}^{-}(v)]=0$.

Now turn to \eqref{e12h3pm}. The defining relations \eqref{gen rel2} give
\begin{gather}\label{l12l33}
\frac{x-1}{qx-q^{-1}}\ell_{12}^{+}(u)\ell_{33}^{-}(v)+
\frac{\big(q-q^{-1}\big)x}{qx-q^{-1}}\ell_{32}^{+}(u)\ell_{13}^{-}(v)=
\frac{1}{\big(y-q^{-2}\big)^2}\sum_{i=1}^{4} a_{i3}(y)\ell_{3i}^{-}(v)\ell_{1i'}^{+}(u).\!\!\!
\end{gather}
By Lemma~\ref{lem:e23=0D}, the left hand side can be written as
\begin{gather*}
\frac{x-1}{qx-q^{-1}}\h_{1}^{+}(u) \e_{12}^{+}(u)\h_{3}^{-}(v)\\
\qquad{} +\left(\frac{x-1}{qx-q^{-1}}\h_{1}^{+}(u)\e_{12}^{+}(u)\f_{31}^{-}(v)
+\frac{\big(q-q^{-1}\big)x}{qx-q^{-1}}\f_{31}^{+}(u)\h_{1}^{+}(u)\e_{12}^{+}(u)\right)
\ts\h_{1}^{-}(v)\e_{13}^{-}(v).
\end{gather*}
Due to \eqref{e12f31}, this expression equals
\begin{gather*}
\frac{x-1}{qx-q^{-1}}\h_{1}^{+}(u) \e_{12}^{+}(u)\h_{3}^{-}(v)\\
 \qquad{} +\left(\frac{x-1}{qx-q^{-1}}\h_{1}^{+}(u)\f_{31}^{-}(v)
+\frac{\big(q-q^{-1}\big)x}{qx-q^{-1}}\f_{31}^{+}(u)\h_{1}^{+}(u)\right)
\ts\e_{12}^{+}(u)\h_{1}^{-}(v)\e_{13}^{-}(v),
\end{gather*}
which simplifies further to
\begin{gather*}
\frac{x-1}{qx-q^{-1}}\h_{1}^{+}(u)\e_{12}^{+}(u)\h_{3}^{-}(v)
+\f_{31}^{-}(v)\ell_{12}^{+}(u)\ell_{13}^{-}(v)
\end{gather*}
by the relation between $\h_{1}^{+}(u)$ and $\f_{31}^{-}(v)$ provided
by Lemma~\ref{lem:nthrootD}. Furthermore, by~\eqref{gen rel2} we also have
\begin{gather*}
\ell_{12}^{+}(u)\ell_{13}^{-}(v)=
\frac{1}{\big(y-q^{-2}\big)^2}\sum_{i=1}^{4} a_{i3}(y)\ell_{1i}^{-}(v)\ell_{1i'}^{+}(u)
\end{gather*}
so that the left hand side of \eqref{l12l33} becomes
\begin{gather*}
\frac{x-1}{qx-q^{-1}}\h_{1}^{+}(u)\e_{12}^{+}(u)\h_{3}^{-}(v)
+\frac{\f_{31}^{-}(v)}{\big(y-q^{-2}\big)^2}\sum_{i=1}^{4} a_{i3}(y)\ell_{1i}^{-}(v)\ell_{1i'}^{+}(u).
\end{gather*}
Using Lemmas~\ref{lem:e23=0D} and \ref{lem:e12=-e34},
in terms of Gaussian generators we get
\begin{gather*}
\frac{x-1}{qx-q^{-1}}\h_{1}^{+}(u)\e_{12}^{+}(u)\h_{3}^{-}(v)\\
\qquad {} =\frac{y-1}{\big(qy-q^{-1}\big)^2}\big( (y-1)\h_{3}^{-}(v)\h_{1}^{+}(u)\e_{12}^{+}(u)
+\big(q-q^{-1}\big)y\h_{3}^{-}(v)\e_{12}^{-}(v)\h_{1}^{+}(u)\big).
\end{gather*}
As a final step, use the relations between $\h_{1}^{+}(u)$ and $\e_{12}^{-}(v)$
and those between $\h_{1}^{\pm}(u)$ and $\h_{3}^{\mp}(v)$
from Lemmas~\ref{lem:n-1throotD} and~\ref{lem:nthrootD}, respectively, to come to the relation
\begin{gather*}
\e_{12}^{+}(u)\h_{3}^{-}(v)=
\frac{q^{-1}y-q}{y-1}\h_{3}^{-}(v)\e_{12}^{+}(u)
+\frac{\big(q-q^{-1}\big)y}{y-1}\h_{3}^{-}(v)\e_{12}^{-}(v),
\end{gather*}
as required.
\end{proof}

\ble\label{lem:e13h2D} In the algebra $U\big(\overline{R}^{\tss[2]}\big)$ we have
\begin{gather*}
\e_{13}^{\pm}(u)\h_{2}^{\mp}(v)=
\frac{q^{-1}u_{\mp}-qv_{\pm}}{u_{\mp}-v_{\pm}}\h_{2}^{\mp}(v)\e_{13}^{\pm}(u)
+\frac{\big(q-q^{-1}\big)u_{\mp}}{u_{\mp}-v_{\pm}}\h_{2}^{\mp}(v)\e_{13}^{\mp}(v),\\
\e_{13}^{\pm}(u)\h_{2}^{\pm}(v)=\frac{q^{-1}u-qv}{u-v}\h_{2}^{\pm}(v)\e_{13}^{\pm}(u)
+\frac{\big(q-q^{-1}\big)u}{u-v}\h_{2}^{\pm}(v)\e_{13}^{\pm}(v),\\
\h_{2}^{\mp}(v)\f_{31}^{\pm}(u)=
\frac{q^{-1}u_{\pm}-qv_{\mp}}{u_{\pm}-v_{\mp}}\f_{31}^{\pm}(u)\h_{2}^{\mp}(v)
+\frac{\big(q-q^{-1}\big)v_{\mp}}{u_{\pm}-v_{\mp}}\f_{31}^{\mp}(v)\h_{2}^{\mp}(v),\\
\h_{2}^{\pm}(v)\f_{31}^{\pm}(u)=
\frac{q^{-1}u-qv}{u-v}\f_{31}^{\pm}(u)\h_{2}^{\pm}(v)
+\frac{\big(q-q^{-1}\big)v}{u-v}\f_{31}^{\pm}(v)\h_{2}^{\pm}(v).
\end{gather*}
\ele

\begin{proof}The arguments for all relations are quite similar so we only give details for one case of the first relation. By \eqref{gen rel2} we have
\begin{gather}\label{l13l22}
\frac{x-1}{qx-q^{-1}}\ell_{13}^{+}(u)\ell_{22}^{-}(v)
+\frac{\big(q-q^{-1}\big)x}{qx-q^{-1}}\ell_{23}^{+}(u)\ell_{12}^{-}(v)
=\frac{1}{\big(y-q^{-2}\big)^2}\sum_{i=1}^4 a_{i2}(y) \ell_{2i}^{-}(v)\ell_{1i'}^{+}(u).\!\!\!\!
\end{gather}
Taking into account Lemma~\ref{lem:e23=0D}, write the left hand side as
\begin{gather*}
\frac{x-1}{qx-q^{-1}}\h_{1}^{+}(u)\e_{13}^{+}(u)\h_{2}^{-}(v)
 +\frac{x-1}{qx-q^{-1}}\h_{1}^{+}(u)\e_{13}^{+}(u)\f_{21}^{-}(v)\ell_{12}^{-}(v)\\
\qquad{} +\frac{\big(q-q^{-1}\big)x}{qx-q^{-1}}\f_{21}^{+}(u)\h_{1}^{+}(u)\e_{13}^{+}(u)\ell_{12}^{-}(v).
\end{gather*}
By \eqref{e13f21} this equals
\begin{gather*}
\frac{x-1}{qx-q^{-1}}\h_{1}^{+}(u)\e_{13}^{+}(u)\h_{2}^{-}(v)
 +\frac{x-1}{qx-q^{-1}}\h_{1}^{+}(u)\f_{21}^{-}(v)\e_{13}^{+}(u)\ell_{12}^{-}(v)\\
\qquad{} +\frac{\big(q-q^{-1}\big)x}{qx-q^{-1}}\f_{21}^{+}(u)\h_{1}^{+}(u)\e_{13}^{+}(u)\ell_{12}^{-}(v).
\end{gather*}
Then by using the relation between $\h_{1}^{+}(u)$ and $\f_{21}^{-}(v)$ from Lemma~\ref{lem:n-1throotD}, we bring the left hand side of \eqref{l13l22} to the form
\begin{gather*}
\frac{x-1}{qx-q^{-1}}\h_{1}^{+}(u)\e_{13}^{+}(u)\h_{2}^{-}(v)
+\f_{21}^{-}(v)\h_{1}^{+}(u)\e_{13}^{+}(u)\ell_{12}^{-}(v).
\end{gather*}
By the defining relations between $\ell_{13}^{+}(u)$ and $\ell_{12}^{-}(v)$ we have
\begin{gather*}
\e_{13}^{+}(u)\ell_{12}^{-}(v)=\frac{1}{\big(y-q^{-2}\big)^2}
\sum_{i=1}^4 a_{i2}(y) \ell_{1i}^{-}(v)\ell_{1i'}^{+}(u)
\end{gather*}
and so the left hand side of \eqref{l13l22} can be written as
\begin{gather*}
\frac{x-1}{qx-q^{-1}}\h_{1}^{+}(u)\e_{13}^{+}(u)\h_{2}^{-}(v)
+\frac{\f_{21}^{-}(v)}{\big(y-q^{-2}\big)^2}\sum_{i=1}^4 a_{i2}(y) \ell_{1i}^{-}(v)\ell_{1i'}^{+}(u).
\end{gather*}
Hence by Lemma~\ref{lem:e23=0D} relation \eqref{l13l22} now reads
\begin{gather}
\frac{x-1}{qx-q^{-1}}\h_{1}^{+}(u)\e_{13}^{+}(u)\h_{2}^{-}(v)\nonumber\\
\qquad{} =\frac{1}{(y-q^{-2})^2}\big(a_{22}(y) \h_{2}^{-}(v)\h_{1}^{+}(u)\e_{13}^{+}(u)
+a_{42}(y) \h_{2}^{-}(v)\e_{24}^{-}(v)\h_{1}^{+}(u)\big).\label{l13l22second}
\end{gather}
Using the equality $\e_{24}^{-}(v)=-\e_{13}^{-}(v)$ from Lemma~\ref{lem:e12=-e34} and
the relations between~$\h_{1}^{\pm}(u)$ and~$\e_{13}^{\mp}(v)$ from Lemma~\ref{lem:nthrootD},
we find that the right hand side of~\eqref{l13l22second} equals
\begin{gather*}
\frac{q^{-1}y-q}{y-1}\frac{x-1}{qx-q^{-1}}\h_1^{+}(u)\h_2^{-}(v)\e_{13}^{+}(u)+
\frac{\big(q-q^{-1}\big)y}{y-1}\frac{x-1}{qx-q^{-1}}\h_1^{+}(u)\h_2^{-}(v)\e_{13}^{-}(v),
\end{gather*}
where we also applied the relations between $\h_1^{+}(u)$ and $\h_2^{-}(v)$. Now~\eqref{l13l22second} turns into one case of the first relation due to the invertibility of $\h_1^{+}(u)$.
\end{proof}

\subsection[Formulas for the series $z^{\pm}(u)$ and $\z^{\pm}(u)$]{Formulas for the series $\boldsymbol{z^{\pm}(u)}$ and $\boldsymbol{\z^{\pm}(u)}$}\label{subsec:fsz}

We will now consider the cases of odd and even $N$ simultaneously, unless stated otherwise.
Recall that the series $z^{\pm}(u)$ and $\z^{\pm}(u)$ were defined in Proposition~\ref{prop:central}. We will now indicate the dependence on $n$ by adding the corresponding superscript.
Write relation~\eqref{DLbarDLbar} in the form
\begin{gather}\label{DLbarD}
D\Lc^{\pm}(u\xi)^{\tra}D^{-1}=\Lc^{\pm}(u)^{-1}\z^{\pm\ts[n]}(u).
\end{gather}
Using the Gauss decomposition for $\Lc^{\pm}(u)$ and taking the $(N,N)$-entry on both sides of~\eqref{DLbarD} we get
\begin{gather}\label{h1h1pr}
\h_1^{\pm}(u\xi)=\h_{1'}^{\pm}(u)^{-1}\z^{\pm\ts[n]}(u).
\end{gather}

\ble\label{lem:eiprei}
The following relations hold in the algebra $U\big(\overline{R}^{\tss[n]}\big)$:
\begin{gather}\label{ei'ei}
\e^{\pm}_{(i+1)'\ts i'}(u)=-\e_{i,i+1}^{\pm}\big(u\xi q^{2i}\big)
\qquad \text{and} \qquad
\f^{\pm}_{i'\ts (i+1)'}(u)=-\f_{i+1,i}^{\pm}\big(u\xi q^{2i}\big)
\end{gather}
for $1\leqslant i\leqslant n-1$.
\ele

\begin{proof}
By Propositions~\ref{prop:central} and~\ref{prop:gauss-consist}, for any $1\leqslant i\leqslant n-1$
we have
\begin{gather}\label{Lbar[n-i+1] z}
\Lc^{\pm\ts[n-i+1]}(u)^{-1}\z^{\pm\ts[n-i+1]}(u)=D^{[n-i+1]}
\Lc^{\pm\ts[n-i+1]}\big(u\xi q^{2i-2}\big)^{\tra}\big(D^{[n-i+1]}\big)^{-1},
\end{gather}
where
\begin{gather*}
D^{[n-i+1]}=\begin{cases}
\diag\tss \big[\ts q^{n-i+1/2},\dots,q^{1/2},\ts 1,\ts q^{-1/2},\dots,q^{-n+i-1/2}\tss\big]
 &\text{for type $B$},\\
\diag\tss \big[\ts q^{n-i},\dots,q,\ts 1,\ts 1,\ts q^{-1},\dots,q^{-n+i}\tss\big]
 &\text{for type $D$}.
 \end{cases}
\end{gather*}
By taking the $(i',i')$ and $((i+1)',i')$-entries on both sides of~\eqref{Lbar[n-i+1] z} we get
\begin{gather}\label{hihipr}
\h_i^{\pm}\big(u\xi q^{2i-2}\big)=\h_{i'}^{\pm}(u)^{-1}\z^{\pm\ts[n-i+1]}(u)
\end{gather}
and
\begin{gather*}
-\e^{\pm}_{(i+1)',i'}(u)\ts\h^{\pm}_{i'}(u)^{-1}\ts\z^{\pm\ts[n-i+1]}(u)
=q\ts \h_i^{\pm}\big(u\xi q^{2i-2}\big)\e_{i,i+1}^{\pm}\big(u\xi q^{2i-2}\big).
\end{gather*}
Due to \eqref{hihipr}, this formula can be written as
\begin{gather}\label{ei'hi}
-\e^{\pm}_{(i+1)',i'}(u)\ts\h_i^{\pm}\big(u\xi q^{2i-2}\big)
=q\ts \h_i^{\pm}\big(u\xi q^{2i-2}\big)\ts\e_{i,i+1}^{\pm}\big(u\xi q^{2i-2}\big).
\end{gather}
Furthermore, by the results of \cite{df:it},
\begin{gather*}
q\ts\h_i^{\pm}(u)\ts\e_{i,i+1}^{\pm}(u)=\e_{i,i+1}^{\pm}\big(uq^2\big)\ts\h_i^{\pm}(u),
\end{gather*}
so that \eqref{ei'hi} is equivalent to
\begin{gather*}
-\e^{\pm}_{(i+1)',i'}(u)\ts\h_i^{\pm}\big(u\xi q^{2i-2}\big)
=\e_{i,i+1}^{\pm}\big(u\xi q^{2i}\big)\ts\h_i^{\pm}\big(u\xi q^{2i-2}\big),
\end{gather*}
thus proving the first relation in \eqref{ei'ei}. The second relation is verified in a similar way.
\end{proof}

\bpr\label{prop:formze}
In the algebras $U\big(\overline{R}^{\tss[n]}\big)$ and $U\big({R}^{[n]}\big)$ we have
the respective formulas:
\begin{gather*}
\z^{\pm\ts[n]}(u) =\prod_{i=1}^{n}\h^{\pm}_{i}\big(u\tss\xi q^{2i}\big)^{-1}
\prod_{i=1}^{n}\h^{\pm}_{i}\big(u\tss\xi q^{2i-2}\big)\ts\cdot \h^{\pm}_{n+1}(u)\h^{\pm}_{n+1}(uq),\\
{z}^{\pm\ts[n]}(u) =\prod_{i=1}^{n}{h}^{\pm}_{i}\big(u\xi q^{2i}\big)^{-1}
\prod_{i=1}^{n}{h}^{\pm}_{i}\big(u\xi q^{2i-2}\big)\ts\cdot {h}^{\pm}_{n+1}(u)\ts{h}^{\pm}_{n+1}(uq)
\end{gather*}
for type $B$, and
\begin{gather*}
\z^{\pm\ts[n]}(u) =\prod_{i=1}^{n-1}\h^{\pm}_{i}\big(u\tss\xi q^{2i}\big)^{-1}
\prod_{i=1}^{n-1}\h^{\pm}_{i}\big(u\tss\xi q^{2i-2}\big)\ts\cdot\h^{\pm}_{n}(u)\h^{\pm}_{n+1}(u),\\
{z}^{\pm\ts[n]}(u) =\prod_{i=1}^{n-1}{h}^{\pm}_{i}\big(u\xi q^{2i}\big)^{-1}
\prod_{i=1}^{n-1}{h}^{\pm}_{i}\big(u\xi q^{2i-2}\big)\ts\cdot {h}^{\pm}_{n}(u)\ts{h}^{\pm}_{n+1}(u)
\end{gather*}
and for type $D$.
\epr

\begin{proof}The arguments for both formulas are quite similar so we only give a proof of the first ones
for types~$B$ and~$D$. Taking the $(2^{\tss\prime},2^{\tss\prime})$-entry on both sides of~\eqref{Lbar[n-i+1] z} and expressing the entries of the matrices $\Lc^{\pm\ts[n]}(u)^{-1}$ and $\Lc^{\pm\ts[n]}(u\xi)^{\tra}$ in terms of the Gauss generators, we get
\begin{gather*}
\h^{\pm}_2(u\xi)+\f^{\pm}_{21}(u\xi)\ts\h^{\pm}_{1}(u\xi)\ts\e^{\pm}_{12}(u\xi)
=\big(\h^{\pm}_{2^{\tss\prime}}(u)^{-1}
+\e^{\pm}_{2^{\tss\prime},1'}(u)\h^{\pm}_{1'}(u)^{-1}
\f^{\pm}_{1',2^{\tss\prime}}(u)^{-1}\big)\z^{\pm\ts[n]}(u).
\end{gather*}
As we pointed out in Remark~\ref{rem:noncent}, the coefficients of the series $\z^{\pm\ts[n]}(u)$ are central in the respective subalgebras generated by the coefficients of $\ell^{\pm\ts[n]}_{ij}(u)$.
Therefore, using~\eqref{h1h1pr}, we can rewrite the above relation as
\begin{gather*}
\h^{\pm}_{2^{\tss\prime}}(u)^{-1}\z^{\pm\ts[n]}(u)=
\h^{\pm}_2(u\xi)+\f^{\pm}_{21}(u\xi)\ts\h^{\pm}_{1}(u\xi)\ts\e^{\pm}_{12}(u\xi)
-\e^{\pm}_{2^{\tss\prime},1'}(u)\ts\h^{\pm}_{1}(u\xi)\ts\f^{\pm}_{1',2^{\tss\prime}}(u).
\end{gather*}
Now apply Lemma~\ref{lem:eiprei} to obtain
\begin{gather*}
\h^{\pm}_{2^{\tss\prime}}(u)^{-1}\z^{\pm\ts[n]}(u)=
\h^{\pm}_2(u\xi)+\f^{\pm}_{21}(u\xi)\ts\h^{\pm}_{1}(u\xi)\ts\e^{\pm}_{12}(u\xi)
-\e^{\pm}_{12}\big(u\xi q^2\big)\ts\h^{\pm}_{1}(u\xi)\ts\f^{\pm}_{21}\big(u\xi q^2\big).
\end{gather*}
On the other hand, by the results of \cite{df:it} we have
\begin{gather*}
\h^{\pm}_{1}(u)\ts\e^{\pm}_{12}(u)=q^{-1}\e^{\pm}_{12}\big(uq^2\big)\ts\h^{\pm}_{1}(u),
\qquad
\h^{\pm}_{1}(u)\ts\f^{\pm}_{21}\big(uq^2\big)=q^{-1}\f^{\pm}_{21}(u)\ts\h^{\pm}_{1}(u),
\end{gather*}
and
\begin{gather*}
[\e^{\pm}_{12}(u),\f^{\pm}_{21}(v)]=
\frac{u\big(q-q^{-1}\big)}{u-v}\ts\big(\h^{\pm}_{2}(v)\h^{\pm}_{1}(v)^{-1}
-\h^{\pm}_{2}(u)\ts\h^{\pm}_{1}(u)^{-1}\big).
\end{gather*}
This leads to the expression
\begin{gather*}
\h^{\pm}_{2^{\tss\prime}}(u)^{-1}\z^{\pm\ts[n]}(u)=
\h^{\pm}_2\big(u\xi q^2\big)\ts\h^{\pm}_1\big(u\xi q^2\big)^{-1}\h^{\pm}_{1}(u\xi).
\end{gather*}
Since $\z^{\pm\ts[n-1]}(u)=\h^{\pm}_{2^{\tss\prime}}(u)\h^{\pm}_2\big(u\xi q^2\big)$, we get a recurrence formula
\begin{gather*}
\z^{\pm\ts[n]}(u)= \h^{\pm}_1\big(u\xi q^2\big)^{-1}\h^{\pm}_{1}(u\xi)\ts\z^{\pm\ts[n-1]}(u).
\end{gather*}
Here we need note that $\xi=q^{2-N}$. To complete the proof, we only need the formulas of $\z^{\pm\ts [1]}(u)$. Working with the algebras $U\big(\overline{R}^{\tss[1]}\big)$ and $U\big(\overline{R}^{\tss[2]}\big)$, respectively, we find by a~similar argument to the above that
\begin{gather*}
\h^{\pm}_{n+1}(u)^{-1}\z^{\pm\ts[1]}(u)=
\h^{\pm}_{n+1}(uq)\ts\h^{\pm}_n(uq)^{-1}\h^{\pm}_{n}\big(uq^{-1}\big)
\end{gather*}
for type $B$, and
\begin{gather*}
\z^{\pm\ts[1]}(u)=
\h^{\pm}_{n}(u)\h^{\pm}_{n+1}(u)
\end{gather*}
for type $D$.
\end{proof}

\subsection[Drinfeld-type relations in the algebras $U\big(\overline{R}^{\tss[n]}\big)$ and $U\big(R^{[n]}_{}\big)$]{Drinfeld-type relations in the algebras $\boldsymbol{U\big(\overline{R}^{\tss[n]}\big)}$ and $\boldsymbol{U\big(R^{[n]}_{}\big)}$}\label{subsec:dtr}

We will now extend the sets of relations produced in Sections~\ref{subsec:typea}, \ref{subsec:lowrankB} and \ref{subsec:lowrankD} to obtain all necessary relations in the algebras $U\big(\overline{R}^{\tss[n]}\big)$ and $U\big(R_{}^{[n]}\big)$ to be able to prove the Main Theorem. We begin by stating three lemmas which are immediate consequences of Corollary~\ref{cor:commu}.

\ble\label{lem:hn+1B}
In the algebra $U\big(\overline{R}^{\tss[n]}\big)$ we have
\begin{gather*}
\h^{\pm}_i(u)\h^{\pm}_{n+1}(v)=\h^{\pm}_{n+1}(v)\h^{\pm}_i(u),\nonumber\\
\frac{u_{\pm}-v_{\mp}}{qu_{\pm}-q^{-1}v_{\mp}}\h^{\pm}_i(u)\h^{\mp}_{n+1}(v)
=\frac{u_{\mp}-v_{\pm}}{qu_{\mp}-q^{-1}v_{\pm}}\h^{\pm}_{n+1}(v)\h^{\pm}_i(u),
\end{gather*}
and
\begin{alignat*}{3}
& \e_{i,i+1}^{\pm}(u)\h_{n+1}^{\mp}(v) =\h_{n+1}^{\mp}(v)\e_{i,i+1}^{\pm}(u),
\qquad &&
\e_{i,i+1}^{\pm}(u)\h_{n+1}^{\pm}(v) =\h_{n+1}^{\pm}(v)\e_{i,i+1}^{\pm}(u),&\\
& \f_{i+1,i}^{\pm}(u)\h_{n+1}^{\mp}(v) =\h_{n+1}^{\mp}(v)\f_{i+1,i}^{\pm}(u),
\qquad&&
\f_{i+1,i}^{\pm}(u)\h_{n+1}^{\pm}(v) =\h_{n+1}^{\pm}(v)\f_{i+1,i}^{\pm}(u),&
\end{alignat*}
where $i=1,\dots,n-1$ for type $B$, and $i=1,\dots,n-2$ for type $D$.
\ele

\ble\label{lem:en+1B}In the algebra $U\big(R^{[n]}\big)$ we have
\begin{alignat*}{3}
&\h_i^{\pm}(u)\e_{n,n+1}^{\pm}(v)=\e_{n,n+1}^{\pm}(v)\h_i^{\pm}(u),\qquad&&
\h_i^{\pm}(u)\e_{n,n+1}^{\mp}(v)=\e_{n,n+1}^{\mp}(v)\h_i^{\pm}(u),&\\
&\h_i^{\pm}(u)\f_{n+1,n}^{\pm}(v)=\f_{n+1,n}^{\pm}(v)\h_i^{\pm}(u),\qquad&&
\h_i^{\pm}(u)\f_{n+1,n}^{\mp}(v)=\f_{n+1,n}^{\mp}(v)\h_i^{\pm}(u)&
\end{alignat*}
for $i=1,\dots,n-1$ in type $B$, while
\begin{alignat*}{3}
&\h_i^{\pm}(u)\e_{n-1,n+1}^{\pm}(v)=\e_{n-1,n+1}^{\pm}(v)\h_i^{\pm}(u),\qquad&&
\h_i^{\pm}(u)\e_{n-1,n+1}^{\mp}(v)=\e_{n-1,n+1}^{\mp}(v)\h_i^{\pm}(u),&\\
&\h_i^{\pm}(u)\f_{n+1,n-1}^{\pm}(v)=\f_{n+1,n-1}^{\pm}(v)\h_i^{\pm}(u),\qquad&&
\h_i^{\pm}(u)\f_{n+1,n-1}^{\mp}(v)=\f_{n+1,n-1}^{\mp}(v)\h_i^{\pm}(u)&
\end{alignat*}
for $i=1,\dots,n-2$ in type $D$.
\ele

\ble\label{lem:eienB}
In the algebra $U\big(R^{[n]}\big)$ we have
\begin{alignat*}{3}
& \e_{i,i+1}^{\pm}(u)\e_{n,n+1}^{\pm}(v)=\e_{n,n+1}^{\pm}(v)\e_{i,i+1}^{\pm}(u),\qquad&&
\e_{i,i+1}^{\pm}(u)\e_{n,n+1}^{\mp}(v)=\e_{n,n+1}^{\mp}(v)\e_{i,i+1}^{\pm}(u),&\\
&\f_{i+1,i}^{\pm}(u)\f_{n+1,n}^{\pm}(v)=\f_{n+1,n}^{\pm}(v)\f_{i+1,i}^{\pm}(u),\qquad&&
\f_{i+1,i}^{\pm}(u)\f_{n+1,n}^{\mp}(v)=\f_{n+1,n}^{\mp}(v)\f_{i+1,i}^{\pm}(u)&
\end{alignat*}
for $i=1,\dots,n-2$ in type $B$, while
\begin{alignat*}{3}
& \e_{i,i+1}^{\pm}(u)\e_{n-1,n+1}^{\pm}(v)=\e_{n-1,n+1}^{\pm}(v)\e_{i,i+1}^{\pm}(u),\qquad&&
\e_{i,i+1}^{\pm}(u)\e_{n-1,n+1}^{\mp}(v)=\e_{n-1,n+1}^{\mp}(v)\e_{i,i+1}^{\pm}(u),&\\
&\f_{i+1,i}^{\pm}(u)\f_{n+1,n-1}^{\pm}(v)=\f_{n+1,n-1}^{\pm}(v)\f_{i+1,i}^{\pm}(u),\qquad&&
\f_{i+1,i}^{\pm}(u)\f_{n+1,n-1}^{\mp}(v)=\f_{n+1,n-1}^{\mp}(v)\f_{i+1,i}^{\pm}(u)&
\end{alignat*}
for $i=1,\dots,n-3$ in type $D$.
\ele

Now we consider the cases $B$ and $D$ separately.

\ble\label{lem:en-1enB}The following relations hold in the algebra $U\big(\overline{R}^{\tss[n]}\big)$ of type~$B$:
\begin{gather*}
\big(qu_{\mp}-q^{-1}v_{\pm}\big)\e_{n-1,n}^{\pm}(u)\e_{n,n+1}^{\mp}(v)
=(u_{\mp}-v_{\pm})\e_{n,n+1}^{\mp}(v)\e_{n-1,n}^{\pm}(u)\\
\qquad{}+\big(q-q^{-1}\big)v_{\pm}\e_{n-1,n+1}^{\pm}(u) -\big(q-q^{-1}\big)u_{\mp}\e_{n-1,n}^{\mp}(v)\e_{n,n+1}^{\mp}(v)\\
\qquad{} -\big(q-q^{-1}\big) u_{\mp}\e_{n-1,n+1}^{\mp}(v),\\
\big(qu-q^{-1}v\big)\e_{n-1,n}^{\pm}(u)\e_{n,n+1}^{\pm}(v)
=(u-v)\e_{n,n+1}^{\pm}(v)\e_{n-1,n}^{\pm}(u)+\big(q-q^{-1}\big)v\e_{n-1,n+1}^{\pm}(u)\\
\qquad{} -\big(q-q^{-1}\big)u\e_{n-1,n}^{\pm}(v)\e_{n,n+1}^{\pm}(v)
-\big(q-q^{-1}\big)u\e_{n-1,n+1}^{\pm}(v),
\end{gather*}
and
\begin{gather*}
(u_{\pm}-v_{\mp})\f_{n,n-1}^{\pm}(u)\f_{n+1,n}^{\mp}(v) =
\big(qu_{\pm}-q^{-1}v_{\mp}\big)\f_{n+1,n}^{\mp}(v)\f_{n,n-1}^{\pm}(u) +\big(q-q^{-1}\big)v_{\mp}\f_{n+1,n-1}^{\mp}(v)\\
\qquad{} -\big(q-q^{-1}v_{\mp}\big)\f_{n+1,n}^{\mp}(v)\f_{n,n-1}^{\mp}(v) -\big(q-q^{-1}\big)u_{\pm}\f_{n+1,n-1}^{\pm}(u),\\
(u-v)\f_{n,n-1}^{\pm}(u)\f_{n+1,n}^{\pm}(v) =
\big(qu-q^{-1}v\big)\f_{n+1,n}^{\pm}(v)\f_{n,n-1}^{\pm}(u)+\big(q-q^{-1}\big)v\f_{n+1,n-1}^{\pm}(v)\\
\qquad{} -\big(q-q^{-1}\big)v\f_{n+1,n}^{\pm}(v)\f_{n,n-1}^{\pm}(v)-\big(q-q^{-1}\big)u\f_{n+1,n-1}^{\pm}(u).
\end{gather*}
\ele

\begin{proof}
We will only prove the first relation.
By \eqref{ELMPj<l} we have
\begin{gather}
\e_{n-1,n}^{\pm}(u)\h^{\mp}_n(v)\e_{n,n+1}^{\mp}(v)-\h^{\mp}_n(v)\e_{n,n+1}^{\mp}(v)\e_{n-1,n}^{\pm}(u)\nonumber\\
\qquad{}=\frac{\big(q-q^{-1}\big)v_{\pm}}{u_{\mp}-v_{\pm}}\h^{\mp}_n(v)\e_{n-1,n+1}^{\pm}(u)
-\frac{\big(q-q^{-1}\big)u_{\mp}}{u_{\mp}-v_{\pm}}\h^{\mp}_n(v)\e_{n-1,n+1}^{\mp}(v).\label{en-1nlnn+1}
\end{gather}
Relation \eqref{ELMPj=l} implies
\begin{gather*}
 \e_{n-1,n}^{\pm}(u)\h^{\mp}_n(v)=
\frac{qu_{\mp}-q^{-1}v_{\pm}}{u_{\mp}-v_{\pm}}\h^{\mp}_n(v)\e_{n-1,n}^{\pm}(u)
+\frac{\big(q-q^{-1}\big)u_{\mp}}{u_{\mp}-v_{\pm}}\h^{\mp}_n(v)\e_{n-1,n}^{\mp}(v)
\end{gather*}
so that \eqref{en-1nlnn+1} can be rewritten as
\begin{gather*}
\frac{qu_{\mp}-q^{-1}v_{\pm}}{u_{\mp}-v_{\pm}}\h^{\mp}_n(v)\e_{n-1,n}^{\pm}(u)\e_{n,n+1}^{\mp}(v)
+\frac{\big(q-q^{-1}\big)u_{\mp}}{u_{\mp}-v_{\pm}}\h^{\mp}_n(v)\e_{n-1,n}^{\mp}(v)\e_{n,n+1}^{\mp}(v)\\
\qquad\quad{}-\h^{\mp}_n(v)\e_{n,n+1}^{\mp}(v)\e_{n-1,n}^{\pm}(u)\\
\qquad{} {}=\frac{\big(q-q^{-1}\big)v_{\pm}}{u_{\mp}-v_{\pm}}\h^{\mp}_n(v)\e_{n-1,n+1}^{\pm}(u)
-\frac{\big(q-q^{-1}\big)u_{\mp}}{u_{\mp}-v_{\pm}}\h^{\mp}_n(v)\e_{n-1,n+1}^{\mp}(v).
\end{gather*}
Since $\h^{\mp}_n(v)$ is invertible, this gives the first relation.
\end{proof}

A similar argument proves the counterpart of Lemma~\ref{lem:en-1enB} for type~$D$.

\ble\label{lem:en-1enD}The following relations hold in the algebra $U(\overline{R}^{\tss[n]})$ of type~$D$:
\begin{gather*}
\big(qu_{\mp}-q^{-1}v_{\pm}\big) \e_{n-2,n-1}^{\pm}(u)\e_{n-1,n+1}^{\mp}(v)\\
\qquad{} =(u_{\mp}-v_{\pm})\e_{n-1,n+1}^{\mp}(v)\e_{n-2,n-1}^{\pm}(u)
+\big(q-q^{-1}\big)v_{\pm}\e_{n-2,n+1}^{\pm}(u)\\
\qquad\quad{} -\big(q-q^{-1}\big)u_{\mp}\e_{n-2,n-1}^{\mp}(v)\e_{n-1,n+1}^{\mp}(v)
-\big(q-q^{-1}\big)u_{\mp}\e_{n-2,n+1}^{\mp}(v),
\\
\big(qu-q^{-1}v\big)\e_{n-2,n-1}^{\pm}(u) \e_{n-1,n+1}^{\pm}(v)\\
\qquad{} =(u-v)\e_{n-1,n+1}^{\pm}(v)\e_{n-2,n-1}^{\pm}(u)+\big(q-q^{-1}\big)v\e_{n-2,n+1}^{\pm}(u)\\
\qquad\quad{} -\big(q-q^{-1}\big)u\e_{n-2,n-1}^{\pm}(v)\e_{n-1,n+1}^{\pm}(v)
-\big(q-q^{-1}\big)u\e_{n-2,n+1}^{\pm}(v),
\end{gather*}
and
\begin{gather*}
(u_{\pm}-v_{\mp})\f_{n-1,n-2}^{\pm}(u) \f_{n+1,n-1}^{\mp}(v)\\
\qquad{} =\big(qu_{\pm}-q^{-1}v_{\mp}\big)\f_{n+1,n-1}^{\mp}(v)\f_{n-1,n-2}^{\pm}(u) +\big(q-q^{-1}\big)v_{\mp}\f_{n+1,n-2}^{\mp}(v)\\
\qquad\quad{} -\big(q-q^{-1}v_{\mp}\big)\f_{n+1,n-1}^{\mp}(v)\f_{n-1,n-2}^{\mp}(v) -\big(q-q^{-1}\big)u_{\pm}\f_{n+1,n-2}^{\pm}(u),
\\
(u-v)\f_{n-1,n-2}^{\pm}(u) \f_{n+1,n-1}^{\pm}(v)\\
 \qquad{} =\big(qu-q^{-1}v\big)\f_{n+1,n-1}^{\pm}(v)\f_{n-1,n-2}^{\pm}(u)+\big(q-q^{-1}\big)v\f_{n+1,n-2}^{\pm}(v)\\
\qquad\quad{} -\big(q-q^{-1}\big)v\f_{n+1,n-1}^{\pm}(v)\f_{n-1,n-2}^{\pm}(v)-\big(q-q^{-1}\big)u\f_{n+1,n-2}^{\pm}(u).
\end{gather*}
\ele

The next lemma is verified by a similar argument
with the use of Corollary~\ref{cor:commu}.

\ble\label{lem:eifnB}
In the algebra $U\big(\overline{R}^{\tss[n]}\big)$ for all $i=1,\dots,n-1$ we have
\begin{alignat*}{3}
& \e_{i,i+1}^{\pm}(u)\f_{n+1,n}^{\pm}(v) =\f_{n+1,n}^{\pm}(v)\e_{i,i+1}^{\pm}(u),
\qquad&&
\e_{i,i+1}^{\pm}(u)\f_{n+1,n}^{\mp}(v) =\f_{n+1,n}^{\mp}(v)\e_{i,i+1}^{\pm}(u),&
\nonumber\\
&\f_{i+1,i}^{\pm}(u)\e_{n,n+1}^{\pm}(v)=\e_{n,n+1}^{\pm}(v)\f_{i+1,i}^{\pm}(u),
\qquad&&
\f_{i+1,i}^{\pm}(u)\e_{n,n+1}^{\mp}(v)=\e_{n,n+1}^{\mp}(v)\f_{i+1,i}^{\pm}(u),&
\end{alignat*}
for type $B$, and for all $i=1,\dots,n-2$ we have
\begin{alignat*}{3}
&\e_{i,i+1}^{\pm}(u)\f_{n+1,n-1}^{\pm}(v) =\f_{n+1,n-1}^{\pm}(v)\e_{i,i+1}^{\pm}(u),
\qquad &&
\e_{i,i+1}^{\pm}(u)\f_{n+1,n-1}^{\mp}(v)=\f_{n+1,n-1}^{\mp}(v)\e_{i,i+1}^{\pm}(u),&\\
&\f_{i+1,i}^{\pm}(u)\e_{n-1,n+1}^{\pm}(v)=\e_{n-1,n+1}^{\pm}(v)\f_{i+1,i}^{\pm}(u),
\qquad&&
\f_{i+1,i}^{\pm}(u)\e_{n-1,n+1}^{\mp}(v)=\e_{n-1,n+1}^{\mp}(v)\f_{i+1,i}^{\pm}(u),&
\end{alignat*}
for type $D$.
\ele

We are now in a position to summarise the results of Sections~\ref{subsec:typea}, \ref{subsec:lowrankB} and \ref{subsec:lowrankD} and give complete lists of relations between the Gaussian generators. The completeness of the relations will be established in Section~\ref{sec:urm}.

\bth\label{thm:relrbar}\quad
\begin{enumerate}\itemsep=0pt
\item[{\rm (i)}] The following relations hold in the algebra $U\big(\overline{R}^{\tss[n]}\big)$ of type~$B$. For the relations invol\-ving~$\h^{\pm}_i(u)$ we have
\begin{gather*}
\h^{\pm}_i(u)\h^{\pm}_j(v) =\h^{\pm}_j(v)\h^{\pm}_i(u),\qquad
\h^{\pm}_i(u)\h^{\mp}_i(v) =\h^{\mp}_i(v)\h^{\pm}_i(u),\qquad i=1,\dots, n,\\
\frac{u_{\pm}-v_{\mp}}{qu_{\pm}-q^{-1}v_{\mp}}\h^{\pm}_i(u)\h^{\mp}_j(v) =
\frac{u_{\mp}-v_{\pm}}{qu_{\mp}-q^{-1}v_{\pm}}\h^{\mp}_j(v)\h^{\pm}_i(u),\qquad i<j,
\end{gather*}
and
\begin{gather*}
\frac{q^{-1}u_{\pm}-qv_{\mp}}{qu_{\pm}-q^{-1}v_{\mp}}\ts
\frac{q^{1/2}u_{\pm}-q^{-1/2}v_{\mp}}{q^{-1/2}u_{\pm}-q^{1/2}v_{\mp}}\ts
\h^{\pm}_{n+1}(u)\h^{\mp}_{n+1}(v)\\
\qquad{} =\frac{q^{-1}u_{\mp}-qv_{\pm}}{qu_{\mp}-q^{-1}v_{\pm}}\ts
\frac{q^{1/2}u_{\mp}-q^{-1/2}v_{\pm}}{q^{-1/2}u_{\mp}-q^{1/2}v_{\pm}}\ts
\h^{\mp}_{n+1}(v)\h^{\pm}_{n+1}(u).
\end{gather*}
The relations involving $\h^{\pm}_i(u)$ and $\Xc_{j}^{\pm}(v)$ are
\begin{gather*}
\h_{i}^{\pm}(u)\Xc_{j}^{+}(v)
=\frac{u-v_{\pm}}{q^{(\ep_i,\alpha_j)}u-q^{-(\ep_i,\alpha_j)}v_{\pm}}
\Xc_{j}^{+}(v)\h_{i}^{\pm}(u),\\
\h_{i}^{\pm}(u)\Xc_{j}^{-}(v) =\frac{q^{-(\ep_i,\alpha_j)}u_{\pm}-q^{(\ep_i,\alpha_j)}v}{u_{\pm}-v}
\Xc_{j}^{-}(v)\h_{i}^{\pm}(u)
\end{gather*}
for $i\neq n+1$, together with
\begin{gather*}
\h_{n+1}^{\pm}(u)\Xc_n^{+}(v) =
\frac{(qu_{\mp}-v)(u_{\mp}-v)}{(u_{\mp}-qv)\big(qu_{\mp}-q^{-1}v\big)}\Xc_n^{+}(v)\h_{n+1}^{\pm}(u),\\
\h_{n+1}^{\pm}(u)\Xc_n^{-}(v) =
\frac{(u_{\pm}-qv)\big(qu_{\pm}-q^{-1}v\big)}{(qu_{\pm}-v)(u_{\pm}-v)}\Xc_n^{-}(v)\h_{n+1}^{\pm}(u),
\end{gather*}
and
\begin{gather*}
\h_{n+1}^{\pm}(u)\Xc_{i}^{+}(v) =\Xc_{i}^{+}(v)\h_{n+1}^{\pm}(u),\qquad
\h_{n+1}^{\pm}(u)\Xc_{i}^{-}(v)=\Xc_{i}^{-}(v)\h_{n+1}^{\pm}(u),
\end{gather*}
for $1\leqslant i\leqslant n-1$. For the relations involving $\Xc^{\pm}_i(u)$ we have
\begin{gather*}
\big(u-q^{\pm (\alpha_i,\alpha_j)}v\big)\Xc_{i}^{\pm}\big(uq^{i}\big)\Xc_{j}^{\pm}\big(vq^{j}\big)
=\big(q^{\pm (\alpha_i,\alpha_j)}u-v\big) \Xc_{j}^{\pm}\big(vq^{j}\big)\Xc_{i}^{\pm}\big(uq^{i}\big)
\end{gather*}
for $i,j=1,\dots,n$; and
\begin{gather*}
\big[\Xc_i^{+}(u),\Xc_j^{-}(v)\big]=
\delta_{ij}\big(q-q^{-1}\big)\\
\hphantom{\big[\Xc_i^{+}(u),\Xc_j^{-}(v)\big]=}{}\times \big(\delta\big(u\ts q^{-c}/v\big)\h_i^{-}(v_+)^{-1}\h_{i+1}^{-}(v_+)
-\delta\big(u\ts q^{c}/v\big)\h_i^{+}(u_+)^{-1}\h_{i+1}^{+}(u_+)\big)
\end{gather*}
together with the Serre relations
\begin{gather}
\sum_{\pi\in \Sym_{r}}\sum_{l=0}^{r}(-1)^l{{r}\brack{l}}_{q_i}\!\!
 \Xc^{\pm}_{i}(u_{\pi(1)})\cdots \Xc^{\pm}_{i}(u_{\pi(l)})
 \Xc^{\pm}_{j}(v)\tss \Xc^{\pm}_{i}(u_{\pi(l+1)})\cdots \Xc^{\pm}_{i}(u_{\pi(r)})=0,\!\!\!\!\!\label{serrex}
\end{gather}
which hold for all $i\neq j$ and we set $r=1-A_{ij}$.

\item[{\rm(ii)}] The following relations hold in the algebra
$U\big(\overline{R}^{\tss[n]}\big)$ of type $D$.

For the relations involving $\h^{\pm}_i(u)$ we have
\begin{gather*}
\h^{\pm}_i(u)\h^{\pm}_j(v) =\h^{\pm}_j(v)\h^{\pm}_i(u),\qquad
\h^{\pm}_i(u)\h^{\mp}_i(v) =\h^{\mp}_i(v)\h^{\pm}_i(u),\qquad i=1,\dots, n+1,\\[0.4em]
\frac{u_{\pm}-v_{\mp}}{qu_{\pm}-q^{-1}v_{\mp}}\h^{\pm}_i(u)\h^{\mp}_j(v) =
\frac{u_{\mp}-v_{\pm}}{qu_{\mp}-q^{-1}v_{\pm}}\h^{\mp}_j(v)\h^{\pm}_i(u)
\end{gather*}
for $i<j$ with $(i,j)\neq (n,n+1)$, and
\begin{gather*}
\frac{q^{-1}u_{\pm}-qv_{\mp}}{qu_{\pm}-q^{-1}v_{\mp}}
\frac{u_{\pm}-v_{\mp}}{u_{\pm}-q^{-1}v_{\mp}}
\h^{\pm}_{n}(u)\h^{\mp}_{n+1}(v)=
\frac{q^{-1}u_{\mp}-qv_{\pm}}{qu_{\mp}-q^{-1}v_{\pm}}
\frac{u_{\mp}-v_{\pm}}{u_{\mp}-q^{-1}v_{\pm}}
\h^{\mp}_{n+1}(v)\h^{\pm}_{n}(u).
\end{gather*}
The relations
involving $\h^{\pm}_i(u)$ and $\Xc_{j}^{\pm}(v)$ are
\begin{gather*}
\h_{i}^{\pm}(u)\Xc_{j}^{+}(v)
=\frac{u-v_{\pm}}{q^{(\ep_i,\alpha_j)}u-q^{-(\ep_i,\alpha_j)}v_{\pm}}
\Xc_{j}^{+}(v)\h_{i}^{\pm}(u),\\
\h_{i}^{\pm}(u)\Xc_{j}^{-}(v)
=\frac{q^{(\ep_i,\alpha_j)}u_{\pm}-q^{-(\ep_i,\alpha_j)}v}{u_{\pm}-v}
\Xc_{j}^{-}(v)\h_{i}^{\pm}(u)
\end{gather*}
for $i\neq n+1$, together with
\begin{gather*}
\h_{n+1}^{\pm}(u)\Xc_n^{+}(v) =
\frac{u_{\mp}-v}{q^{-1}u_{\mp}-qv}\Xc_n^{+}(v)\h_{n+1}^{\pm}(u),\\
\h_{n+1}^{\pm}(u)\Xc_n^{-}(v) =
\frac{q^{-1}u_{\pm}-qv}{u_{\pm}-v}\Xc_n^{-}(v)\h_{n+1}^{\pm}(u),
\end{gather*}
and
\begin{gather*}
\h_{n+1}^{\pm}(u)\Xc_{n-1}^{+}(v) =
\frac{u_{\mp}-v}{qu_{\mp}-q^{-1}v}\Xc_{n-1}^{+}(v)\h_{n+1}^{\pm}(u),\\
\h_{n+1}^{\pm}(u)\Xc_{n-1}^{-}(v) =
\frac{qu_{\pm}-q^{-1}v}{u_{\pm}-v}\Xc_{n-1}^{-}(v)\h_{n+1}^{\pm}(u),
\end{gather*}
while
\begin{gather*}
\h_{n+1}^{\pm}(u)\Xc_{i}^{+}(v)
 =\Xc_{i}^{+}(v)\h_{n+1}^{\pm}(u),\qquad
\h_{n+1}^{\pm}(u)\Xc_{i}^{-}(v)
 =\Xc_{i}^{-}(v)\h_{n+1}^{\pm}(u),
\end{gather*}
for $1\leqslant i\leqslant n-2$. For the relations involving $\Xc^{\pm}_i(u)$ we have
\begin{gather*}
\big(u-q^{\pm (\alpha_i,\alpha_j)}v\big)\Xc_{i}^{\pm}\big(uq^{i}\big)\Xc_{j}^{\pm}\big(vq^{j}\big)
=\big(q^{\pm (\alpha_i,\alpha_j)}u-v\big) \Xc_{j}^{\pm}\big(vq^{j}\big)\Xc_{i}^{\pm}\big(uq^{i}\big)
\end{gather*}
for $i,j=1,\dots,n-1$;
\begin{gather*}
\big(u-q^{\pm (\alpha_i,\alpha_n)}v\big)\Xc_{i}^{\pm}\big(uq^{i}\big)\Xc_{n}^{\pm}\big(vq^{n-1}\big)
=\big(q^{\pm (\alpha_i,\alpha_n)}u-v\big) \Xc_{n}^{\pm}\big(vq^{n-1}\big)\Xc_{i}^{\pm}\big(uq^{i}\big)
\end{gather*}
for $i=1,\dots,n-1$;
\begin{gather*}
\big(u-q^{\pm (\alpha_n,\alpha_n)}v\big)\Xc_{n}^{\pm}(u)\Xc_{n}^{\pm}(v)
=\big(q^{\pm (\alpha_n,\alpha_n)}u-v\big) \Xc_{n}^{\pm}(v)\Xc_{n}^{\pm}(u)
\end{gather*}
and
\begin{gather*}
\big[\Xc_i^{+}(u),\Xc_j^{-}(v)\big]=
\delta_{ij}\big(q-q^{-1}\big)\\
\hphantom{\big[\Xc_i^{+}(u),\Xc_j^{-}(v)\big]=}{}\times \big(\delta\big(u\ts q^{-c}/v\big)\h_i^{-}(v_+)^{-1}\h_{i+1}^{-}(v_+)
-\delta\big(u\ts q^{c}/v\big)\h_i^{+}(u_+)^{-1}\h_{i+1}^{+}(u_+)\big)
\end{gather*}
together with the Serre relations
\begin{gather} \label{serrexD}
\sum_{\pi\in \Sym_{r}}\sum_{l=0}^{r}(-1)^l{{r}\brack{l}}_{q_i}\!\!
 \Xc^{\pm}_{i}(u_{\pi(1)})\cdots \Xc^{\pm}_{i}(u_{\pi(l)})
 \Xc^{\pm}_{j}(v)\tss \Xc^{\pm}_{i}(u_{\pi(l+1)})\cdots \Xc^{\pm}_{i}(u_{\pi(r)})=0,\!\!\!\!
\end{gather}
which hold for all $i\neq j$ and we set $r=1-A_{ij}$.
\end{enumerate}
\eth

\begin{proof}
All relations except for \eqref{serrex} and \eqref{serrexD} follow from
the corresponding results in Sections~\ref{subsec:typea},
\ref{subsec:lowrankB} and \ref{subsec:lowrankD} by applying
Theorem~\ref{thm:red} and Proposition~\ref{prop:gauss-consist} and
recalling the definition~\eqref{Xibar}. The remaining Serre
relations are verified in the same way as for type~$C$ \cite[Section~4.6]{jlm:ib-c}
by adapting the Levendorski argument~\cite{l:gd} to the quantum affine algebras.
\end{proof}

By using Theorem~\ref{thm:relrbar} and Proposition~\ref{prop:corrgauss} we arrive at the following homomorphism theorem for the extended quantum affine algebra $U^{\ext}_q(\wh{\oa}_{N})$ introduced in Definition~\ref{def:eqaa}.

\bth\label{thm:homom}
The mapping
\begin{alignat*}{3}
& X^{+}_i(u) \mapsto X^{+}_i(u), \qquad&&\text{for}\quad i=1,\dots,n,& \\
& X^{-}_i(u) \mapsto X^{-}_i(u), \qquad&&\text{for}\quad i=1,\dots,n,& \\
& h^{\pm}_j(u)\mapsto h^{\pm}_j(u), \qquad&&\text{for}\quad j=1,\dots,n+1,&
\end{alignat*}
defines a homomorphism $DR\colon  U^{\ext}_q(\wh{\oa}_{N})\rightarrow U(R)$, where $X_i^{\pm}(u)$ on the right hand side is given by~\eqref{Xi}, \eqref{Xnbarp} and~\eqref{Xnbarm}.
\eth

We will show in the next section that the homomorphism $DR$ provided by Theorem~\ref{thm:homom}
is an isomorphism by constructing the inverse map with the use of
the universal $R$-matrix for the algebra $U_q(\wh{\oa}_{N})$ in a way similar to
types $A$ and $C$; see~\cite{fm:ha} and~\cite[Section~5]{jlm:ib-c}.

\section[The universal $R$-matrix and inverse map]{The universal $\boldsymbol{R}$-matrix and inverse map}\label{sec:urm}

We will need explicit formulas for the
universal $R$-matrix for the quantum affine algebras
obtained by Khoroshkin and Tolstoy~\cite{kt:ur}
and Damiani~\cite{d:rm,da:Un}.

Recall that the Cartan matrix for the Lie algebra $\oa_{N}$ is defined in~\eqref{cartan}
and consider the diagonal matrix $C=\diag[r_1,r_2,\dots,r_n]$ with
$r_i=(\al_i,\al_i)/2$. The matrix $B=[B_{ij}]:=CA$ is symmetric with $B_{ij}=(\al_i,\al_j)$.
We will use the notation $\tilde{B}=[\tilde{B}_{ij}]$ for the inverse matrix $B^{-1}$.
We will also need the $q$-deformed matrix $B(q)=[B_{ij}(q)]$ with $B_{ij}(q)=[B_{ij}]_q$
and its inverse $\tilde{B}(q)=[\tilde{B}_{ij}(q)]$; see~\eqref{kq}.
Both $n\times n$ matrices $\tilde{B}$ and $\tilde{B}(q)$ are symmetric and
for $N=2n+1$ (type $B$) we have
\begin{gather*}
\tilde{B}_{ij}=j\qquad\text{for}\quad j\leqslant i
\end{gather*}
and
\begin{gather}\label{Bijqk}
\tilde{B}_{ij}(q)=
\begin{cases}\dfrac{[j]_{q}}{[n]_{q}-[n-1]_{q}} &\text{for}\quad i=n,\vspace{1mm}\\
 \dfrac{[j]_{q}\ts\big([n-i]_{q}-[n-i-1]_{q}\big)}{[n]_{q}-[n-1]_{q}}
 &\text{for}\quad j\leqslant i<n,
\end{cases}
\end{gather}
whereas for $N=2n$ (type $D$) the entries are given by
\begin{gather*}
\tilde{B}_{ij}=\begin{cases}
 j &\text{for}\quad j\leqslant i\leqslant n-2, \\
 j/2 &\text{for}\quad j\leqslant n-2,\quad i=n-1,n, \\
 n/4 &\text{for}\quad i=j\geqslant n-1,\\
 (n-2)/4 &\text{for}\quad i=n,\quad j=n-1,
 \end{cases}
\end{gather*}
and
\begin{gather}\label{BijqktypeD}
\tilde{B}_{ij}(q)=\begin{cases}
 \dfrac{[j]_{q}\ts [2]_{q^{n-1-i}}}{[2]_{q^{n-1}}}
 &\text{for}\quad j\leqslant i\leqslant n-2, \vspace{1mm}\\
 \dfrac{[j]_{q}}{[2]_{q^{n-1}}}
 &\text{for}\quad j\leqslant n-2,\quad i=n-1,n, \vspace{1mm}\\
 \dfrac{[n]_{q}}{[2]_{q}\ts [2]_{q^{n-1}}}
 &\text{for}\quad i=j\geqslant n-1,\vspace{1mm}\\
 \dfrac{[n-2]_{q}}{[2]_{q}\ts [2]_{q^{n-1}}}
 &\text{for}\quad i=n,\quad j=n-1.
 \end{cases}
\end{gather}

As with type $C$
\cite[Section~5]{jlm:ib-c}, we will use the parameter-dependent
universal $R$-matrix defined in terms of the presentation of the quantum affine algebra
used in Section~\ref{subsec:isoDJD}. The formula for the $R$-matrix uses
the $\hbar$-adic settings so we will regard
the algebra over $\CC[[\hbar]]$ and set $q=\exp(\hbar)\in \CC[[\hbar]]$. It is well-known that
the $\mathbb C\big(q^{1/2}\big)$-algebra $U_q(\mathfrak g_
N)$ actually embeds inside the $\CC[[\hbar]]$-algebra $U_{\hbar}(\mathfrak g_N)$ due to the flatness of the latter as a deformation of
$U(\mathfrak g_N)$.
Define elements $h_1,\dots,h_n$ by
setting $k_i=\exp(\hbar h_i)$.
The universal $R$-matrix is given by
\begin{gather}\label{rdec}
\Rc(u)=\Rc^{>0}(u)\Rc^{0}(u)\Rc^{<0}(u),
\end{gather}
where
\begin{gather*}
 \Rc^{>0}(u)  = \prod_{\alpha\in \Delta_+}\prod_{k\geqslant 0}
 \exp_{q_i}\big(\big(q_i^{-1}-q_i\big)u^{k}E_{\alpha+k\delta}\otimes F_{\alpha+k\delta}\big), \\
 \Rc^{<0}(u)  = T^{-1}\prod_{\alpha\in \Delta_+}\prod_{k> 0}
 \exp_{q_i}\big(\big(q_i^{-1}-q_i\big)u^{k}E_{-\alpha+k\delta}\otimes F_{-\alpha+k\delta}\big)\ts T
\end{gather*}
with $T=\exp(-\hbar \tilde{B}_{ij}\tss h_i\otimes h_j)$ and
\begin{gather*}
 \Rc^{0}(u)=\exp\left(
 \sum\limits_{k>0}\sum\limits_{i,j=1}^{n}\frac{\big(q_i^{-1}-q_i\big)\big(q_j^{-1}-q_j\big)}{q^{-1}-q}
 \frac{k}{[k]_q}\tilde{B}_{ij}(q^k)u^k q^{kc/2}a_{i,k}\otimes a_{j,-k}q^{-kc/2}\right)\ts T
\end{gather*}
(see \cite[Definition~4]{d:rm} for the description of the order of the products in $\Rc^{>0}$ and $\Rc^{<0}$).
It satisfies the Yang--Baxter equation in the form
\begin{gather}\label{MYBE}
\Rc_{12}(u)\Rc_{13}\big(uvq^{-c_2}\big)\Rc_{23}(v) =\Rc_{23}(v)\Rc_{13}\big(uvq^{c_2}\big)\Rc_{12}(u)
\end{gather}
where $c_2=1\otimes c\otimes 1$; cf.~\cite{fri:qa}.

A straightforward calculation verifies the following formulas for the vector representation of the quantum affine algebra. As before, we denote by $e_{ij}\in\End\CC^{N}$ the standard matrix units.

\bpr\label{prop:FFRep}
The mappings $q^{\pm c/2}\mapsto 1$,
\begin{gather*}
 x^{+}_{ik}  \mapsto -q^{-ik}e_{i+1,i}+q^{-(2n-1-i)k}e_{i',(i+1)'}, \\
 x^{-}_{ik}  \mapsto -q^{-ik}e_{i,i+1}+q^{-(2n-1-i)k}e_{(i+1)',i'},\\
 a_{ik}  \mapsto \frac{[k]_{q_i}}{k}\big(q^{-ik}\big(q^{-k}e_{i+1,i+1}-q^{k}e_{ii}\big)
 +q^{-(2n-1-i)k}\big(q^{-k}e_{i'i'}-q^ke_{(i+1)'(i+1)'}\big)  \big),\\
 k_i  \mapsto q(e_{i+1,i+1}+e_{i',i'})+q^{-1}(e_{ii}+e_{(i+1)',(i+1)'})
 +\sum\limits_{j\neq i,i+1,i',(i+1)'}e_{jj},
\end{gather*}
 for $i=1,\dots,n-1$, and
\begin{gather*}
 x^{+}_{nk}  \mapsto [2]_{q_n}^{1/2}\big({-}q^{-nk}e_{n+1,n}+q^{-(n-1)k}e_{n',n+1}\big), \\
 x^{-}_{nk}  \mapsto [2]_{q_n}^{1/2}\big({-}q^{-nk}e_{n,n+1}+q^{-(n-1)k}e_{n+1,n'}\big),\\
 a_{nk}  \mapsto \frac{[2k]_{q_n}}{k}\big({-}q^{-(n-1)k}e_{nn}+\big(q^{-nk}-q^{-(n-1)k}\big)e_{n+1,n+1}\big)+q^{-nk}e_{n'n'}
 \big),\\
 k_n  \mapsto qe_{n',n'}+q^{-1}e_{nn}  +\sum\limits_{j\neq n,n'}e_{jj},
\end{gather*}
in type $B$, and the mappings $q^{\pm c/2}\mapsto 1$,
\begin{gather*}
 x^{+}_{ik}  \mapsto -q^{-ik}e_{i+1,i}+q^{-(2n-2-i)k}e_{i',(i+1)'}, \\
 x^{-}_{ik}  \mapsto -q^{-ik}e_{i,i+1}+q^{-(2n-2-i)k}e_{(i+1)',i'},\\
 a_{ik}  \mapsto \frac{[k]_{q_i}}{k}\big(q^{-ik}\big(q^{-k}e_{i+1,i+1}-q^{k}e_{ii}\big)
 +q^{-(2n-2-i)k}\big(q^{-k}e_{i'i'}-q^ke_{(i+1)'(i+1)'}\big) \big),\\
 k_i  \mapsto q(e_{i+1,i+1}+e_{i',i'})+q^{-1}(e_{ii}+e_{(i+1)',(i+1)'})
 +\sum\limits_{j\neq i,i+1,i',(i+1)'}e_{jj},
\end{gather*}
 for $i=1,\dots,n-1$, and
\begin{gather*} x^{+}_{nk}  \mapsto q^{-(n-1)k}(-e_{n+1,n-1}+e_{n+2,n}), \\
 x^{-}_{nk}  \mapsto q^{-(n-1)k}(-e_{n-1,n+1}+e_{n,n+2}),\\
 a_{nk}  \mapsto \frac{[k]_{q_n}}{k}q^{-(n-1)k}\big(q^{-k}e_{n+1,n+1}-q^{k}e_{n-1,n-1}
+q^{-k}e_{n+2,n+2}-q^{k}e_{n,n}
 \big),\\
 k_n  \mapsto q(e_{n+1,n+1}+e_{n+2,n+2})+q^{-1}(e_{n-1,n-1}+e_{n,n})
 +\sum\limits_{j\neq n-1,n,n+1,n+2}e_{jj},
\end{gather*}
in type $D$, define a representation $\pi_V\colon U_{\hbar}(\wh{\oa}_{N})\to\End V$
of the algebra $U_{\hbar}(\wh{\oa}_{N})$ on the vector space $V=\mathbb{C}^{N}[\![\hbar]\!]$.
\epr

It follows from the results of \cite[Theorem~4.2]{fri:qa} that the $R$-matrix defined in \eqref{ru}
coincides with the image of the universal $R$-matrix:
\begin{gather*}
R(u)=(\pi_{V}\otimes \pi_{V})\ts\Rc(u).
\end{gather*}
Introduce the $L$-{\em operators} in $U_q(\wh{\oa}_{N})$ by the formulas
\begin{gather*}
 \tilde{L}^{+}(u)  = (\id\otimes \pi_{V})\ts \Rc_{21}\big(uq^{c/2}\big),\qquad
 \tilde{L}^{-}(u)  = (\id\otimes \pi_{V})\ts \Rc_{12}\big(u^{-1}q^{-c/2}\big)^{-1}.
\end{gather*}
Recall the series $z^{\pm}(u)$ defined in \eqref{zpm}. Their coefficients
are central in the algebra $U^{\ext}_q(\wh{\oa}_{N})$; see Proposition~\ref{prop:zu}.
Therefore, the Yang--Baxter equation~\eqref{MYBE} implies the relations for the $L$-operators:
\begin{gather*}
R(u/v)L^{\pm}_1(u)L^{\pm}_2(v)  = L^{\pm}_2(v)L^{\pm}_1(u)R(u/v), \\
 R(u_{+}/v_{-})L^{+}_1(u)L^{-}_2(v)  = L^{-}_2(v)L^{+}_1(u)R(u_{-}/v_{+}),
\end{gather*}
where we set
\begin{gather*}
L^{+}(u) = \tilde{L}^{+}(u)\prod_{m=0}^{\infty}
z^{+}\big(u\tss\xi^{-2m-1}\big)\tss z^{+}\big(u\tss\xi^{-2m-2}\big)^{-1},\\
L^{-}(u) = \tilde{L}^{-}(u)\prod_{m=0}^{\infty}
z^{-}\big(u\tss\xi^{-2m-1}\big)\tss z^{-}\big(u\tss\xi^{-2m-2}\big)^{-1}.
\end{gather*}
Note that although these formulas for the entries of the matrices $L^{\pm}(u)$
involve a completion of the center of the algebra $U^{\ext}_q(\wh{\oa}_{N})$,
it will turn out that the coefficients of the series in $u^{\pm1}$ actually
belong to~$U^{\ext}_q(\wh{\oa}_{N})$; see the proof of
Proposition~\ref{prop:Hp} below. Thus, we may conclude that the mapping
\begin{gather}\label{inversemap}
RD\colon \ L^{\pm}(u)\mapsto L^{\pm}(u)
\end{gather}
defines a homomorphism $RD$ from the algebra $U(R)$
to a completed algebra $U^{\ext}_q(\wh{\oa}_{N})$,
where we use the same notation for the corresponding elements of the algebras.

By using the vector representation $\pi_{V}$ defined in Proposition~\ref{prop:FFRep},
introduce the matrices $F^{\pm}(u)$, $E^{\pm}(u)$ and $H^{\pm}(u)$ by setting
\begin{gather*}
F^{+}(u) =(\id\otimes \pi_{V})\ts \Rc^{>0}_{21}\big(uq^{c/2}\big),\qquad
E^{+}(u) =(\id\otimes \pi_{V})\ts \Rc^{<0}_{21}\big(uq^{c/2}\big),\\
H^{+}(u) =(\id\otimes \pi_{V})
\Rc^{0}_{21}(uq^{c/2})\ts\prod_{m=0}^{\infty} z^{+}\big(u\xi^{-2m-1}\big)z^{+}\big(u\xi^{-2m-2}\big)^{-1},
\end{gather*}
and
\begin{gather*}
E^{-}(u)=(\id\otimes \pi_V)\Rc^{>0}\big(u_+^{-1}\big)^{-1},\qquad
F^{-}(u)=(\id\otimes \pi_V)\Rc^{<0}\big(u_+^{-1}\big)^{-1},\\
H^{-}(u)=(\id\otimes \pi_V)\big(\Rc^{0}\big(u_+^{-1}\big)\big)^{-1}
\prod_{m=0}^{\infty}z^{-}\big(u\tss\xi^{2m-1}\big)^{-1}z^{-}\big(u\tss\xi^{2m-2}\big).
\end{gather*}
The decomposition \eqref{rdec} implies the corresponding decomposition for the matrix $L^{\pm}(u)$:
\begin{gather*}
L^{\pm}(u)=F^{\pm}(u)H^{\pm}(u)E^{\pm}(u).
\end{gather*}

Recall the Drinfeld generators $x^{\pm}_{i,k}$ of the algebra $U_q(\wh{\oa}_{N})$,
as defined in the Introduction, and combine them into the formal series
\begin{alignat*}{3}
& x_i^{-}(u)^{\geqslant 0}=\sum\limits _{k\geqslant 0}x^{-}_{i,-k}u^{k},
\qquad&&  x^{+}_{i}(u)^{> 0}=\sum\limits _{k> 0}x^{+}_{i,-k}u^{k},&\\
& x_i^{-}(u)^{< 0}=\sum\limits _{k> 0}x^{-}_{i,k}u^{-k},
\qquad &&  x^{+}_{i}(u)^{\leqslant 0}=\sum\limits _{k\geqslant 0}x^{+}_{i,k}u^{-k}.&
\end{alignat*}
Furthermore, for all $i=1,\dots, n-1$ set
\begin{alignat*}{3}
&f_i^{+}(u)=\big(q_i -q_i^{-1}\big)x_i^{-}\big(u_+q^{-i}\big)^{\geqslant 0},\qquad&&
e_i^{+}(u)=\big(q_i -q_i^{-1}\big)x^{+}_{i}\big(u_{-}q^{-i}\big)^{> 0},& \\
& f_{i}^{-}(u)=\big(q_i^{-1}-q_i\big)x^{-}_{i}\big(u_{-}q^{-i}\big)^{<0},\qquad &&
e_i^{-}(u)=\big(q_i^{-1}-q_i\big)x^{+}_{i}\big(u_{+}q^{-i}\big)^{\leqslant 0},&
\end{alignat*}
whereas
\begin{alignat*}{3}
& f_n^{+}(u)=\big(q_n -q_n^{-1}\big)\ts [2]_{q_n}^{1/2}\ts x_n^{-}\big(u_+q^{-n}\big)^{\geqslant 0},\qquad&&
e_n^{+}(u)=\big(q_n -q_n^{-1}\big)\ts [2]_{q_n}^{1/2}\ts x_n^{+}\big(u_{-}q^{-n}\big)^{> 0},& \\
& f_n^{-}(u)=\big(q_n^{-1}-q_n\big)\ts [2]_{q_n}^{1/2}\ts x^{-}_{n}\big(u_{-}q^{-n}\big)^{<0},\qquad&&
e_n^{-}(u)=\big(q_n^{-1}-q_n\big)\ts [2]_{q_n}^{1/2}\ts x_n^{+}\big(u_{+}q^{-n}\big)^{\leqslant 0}.&
\end{alignat*}
for type $B$, and
\begin{alignat*}{3}
& f_n^{+}(u)=\big(q_n -q_n^{-1}\big)x_n^{-}\big(u_+q^{-(n-1)}\big)^{\geqslant 0},\qquad&&
e_n^{+}(u)=\big(q_n -q_n^{-1}\big)x_n^{+}\big(u_{-}q^{-(n-1)}\big)^{> 0},& \\
& f_n^{-}(u)=\big(q_n^{-1}-q_n\big)x^{-}_{n}\big(u_{-}q^{-(n-1)}\big)^{<0},\qquad&&
e_n^{-}(u)=\big(q_n^{-1}-q_n\big)x_n^{+}\big(u_{+}q^{-(n-1)}\big)^{\leqslant 0}.&
\end{alignat*}
for type $D$.

\bpr\label{prop:Fp}
The matrix $F^{\pm}(u)$ is lower unitriangular and has the form
\begin{gather*}
F^{\pm}(u)=
\begin{bmatrix}
 1 & & & & & & & \\
 f_{1}^{\pm}(u) & \qquad 1 & & & & &\bigcirc & \\
 & \qquad \ddots& & \ddots& & & & \\
 & & & f^{\pm}_n(u)& 1& & & \\
 & & & & -f^{\pm}_{n-1}\big(u\xi q^{2(n-1)}\big) & 1 & & \\
 & \qquad \bigstar & & & & \ddots& \ddots & \\
 & & & & & & -f_{1}^{\pm}\big(u \xi q^{2}\big)& 1
\end{bmatrix}
\end{gather*}
for type $B$, and
\begin{gather*}
F^{\pm}(u)=
\begin{bmatrix}
 1 & & & & & & & \\
 f_{1}^{\pm}(u) & 1 & & & & &\bigcirc & \\
 & \ddots& \ddots& & & & & \\
 & & f_{n-1}^{\pm}(u) & 1 & & & & \\
 & & f_{n}^{\pm}(u) & 0 & 1 & & & \\
 & & & -f_{n}^{\pm}(u) & -f_{n-1}^{\pm}\big(u\xi q^{2(n-1)}\big) & 1 & & \\
 & & \bigstar & & & \ddots & \ddots & \\
 & & & & & & -f_{1}^{\pm}\big(u\xi q^{2}\big)& 1
\end{bmatrix}
\end{gather*}
for type $D$.
\epr

\begin{proof}The argument is a straightforward verification relying on the formulas
of Proposition~\ref{prop:FFRep}; cf.~\cite[Proposition~5.2]{jlm:ib-c}.
\end{proof}

As in Section~\ref{subsec:eqaa}, we will assume that the algebra $U_q(\wh{\oa}_{2n})$
is extended by adjoining the square roots $(k_{n-1}k_n)^{\pm 1/2}$
(no extension is necessary in type $B$).

\ble\label{lem:K} The image $(\id\ot \pi_V)\big(T_{21}\big)$ is the diagonal matrix
\begin{gather*}
\diag\left[\prod_{b=1}^{n}k_b,\tss \prod_{b=2}^{n}k_b, \tss \dots, \tss
\prod_{b=i}^{n}k_b, \tss\dots, \tss k_n, \tss 1,
\tss k_n^{-1},\tss \dots, \prod_{b=1}^{n}k_b^{-1}\right]
\end{gather*}
for type $B$, and
\begin{gather*}
\diag\left[\prod_{b=1}^{n-2}k_b(k_{n-1}k_n)^{1/2},
\prod_{b=2}^{n-2}k_b(k_{n-1}k_n)^{1/2}, \tss\dots,
\tss (k_{n-1}k_n)^{1/2}, \tss (k_{n-1}^{-1}k_n)^{1/2},\right.\\
\left.\qquad{} \tss (k_{n-1}^{-1}k_n)^{-1/2},\tss (k_{n-1}k_n)^{-1/2}, \dots,
\prod_{b=1}^{n-2}k_b^{-1}(k_{n-1}k_n)^{-1/2}\right]
\end{gather*}
for type $D$.
\ele

\begin{proof}
The calculation is the same as in type $C$; see \cite[Lemma~5.3]{jlm:ib-c}.
\end{proof}

\bpr\label{prop:Ep}
The matrix $E^{\pm}(u)$ is upper unitriangular and has the form
\begin{gather*}
E^{\pm}(u)=
\begin{bmatrix}
 1 & \!e_1^{\pm}(u)\! & & & & & & & \\
 & 1 & \!e_2^{\pm}(u)\! & & & & \bigstar& &\\
 & & \ddots & \ddots & & & & & \\
 & & & 1 & \!e_n^{\pm}(u)\! & & & & \\
 & & & & 1 & \!-e_{n-1}^{\pm}\big(u\xi q^{2n-2)}\big)\! & & & \\
 & & & & & \ddots & \ddots & &\\
 & & & & & & 1 & \!-e_2^{\pm}\big(u\xi q^4\big)\!&\\
 & & \bigcirc & & & & & 1& \! -e_1^{\pm}\big(u\xi q^2\big)\\
 & & & & & & & & 1
\end{bmatrix}
\end{gather*}
for type $B$, and
\begin{gather*}
E^{\pm}(u)=
\begin{bmatrix}
 1 & e_1^{\pm}(u) & & & & & & & \\
 & \ddots & \ddots & & & & & \bigstar & \\
 & & 1 & e_{n-1}^{\pm}(u) & e_n^{\pm}(u) & & & & \\
 & & & 1 & 0 & -e_n^{\pm}(u) & & & \\
 & & & & 1 & -e_{n-1}^{\pm}(u) & & & \\
 & & & & & \ddots & \ddots & & \\
 & &\bigcirc & & & & 1 & -e_2^{\pm}\big(u\xi q^4\big) & \\
 & & & & & & & 1 & -e_1^{\pm}\big(u\xi q^2\big) \\
 & & & & & & & & 1
\end{bmatrix}
\end{gather*}
for type $D$.
\epr

\begin{proof}
By the construction of the root vectors $E_{-\alpha+k\delta}$ and
the formulas for the representation~$\pi_V$ provided by Proposition~\ref{prop:FFRep},
it is sufficient to evaluate the image of the product
\begin{gather}\label{epleva}
T_{21}^{-1}\prod_{k> 0}
 \exp_{q_i}\big(\big(q_i^{-1}-q_i\big)u^{k}q^{kc/2}F_{-\alpha_i+k\delta}
 \otimes E_{-\alpha_i+k\delta}\big)\ts T_{21}
\end{gather}
with respect to $\id\otimes \pi_{V}$ for simple roots $\al_i$ with $i=1,\dots,n$. Using the isomorphism of Section~\ref{subsec:isoDJD}, we can rewrite the internal product
in terms of Drinfeld generators as
\begin{gather*}
\prod_{k> 0}
 \exp_{q_i}\big(\big(q_i^{-1}-q_i\big)\big(uq^{c/2}\big)^{k}q^{-kc}x^{+}_{i,-k}k_i
 \otimes q^{kc}k_i^{-1}x^{-}_{i,k}\big).
\end{gather*}
The calculation breaks into a few cases depending on the type ($B$ and $D$)
and the value of~$i$, but it is quite similar in all cases;
cf.\ \cite[Proposition~5.2]{jlm:ib-c}. We will only give details in the case $i=n$ in type $D$
for the matrix $E^{+}(u)$. Note that $q_n=q$ and so by Proposition~\ref{prop:FFRep},
\begin{gather*}
(\id\otimes \pi_V)\prod_{k> 0}
 \exp_{q}\big(\big(q^{-1}-q\big)(u_{-})^{k}x^{+}_{n,-k}k_n\otimes q^{kc}k_n^{-1}x^{-}_{n,k}\big)\\
\qquad{} =\prod_{k> 0}\exp_{q}\big(q\big(q^{-1}-q\big)(u_{-})^{k}x^{+}_{n,-k}k_n\otimes q^{-(n-1)k}(-e_{n-1,n+1} +e_{(n+1)',(n-1)'})\big).
\end{gather*}
Hence, expanding the $q$-exponent and applying Lemma~\ref{lem:K}, we find that the image of the expression \eqref{epleva} with $i=n$ with respect to the operator $\id\otimes \pi_{V}$ is found by
\begin{gather*}
 1-q\big(q^{-1}-q\big)
\sum\limits_{k> 0}\big((k_{n-1}k_n)^{-1/2}x^{+}_{n,-k}(k_{n-1}k_n)^{1/2}\big)\big(u_{-}q^{-(n-1)}\big)^{k}
 \otimes e_{n-1,n+1}\\
\qquad{} +q\big(q^{-1}-q\big)
\sum\limits_{k\geqslant 0}\big(\big(k_{n-1}^{-1}k_n\big)^{-1/2}x^{+}_{n-1,-k}\big(k_{n-1}^{-1}k_n\big)^{1/2}
\big)\big(u_{-}q^{-(n-1)}\big)^{k}\otimes e_{(n+1)',(n-1)'}.
\end{gather*}
By using the relations $k_i\ts x_{j,k}^{\pm}\tss k_i^{-1}=q_i^{\pm A_{ij}}\ts x_{jk}^{\pm}$,
we can write this expression as
\begin{gather*}
1-\big(q^{-1}-q\big)x^{+}_{n}\big(u_{-}q^{-(n-1)}\big)^{> 0}\otimes e_{n-1,n+1}
 +\big(q^{-1}-q\big)x^{+}_{n}\big(u_{-} \xi q^{n-1}\big)^{> 0}\otimes e_{(n+1)',(n-1)'}\\
\qquad{} =1+e^{+}_{n}(u)\otimes e_{n-1,n+1} -e^{+}_{n}(u)\otimes e_{(n+1)',(n-1)'}.
\end{gather*}
This proves that the $(n-1,n+1)$ entry of $E^{+}(u)$ is $e^{+}_{n}(u)$, while
the $((n+1)',(n-1)')$ entry is $-e^{+}_{n}(u)$, as required.
\end{proof}

In the next proposition we use the series $z^{\pm}(u)$ introduced in~\eqref{zpm}. Their coefficients belong to the center of the algebra $U^{\ext}_q(\wh{\oa}_{N})$; see Proposition~\ref{prop:zu}.
For a nonnegative integer $m$ with $m< n$ we will denote by $z^{\pm\ts[n-m]}(u)$ the respective series
for the subalgebra of $U^{\ext}_q(\wh{\oa}_{N})$, whose generators
are all elements $X^{\pm}_{i,k}$, $h^{\pm}_{j,k}$ and $q^{c/2}$ such that $i,j\geqslant m+1$;
see Definition~\ref{def:eqaa}. We also denote by $\xi^{[n-m]}$ the parameter
$\xi$ for this subalgebra so that
\begin{gather*}
\xi^{[n-m]}=\begin{cases}
 q^{-2n+2m+1}& \text{for type $B$},\\
 q^{-2n+2m+2}& \text{for type $D$}.
 \end{cases}
\end{gather*}

\bpr\label{prop:Hp}
The matrix $H^{\pm}(u)$ is diagonal and has the form
\begin{gather*}
H^{\pm}(u)=\diag\tss\big[h_1^{\pm}(u),\dots,h^{\pm}_{n}(u),h^{\pm}_{n+1}(u),
z^{\pm\tss [1]}(u)\tss h_{n}^{\pm}\big(u\xi^{[1]}\big)^{-1},\dots,
z^{\pm\tss [n]}(u)\tss h_1^{\pm}\big(u\xi^{[n]}\big)^{-1}\big]
\end{gather*}
for type $B$, and
\begin{gather*}
H^{\pm}(u)=\diag\tss\big[h_1^{\pm}(u),\dots,h^{\pm}_{n}(u),
z^{\pm\tss [1]}(u)\tss h_{n}^{\pm}\big(u\xi^{[1]}\big)^{-1},\dots,
z^{\pm\tss [n]}(u)\tss h_1^{\pm}\big(u\xi^{[n]}\big)^{-1}\big]
\end{gather*}
for type $D$.
\epr

\begin{proof}
The starting point is a universal expression for $H^{\pm}(u)$ which is valid for
all three types~$B$,~$D$ and~$C$ (the latter was considered in
\cite[Section~5]{jlm:ib-c}) and is implied by the definition. In particular, for
$H^{+}(u)$ we have:
\begin{gather*}
H^{+}(u) =\exp\left(\sum\limits_{k>0}\sum\limits _{i,j=1}^{n}\frac{\big(q_i^{-1}-q_i\big)
\big(q_j^{-1}-q_j\big)}{q^{-1}-q}\frac{k}{[k]_q}
\tilde{B}_{ij}(q^k)u^k a_{j,-k}\otimes \pi_V(a_{i,k}) \right)\\
\hphantom{H^{+}(u) =}{} \times (\id\ot \pi_V)(T_{21})\prod_{m=0}^{\infty} z^{+}\big(u\xi^{-2m-1}\big)z^{+}\big(u\xi^{-2m-2}\big)^{-1},
\end{gather*}
where the matrix elements $\tilde{B}_{ij}(q)$ are defined in
\eqref{Bijqk} and \eqref{BijqktypeD}. The calculation is then performed in the same way
as for type $C$ with the use of Propositions~\ref{prop:embed}, \ref{prop:FFRep}
and Lemma~\ref{lem:K}; see \cite[Proposition~5.5]{jlm:ib-c}.
\end{proof}

Taking into account Propositions~\ref{prop:Fp}, \ref{prop:Ep} and \ref{prop:Hp}
we arrive at the following result.
\bco\label{cor:isom}
The homomorphism
\begin{gather*}
RD\colon \ U(R)\rightarrow U^{\ext}_q(\wh{\oa}_{N})
\end{gather*}
defined in \eqref{inversemap} is the inverse map to the homomorphism $DR$ defined in Theorem~{\rm \ref{thm:homom}}. Hence the algebra~$U(R)$ is isomorphic to $U^{\ext}_q(\wh{\oa}_{N})$.
\eco

Corollary~\ref{cor:isom} together with the results of Sections~\ref{subsec:eqaa} and \ref{subsec:fsz}
complete the proof of the Main Theorem.

\subsection*{Acknowledgements}

Jing acknowledges the National Natural Science Foundation of
China grant 11531004 and Simons Foundation grant
523868. Liu acknowledges the National Natural Science Foundation
of China grant 11531004, 11701182
and the Guangdong Natural Science Foundation
grant 2019A1515012039.
Liu and Molev acknowledge the support of
the Australian Research Council, grant DP180101825.

\pdfbookmark[1]{References}{ref}
\LastPageEnding

\end{document}